\documentclass{ownamsart}

\input{diagrams}

\newarrow{Dashto}{}{dash}{}{dash}{>}
\newarrow{Eqto}{=}{=}{=}{=}{=}

\usepackage{amssymb}
\usepackage{amscd}
\usepackage{amsmath}
\usepackage{latexsym}
\usepackage{verbatim}
\usepackage[latin1]{inputenc}
\usepackage[dvips]{graphics}

\usepackage{epsfig}

\newtheorem{Thm}{Theorem}[section]
\newtheorem{Prop}[Thm]{Proposition}
\newtheorem{Lem}[Thm]{Lemma}
\newtheorem{Cor}[Thm]{Corollary}

\theoremstyle{definition}
\newtheorem{Def}[Thm]{Definition}
\newtheorem{Ex}[Thm]{Example}

\theoremstyle{remark}
\newtheorem{Rem}[Thm]{Remark}
\newtheorem{Rems}[Thm]{Remarks}
\newtheorem{Obs}[Thm]{Observations}

\newcommand{\Prf}{\noindent\textit{Proof. }}
\newcommand{\PrfOf}[1]{\noindent\textit{Proof of #1.}}
\newcommand{\TopPrfOf}[1]{\noindent\textit{Topological proof of #1.}}
\newcommand{\AnPrfOf}[1]{\noindent\textit{Analytical proof of #1.}}

\newcommand{\bbZ}{\mathbb{Z}}
\newcommand{\bbQ}{\mathbb{Q}}
\newcommand{\bbR}{\mathbb{R}}
\newcommand{\bbC}{\mathbb{C}}
\newcommand{\bbF}{\mathbb{F}}

\newcommand{\Spinc}{\operatorname{Spin^{\hspace*{-.1em}\rm c}}}

\newcommand{\bEG}{\underline{E\Gamma}}

\newcommand{\bE}{\underline{E}}

\newcommand{\betat}[1]{\beta^{({\rm t})}_{#1}}
\newcommand{\betaa}[1]{\beta^{({\rm a})}_{#1}}
\newcommand{\tbetat}[1]{\bar{\beta}^{({\rm t})}_{#1}}
\newcommand{\tbetaa}[1]{\bar{\beta}^{({\rm a})}_{#1}}
\newcommand{\betatev}{\beta_{\even}^{({\rm t})}}
\newcommand{\betaaev}{\beta_{\even}^{({\rm a})}}

\newcommand{\iotaX}[1]{\iota^{\nnspace{\scriptscriptstyle #1}}}

\newcommand{\Ker}{\mathop{\rm Ker}\nolimits}
\newcommand{\Ima}{\mathop{\rm Im}\nolimits}
\newcommand{\diag}{\mathop{\rm Diag}\nolimits}
\newcommand{\id}{\mathop{\rm i\hspace*{-.03em}d}\nolimits}
\newcommand{\incl}{\mathop{\rm incl}\nolimits}
\newcommand{\dom}{\mathop{\rm dom}\nolimits}
\newcommand{\Index}{\mathop{\rm Index}\nolimits}
\newcommand{\Wind}{\mathop{\rm Wind}\nolimits}
\newcommand{\Area}{\mathop{\rm Area}\nolimits}
\newcommand{\Length}{\mathop{\rm Length}\nolimits}
\newcommand{\ev}{\mathop{\rm ev}\nolimits}
\newcommand{\Sign}{\mathop{\rm Sign}\nolimits}
\newcommand{\Lip}{\mathop{\rm Lip}\nolimits}
\newcommand{\Ext}{\mathop{\rm Ext}\nolimits}

\newcommand{\ie}{{\sl i.e.\ }}

\newcommand{\ab}{{\rm a\hspace*{-.05em}b}}

\newcommand{\KK}{K\!K}
\newcommand{\RKK}{R\hspace*{-.03em}K\!K}

\newcommand{\calT}{\mathcal{T}\hspace*{.07em}}

\newcommand{\onto}{\,-\!\!\!\!\twoheadrightarrow}

\newarrow{Eqto}{=}{=}{=}{=}{=}

\newcommand{\itemspace}{\qquad\quad\;\;}

\newcommand{\nnspace}{\hspace*{-.05em}}
\newcommand{\alg}{{\rm a}\nnspace{\rm l}\nnspace{\rm g}}
\newcommand{\topo}{{\rm t}\nnspace{\rm o}\nnspace{\rm p}}

\newcommand{\nnnspace}{\hspace*{-.03em}}
\newcommand{\even}{{\rm e}\nnnspace{\rm v}}
\newcommand{\odd}{{\rm o}\nnnspace{\rm d}\nnspace{\rm d}}

\begin{document}


\title[Unbounded symmetric operators in $K$-homology]{Unbounded symmetric operators in $K$-homology
and the Baum-Connes Conjecture}


\author{Hela Bettaieb}

\address{\'{E}cole d'ing\'{e}nieurs du Canton de Vaud, Route de Cheseaux 1, CH-1401 Yverdon, Switzerland}

\email{hela.bettaieb@eivd.ch}

\author{Michel Matthey}

\address{University of Lausanne, IGAT (Institute for Geometry, Algebra and Topology), BCH,
EPFL, CH-1015 Lausanne, Switzerland}

\email{michel.matthey@ima.unil.ch}

\author{Alain Valette}

\address{Institut de Mathématiques, Université de Neuchâtel, Rue Emile Argand 13,
CH-2000 Neuchâtel, Switzerland}

\email{valette@maths.unine.ch}

\thanks{Research of the first and second named authors partially supported by
Swiss National Fund for Scientific Research grant no. 20-50575.97.}


\subjclass{Primary 19K33, 19K35; Secondary 19K56, 19L64}


\keywords{Unbounded symmetric operators, $K$-homology, $\KK$-theory, Group homology,
Baum-Connes Conjecture}


\begin{abstract}
Using the unbounded picture of analytical $K$-homology, we associate a well-defined
$K$-homology class to an unbounded symmetric operator satisfying certain mild technical
conditions. We also establish an ``addition formula'' for the Dirac operator on the circle and
for the Dolbeault operator on closed surfaces. Two proofs are provided, one using topology
and the other one, surprisingly involved, sticking to analysis, on the basis of the previous
result. As a second application, we construct, in a purely analytical language, various
homomorphisms linking the homology of a group in low degree, the $K$-homology of its
classifying space and the analytic $K$-theory of its $C^{*}$-algebra, in close
connection with the Baum-Connes assembly map. For groups classified by a $2$-complex,
this allows to reformulate the Baum-Connes Conjecture.
\end{abstract}


\maketitle


\part{Introduction}

\label{p-I}


\section{Statement of the main results}

\label{s-intro}


%
The non-commutative geometry approach to $K$-homology rests on the
concept of \emph{unbounded Fredholm module}, due to Connes
(\cite{NCDG}, Chap.~I, Section~6). Subsequently, this object was
renamed \emph{$K$-cycle} (\cite{NCG}, Def.~11 in Section~IV.2.$\gamma$)
and then, quite conveniently, \emph{spectral triple}
(see \cite{ConMos}) to emphasize the connection with spectral
geometry. Recall that, if $\mathcal{A}$ is an involutive algebra
represented on the Hilbert space $\mathcal{H}$\,, and $D$ is a
self-adjoint operator on $\mathcal{H}$ with compact resolvent, the
triple $(\mathcal{A},\mathcal{H},D)$ is \emph{spectral} if $D$ almost commutes
with any $a\in \mathcal{A}$\,, \ie if $[D,a]$ is bounded for every $a\in
\mathcal{A}$\,.

\smallskip

Given a separable $C^{*}$-algebra $A$\,, our goal is to define certain classes in the $K$-homology
group $K^{*}(A):=\KK_{*}(A,\bbC)$\,, using unbounded symmetric operators (that are definitely
not assumed to be self-adjoint). Not surprisingly, this will force the operators considered
to fulfill some technical conditions. The following definition lists the properties
we need (more details are provided in Section~\ref{s-KK} below, in particular concerning
deficiency indices and invertibility of $T^{*}T+1$).

\begin{Def}
\label{def-triple}
We call a triple $(\mathcal{H},\pi,T)$ a \emph{symmetric unbounded
Fredholm module of degree $i$}
over the separable $C^{*}$-algebra $A$ if it consists of the following data\,:
\begin{itemize}
    \item [(a)] an integer $i\in\{0,1\}$\,;
    \item [(b)] a Hilbert space $\mathcal{H}$\,;
    \item [(c)] a $*$-representation $\pi\colon A\longrightarrow B(\mathcal{H})$ of $A$\,;
    \item [(d)] a densely defined closed symmetric operator $T$ on $\mathcal{H}$ with domain
    $\dom(T)$\,.
\end{itemize}
These data are required to fulfill the following conditions\,:
\begin{itemize}
    \item [{\bf (i)}] the deficiency indices of $T$ coincide and are finite;
    \item [{\bf (ii)}] the operator $\pi(a)(T^{*}T+1)^{-1}$ is compact for every $a\in A$\,;
    \item [{\bf (iii)}] the operators $[\pi(b),T]$ and $[\pi(b),T^{*}]$ are densely defined
    and $[\pi(b),T^{*}]$ is bounded for every $b$ in some dense \emph{subspace} $\mathcal{B}$
    of $A$\,.
\end{itemize}
(The subspace $\mathcal{B}$ is not required to be a subalgebra.) If $i=0$\,, we moreover require $\mathcal{H}$
to be $\bbZ/2$-graded, $\pi$ to preserve the grading, and $T$ to reverse it.
\end{Def}

Of course, this would define an unbounded Fredholm module $[\mathcal{H},\pi,T]\in \KK_{i}(A,\bbC)$
in the sense of Connes, if $T$ would moreover be self-adjoint. In fact, our two main
results read as follows.

\begin{Thm}
\label{thm-main}
Let $(\mathcal{H},\pi,T)$ be a symmetric unbounded Fredholm module of degree $i$
over a separable $C^{*}$-algebra $A$\,, as defined above. Then, there exists, in the
unbounded picture of analytical $K$-homology, a well-defined $K$-homology class
$$
[\mathcal{H},\pi,T]\in K^{i}(A)\,,
$$
that is canonical, and coincides with the usual class in case $T$ is self-adjoint.
More precisely, given an arbitrary self-adjoint extension $\widetilde{T}$ of $T$
-- and at least one such extension exists --, one has
$$
[\mathcal{H},\pi,T]=[\mathcal{H},\pi,\widetilde{T}]\;\in\;K^{i}(A)\,,
$$
independently of the choice of $\widetilde{T}$\,.
\end{Thm}

As a consequence of this first result (more precisely of a generalization of it), we will
derive the following theorem.

\begin{Thm}
\label{thm-main2}
Let $(\mathcal{H},\pi,T_{1})$ and $(\mathcal{H},\pi,T_{2})$ be two unbounded
Fredholm modules of degree $i$ for a separable $C^{*}$-algebra $A$ (in the usual
sense of Connes), with the same Hilbert space $\mathcal{H}$ and the same $*$-representation
$\pi$\,. Suppose that $T_{1}$ and $T_{2}$ admit, as a common restriction, an operator $T$
satisfying the two conditions
\begin{itemize}
    \item [(a)] $T$ is densely defined;
    \item [(b)] $[\pi(b),T^{**}]$ is densely defined and bounded for every $b$ in some dense $*$-closed
    subspace $\mathcal{B}$ of $A$\,.
\end{itemize}
Then, one has the following equality of $K$-homology classes\,:
$$
[\mathcal{H},\pi,T_{1}]=[\mathcal{H},\pi,T_{2}]\;\in\;K^{i}(A)\,.
$$
Furthermore, since $T_{1}$ enters in an unbounded Fredholm module,
there exists a $*$-closed dense subspace $\mathcal{B}'$ of $A$ such that $[\pi(b'),T_{1}]$
is densely defined and bounded for every $b'\in\mathcal{B'}$\,, and condition {\rm (b)} above
can be replaced by the next one while keeping the same conclusion\,:
\begin{itemize}
\item [(b$^{\prime}$)] $[\pi(b'),T^{**}]$ is densely defined for every $b'$ in $\mathcal{B'}$\,.
\end{itemize}
\end{Thm}

In both theorems, the separability assumption is needed to apply Baaj-Julg's results~\cite{BJ} (see the
proof of Proposition~2.3 therein). Recall, for later applications, that for $X$ a compact Hausdorff space,
the $C^{*}$-algebra $C(X)$ is separable if and only if $X$ is metrizable (or equivalently, second-countable).
For example, as is well-known, a CW-complex is metrizable if and only if it is locally finite~\cite[Prop.~1.5.17]{FriPic}.
At this point, let us mention that \emph{throughout the paper, we assume that all spaces and maps between
them are pointed}.

\medskip

Here is a description of the content of the paper.

\smallskip

One of our goals is to apply our results to establish an ``addition formula'' concerning certain
differential operators on closed manifolds of dimension one and two. More precisely, we would
like to study the behaviour under connected sum of the $K$-homology class given by the Dirac operator
in dimension one and by the Dolbeault operator in dimension two. In fact, in dimension one, the situation
is well-behaved for the usual connected sum, but in dimension two, this leads to the introduction
of a variation of the connected sum. For both considered dimensions, we will present two proofs of each
``addition formula'', one using standard and well-established tools from algebraic topology (in particular,
``topological index theory''), and the other one in a purely analytical framework on the basis of Theorem~\ref{thm-main2}.
One of the interest of this latter approach is that the analytical proof is astonishingly involved. The
context will be explained in detail in Section~\ref{s-appl-mfld}, and the proofs are presented later,
in Section~\ref{s-mfld-top-proof} for the topological proof, in Section~\ref{s-mfld-proofs-1} for the
analytical proof in dimension one, and in Section~\ref{s-mfld-proofs-2} for the analytical proof in
dimension two.

\smallskip

In Section~\ref{s-appl-BC}, we explain in detail the framework of our application of these results in connection
with the Baum-Connes Conjecture, which aims at computing the $K$-group $K_{*}(C^{*}_{r}\Gamma)$ for a group $\Gamma$
(see there for the notation and definitions). The $K$-homology group $K_{j}(B\Gamma)$ and the homology group
$H_{j}(\Gamma;\bbZ)$ will be involved for $j=1,2$\,. More precisely, for both values of $j$\,,
we will construct two maps
$$
\betat{j}\colon H_{j}(\Gamma;\bbZ)\longrightarrow K_{j}(B\Gamma)\qquad\mbox{and}\qquad\betaa{j}\colon H_{j}(\Gamma;
\bbZ)\longrightarrow K_{j}(C^{*}_{r}\Gamma)\,.
$$
Our main concern will be to define these maps in a purely analytical language, \ie using the
unbounded picture
of analytical $K$-homology, and, as a major difficulty, to prove that they are well-defined group homomorphisms, while sticking to
this analytical language. It turns out that the proof of this property will precisely amount to the ``addition formulae''
for suitable differential operators as in Section~\ref{s-appl-mfld}, hence the close connection with Theorem~\ref{thm-main2}.

\smallskip

In Section~\ref{s-KK}, we state a general theorem, namely Theorem~\ref{thm-KK}, that allows to associate to an unbounded
symmetric Fredholm module $(\mathcal{H},\pi,T)$ a ``usual'' $\KK$-theory class in some group $\KK_{i}(A,C(\mathcal{U}))$\,,
where $\mathcal{U}$ is a suitable non-empty compact Hausdorff space depending on $T$\,. Evaluation at an arbitrary point
$u$ of $\mathcal{U}$ will provide the $K$-homology class $[\mathcal{H},\pi,T]\in K^{*}(A)$ we are looking for. The punch-line
is that this will \emph{not} depend on the choice of $u$ (the point being path-connectedness of $\mathcal{U}$). As a consequence,
Theorem~\ref{thm-main} follows from this.

\smallskip

The proof of Theorem~\ref{thm-KK}, and hence of Theorem~\ref{thm-main}, is presented in Section~\ref{s-proofs}.

\smallskip

In Section~\ref{s-general}, we address a generalization of Theorem~\ref{thm-main}, where we reduce the
assumptions on the triple $(\mathcal{H},\pi,T)$ to the strict minimum (according to our proof); as the
main relaxation of assumptions, finiteness of deficiency indices will be dropped. Using this generalization,
we then establish Theorem~\ref{thm-main2}. After this section, we move, for the rest of the paper, to the
applications, namely on the ``addition formulae'' and around the Baum-Connes Conjecture.

\smallskip

As we have said, we will present the proofs of the ``addition formulae'' for the Dirac and the Dolbeault
operators in Section~\ref{s-mfld-top-proof} (topological in both dimensions), in Section~\ref{s-mfld-proofs-1}
(analytical in dimension one) and in Section~\ref{s-mfld-proofs-2} (analytical in dimension two).

\smallskip

We treat the case $j=1$ of our application to the Baum-Connes Conjecture in Section~\ref{s-j1}. In this case,
$H_{1}(\Gamma;\bbZ)$ is $\Gamma^{\ab}$\,, the abelianization of $\Gamma$\,. We will see that $\betaa{1}$ is exactly
the map $\Gamma^{\ab}\longrightarrow K_{1}(C^{*}_r\Gamma)$ induced by the canonical inclusion of $\Gamma$ in the
group of invertibles of $C^{*}_{r}\Gamma$\,. It was proved by Elliott and Natsume~\cite{ElN} (and reproved in~\cite{BV})
that $\betaa{1}$ is rationally injective.

\smallskip

In Section~\ref{s-j2}, for $j=2$\,, Zimmermann's description of $H_{2}(\Gamma;\bbZ)$ in~\cite{Zim}
allows us to define $\betat{2}$ and $\betaa{2}$\,. We were not able to prove rational injectivity
of $\betaa{2}$\,.

\smallskip

In Section~\ref{s-two-dim}, we draw consequences of our constructions for groups which admit a $2$-dimensional
classifying space; we call these groups \emph{$2$-dimensional}. We use our maps $\betaa{j}$ to propose, for these
groups, an equivalent formulation of the Baum-Connes Conjecture with the left hand side replaced by integral group
homology.

\smallskip

We point out that~\cite{MM2,MM1} contain closely related results; the relation of
the maps $\betaa{j}$ with algebraic $K$-theory (and Steinberg symbols for $j=2$) is
studied in~\cite{MOO}.


\tableofcontents


\section{Description of the application to analysis on manifolds}

\label{s-appl-mfld}


%
We explain here two applications of Theorem~\ref{thm-main2} in the context of
differential operators on manifolds. One of the applications is for the circle,
\ie in dimension one, and the other is in dimension two, more precisely for
Riemann surfaces. Explicitly, we will state two ``addition formulae'' for
suitable differential operators. The topological proof is presented in
Section~\ref{s-mfld-top-proof}, and the analytical proof in Section~\ref{s-mfld-proofs-1}
for the one-dimensional case, and in Section~\ref{s-mfld-proofs-2} for
the two-dimensional case.

\smallskip

We first recall that for a $\sigma$-compact Hausdorff topological space $X$\,, for instance
a CW-complex, one has a canonical and natural isomorphism
$$
K_{*}(X)\cong\RKK_{*}(X,\bbC)\,,
$$
where $K_{*}$ is $K$-homology with compact supports, and $\RKK_{*}$ is Kasparov's $\KK$-theory
with compact supports, that we will see in the unbounded picture of
$K$-homology (more on
this is contained in Sections \ref{s-KK}, \ref{s-mfld-proofs-1} and \ref{s-mfld-proofs-2}). If $X$
is compact Hausdorff, one further has $\RKK_{*}(X,\bbC)=\KK_{*}(C(X),\bbC)$\,.

\smallskip

Now, we start with the one-dimensional situation. Consider the Dirac operator on the circle $S^{1}$
and the corresponding $K$-homology class, namely
$$
D:=\frac{1}{i}\!\cdot\!\frac{d}{d\theta}\qquad\mbox{and}\qquad[D]\in K_{1}(S^{1})
\cong\KK_{1}(C(S^{1}),\bbC)\,;
$$
details are provided in Section~\ref{s-mfld-proofs-1}. Now, the ``addition formula'' reads as
follows.

\begin{Thm}
\label{thm-add-formula-1}
Let $X$ be a pointed CW-complex, and let $f_{1},f_{2}\colon S^{1}\longrightarrow X$ be
two pointed continuous maps, that are constant in a small neighbourhood of the
base-point of $S^{1}$\,. Consider the connected sum of these two maps (along
a closed interval sitting inside the given neighbourhood for both copies of $S^{1}$)
$$
f_{1}\# f_{2}\colon S^{1}\# S^{1}\longrightarrow X\,,
$$
and identify the closed oriented manifold $S^{1}\# S^{1}$ with $S^{1}$ as usual. Then, in
$K$-homology, one has
$$
(f_{1}\# f_{2})_{*}[D]=(f_{1})_{*}[D]+(f_{2})_{*}[D]\;\in\;K_{1}(X)\cong\RKK_{1}(X,\bbC)\,.
$$
\end{Thm}

We pass to the two-dimensional setting. Let $\Sigma_{g}$ be a closed oriented surface of genus
$g\geq 0$ (in particular, without boundary). We fix an auxiliary K\"{a}hler structure on
$\Sigma_{g}$\,, \ie we view $\Sigma_{g}$ as a complex curve equipped with a suitably compatible
Riemannian metric. Consider the Dolbeault operator and its $K$-homology class
$$
\bar{\partial}_{\Sigma_{g}}:=\bar{\partial}\oplus\bar{\partial}^{*}\qquad\mbox{and}\qquad[\bar
{\partial}_{g}]:=[\bar{\partial}_{\Sigma_{g}}]\in K_{0}(\Sigma_{g})\cong\KK_{0}(C(\Sigma_{g}),\bbC)\,;
$$
again, we will be more explicit in Section~\ref{s-mfld-proofs-2}. As we will explain is that section,
the connected sum for surfaces does \emph{not} satisfy the ``same addition formula'' as in Theorem~\ref{thm-add-formula-1}.
We will explain that the exact source of the problem is the non-additivity of the Euler characteristic
under the connected sum. We therefore introduce a modified version of the usual connected sum.

\smallskip

Thus, let $\Sigma_{g_{1}}$ and $\Sigma_{g_{2}}$ be surfaces of genus $g_{1}$ and $g_{2}$ respectively.
By cutting out a handle in $\Sigma_{g_{1}}$ (resp.\ in $\Sigma_{g_{2}}$), see Figure 1, and gluing along
the boundary circles in an orientation preserving way, we get a closed oriented surface $\Sigma_{g_{1}}
\natural\Sigma_{g_{2}}$ of genus $g_{1}+g_{2}-1$\,, as in Figure 2.

\epsfig{file=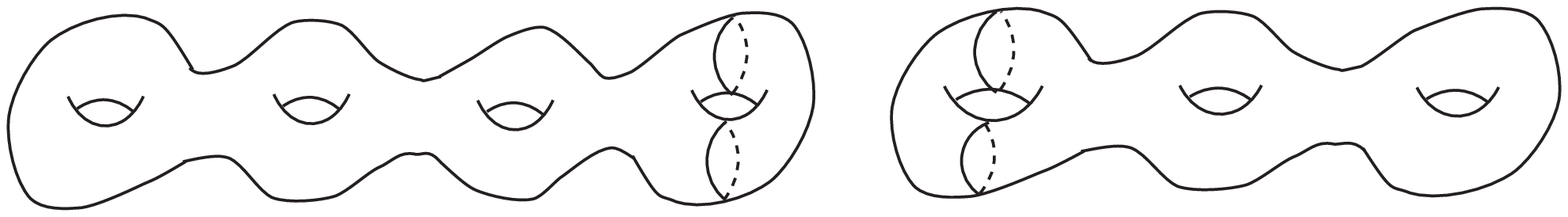,width=33em}

\vspace*{.2em}

\centerline{Figure 1}

\vspace*{-.3em}

\epsfig{file=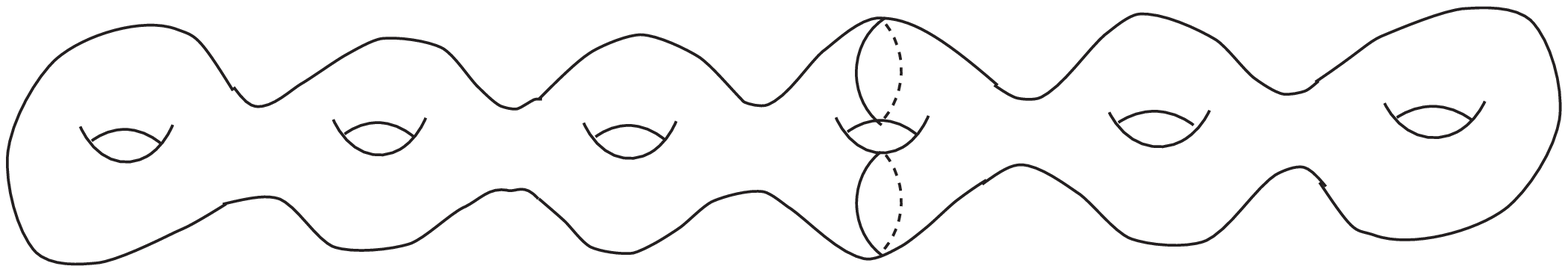,width=33em}

\vspace*{-2.1em}

\centerline{Figure 2}

\vspace*{1em}

We assume that the base-point of $\Sigma_{g_{1}}$ is identified with the base-point of $\Sigma_{g_{2}}$ in this
operation (in particular, both base-points sit on two corresponding circles among the four boundary circles). We
single out that the Euler characteristic is additive for this modified connected sum, \ie
$$
\chi(\Sigma_{g_{1}}\natural\Sigma_{g_{2}})=\chi(\Sigma_{g_{1}})+\chi(\Sigma_{g_{2}})\,.
$$
In this situation, the ``addition formula'' reads as follows.

\begin{Thm}
\label{thm-add-formula-2}
Let $X$ be a pointed CW-complex. For $i=1,2$\,, let $f_{i}\colon\Sigma_{g_{i}}\longrightarrow X$
be a pointed continuous map, that is constant in a small neighbourhood of a handle of $\Sigma_{g_{i}}$
($g_{i}\geq 0$). Consider the modified connected sum (along the two given handles)
$$
f_{1}\natural f_{2}\colon\Sigma_{g_{1}}\natural\Sigma_{g_{2}}\longrightarrow X
$$
of these two maps, and identify the closed oriented manifold $\Sigma_{g_{1}}\natural\Sigma_{g_{2}}$ with
$\Sigma_{g_{1}+g_{2}-1}$ in the usual way. Then, in $K$-homology, one has
$$
(f_{1}\natural f_{2})_{*}[\bar{\partial}_{g_{1}+g_{2}-1}]=(f_{1})_{*}[\bar{\partial}_{g_{2}}]+(f_{2})_{*}[\bar{
\partial}_{g_{2}}]\;\in\;K_{0}(X)\cong\RKK_{0}(X,\bbC)\,.
$$
\end{Thm}

For the analytical proofs of Theorems \ref{thm-add-formula-1} and \ref{thm-add-formula-2}, we have
to impose suitable boundary conditions on the glued parts of the considered manifolds; as a consequence,
we must deal with symmetric non-self-adjoint operators. Theorems \ref{thm-main} and \ref{thm-main2}
then ensure the well definiteness of the corresponding ``glued'' $K$-homology classes.

\smallskip

Another occurrence of symmetric non-self-adjoint operators arises in analytic $K$-homology, in
the discussion of excision in that framework, see~\cite[Section~10.8]{HigRoe}.


\section{Description of the application to the Baum-Connes Conjecture}

\label{s-appl-BC}


%
We describe here our second application of Theorem~\ref{thm-main}, namely, in the
framework of the celebrated Baum-Connes Conjecture, that we also introduce
with some explanations.

\smallskip

Let $\Gamma$ be a countable discrete group. The \emph{Baum-Connes Conjecture} for $\Gamma$ is the
statement that the Baum-Connes assembly map, or analytical index map,
$$
\mu^{\Gamma}_{i}\colon K^{\Gamma}_{i}(\bEG)\longrightarrow K_{i}(C^{*}_{r}\Gamma)\qquad(i=0,1)
$$
is an isomorphism. Here, $K^{\Gamma}_{*}(\bEG)$ is the $\Gamma$-equivariant $K$-homology with $\Gamma$-compact
supports of $\bEG$\,, the classifying space for proper $\Gamma$-actions, and $K_{*}(C^{*}_{r}\Gamma)$ is the
analytical $K$-theory of $C^{*}_{r}\Gamma$\,, the reduced $C^{*}$-algebra of $\Gamma$\,. For precise definitions
of the objects involved, various examples and the relevance of the conjecture to questions in topology, geometry,
algebra and analysis, we refer to~\cite{BCH,Val,Sch,Val99}; see also~\cite{Julg,GSk} for excellent
surveys of progresses on the conjecture up to 1999. Recall also that both $K^{\Gamma}_{*}(\bEG)$ and $K_{*}
(C^{*}_{r}\Gamma)$ are $2$-periodic by virtue of Bott periodicity. For this reason, we will stick to the groups
$K_{0}$ and $K_{1}$\,.

\smallskip

Denote by $F\Gamma$ the vector space of $\bbC$-valued functions on $\Gamma$\,, with finite support contained in the
set of finite-order elements of $\Gamma$\,. The space $F\Gamma$ becomes a $\Gamma$-module by letting $\Gamma$ act by
conjugation; $H_{*}(\Gamma;F\Gamma)$ denotes the corresponding homology group. Baum and Connes defined in~\cite{BC}
a Chern character
$$
ch^{\Gamma}_{*}\colon K^{\Gamma}_{i}(\bEG)\longrightarrow\bigoplus_{n=0}^{\infty}H_{2n+i}(\Gamma;F\Gamma)\qquad(i=0,1)
$$
that becomes an isomorphism after tensoring with $\bbC$ (see also~\cite{MM1}).

\smallskip

Now, let $B\Gamma$ be the classifying space of $\Gamma$\,, and let $K_{*} (B\Gamma)$ denote its $K$-homology with compact
supports, which is also $2$-periodic by Bott periodicity. There is a canonical map $\varphi_{*}^{\Gamma}\colon K_{*}
(B\Gamma)\longrightarrow K^{\Gamma}_{*}(\bEG)$\,, which is an isomorphism for $\Gamma$ torsion-free. Indeed, denote
by $E\Gamma$ the universal cover of $B\Gamma$\,. This map $\varphi_{*}^{\Gamma}$ is the composition of the canonical
isomorphism $K_{*}(B\Gamma)\cong K^{\Gamma}_{*}(E\Gamma)$ with the forgetful map $K^{\Gamma}_{*}(E\Gamma)\longrightarrow
K_{*}^{\Gamma}(\bEG)$ obtained by noticing that any free and proper $\Gamma$-action is proper. Of course, for $\Gamma$
torsion-free, the spaces $E\Gamma$ and $\bEG$ coincide (up to $\Gamma$-equivariant homotopy). Together with the Chern
character $ch_{*}$ in $K$-homology, these maps fit into the commutative diagram (see~\cite{BCH} and~\cite{MM1})
\begin{displaymath}
 \begin{diagram}[height=2em,width=3em]
   & & K_{i}(B\Gamma) & \rTo^{\varphi_{i}^{\Gamma}} & K_{i}^{\Gamma}(\bEG) &
   \rTo^{\!\!\!\!\!\!\!\!\!\!\mu_{i}^{\Gamma}} & K_{i}(C^{*}_{r}\Gamma) \\
   & & \dTo^{ch_{*}} & & \dTo_{ch_{*}^{\Gamma}} & & \\
  \bigoplus_{n=0}^{\infty}H_{2n+i}(\Gamma;\bbZ) & \rTo & \smash
  {\bigoplus_{n=0}^{\infty}}H_{2n+i}(\Gamma;\bbQ) & \rInto & \smash
  {\bigoplus_{n=0}^{\infty}}H_{2n+i}(\Gamma;F\Gamma) & & \\
 \end{diagram}
\end{displaymath}
for $i=0$ and $1$\,. (Throughout the paper, we identify the integral (resp.\ rational) homology of $\Gamma$ with that
of $B\Gamma$\,.) This shows in particular that $\varphi_{*}^{\Gamma}$ is rationally injective. We let $\nu^{\Gamma}_
{*}:=\mu^{\Gamma}_{*}\circ\varphi_{*}^{\Gamma}\colon K_{*}(B\Gamma)\longrightarrow K_{*}(C^{*}_{r}\Gamma)$ be the
Novikov assembly map. The reason for this terminology is that rational injectivity of $\nu^{\Gamma}_{*}$ implies
the Novikov Conjecture on higher signatures for the group $\Gamma$\,.

\smallskip

At the very beginning, this paper started out from a desire to exploit the bottom line of this diagram, in order
to better understand the top line. Since, in favorable cases, geometry and topology provide explicit models for
$B\Gamma$\,, from which group homology $H_{*}(\Gamma;\bbZ)$ can be computed, or at least well-understood, it seems
interesting to try to construct directly, out of integral homology classes, elements in $K_{*}(B\Gamma)$ and
$K_{*}(C^{*}_{r}\Gamma)$\,. In other words, we are looking for maps
$$
\betat{j}\colon H_{j}(\Gamma;\bbZ)\longrightarrow K_{i}(B\Gamma)\qquad\mbox{and}\qquad\betaa{j}\colon H_{j}(\Gamma;\bbZ)
\longrightarrow K_{i}(C^{*}_{r}\Gamma)\,,
$$
where $i\equiv j\pmod{2}$\,, such that the diagram
\begin{displaymath}
 \begin{diagram}[height=2em,width=3em]
   K_{i}(B\Gamma) & & \rTo^{\nu_{i}^{\Gamma}} & & K_{i}(C^{*}_{r}\Gamma) \\
    & \luTo_{\betat{j}} & & \ruTo_{\betaa{j}} & \\
    & & H_{j}(\Gamma;\bbZ) & & \\
 \end{diagram}
\end{displaymath}
commutes. To ensure non-triviality, $\betat{j}$ should be rationally a right-inverse of the Chern character in degree
$j$\,, \ie
$$
(ch_{j}\otimes\id_{\bbQ})\circ(\betat{j}\otimes\id_{\bbQ})=\id_{H_{j}(\Gamma;\bbQ)}\,.
$$
Moreover, we do not want to define $\betaa{j}$ merely as $\nu^{\Gamma}_{j}\circ\betat{j}$\,, but look instead
for a direct and explicit construction. Indeed, it would follow from the Baum-Connes Conjecture that $\betaa{j}$
is rationally injective; one may then try to prove this directly.

\smallskip

To illustrate this program, let us consider the easy case where $j=0$\,. Of course $H_{0}(\Gamma;\bbZ)\cong\bbZ$\,,
and we define
$$
\betat{0}\colon\bbZ\longrightarrow K_{0}(B\Gamma)\,,\quad n\longmapsto n\cdot\iotaX{B\Gamma}_{*}[1]\,,
$$
where $\iotaX{B\Gamma}\colon pt\longrightarrow B\Gamma$ is the inclusion of the base-point, and the class $[1]$
is a prescribed generator of $K_{0}(pt)\cong\bbZ$\,. It is obvious that $\betat{0}$ is a right-inverse of the
map $ch_{0}^{\bbZ}\colon K_{0}(B\Gamma)\longrightarrow H_{0}(\Gamma;\bbZ)$\,, \ie the integral Chern character
in degree zero (compare with Lemma~\ref{lem4} below, and with~\cite{MM2}). On the other hand, we define
$$
\betaa{0}\colon\bbZ\longrightarrow K_{0}(C^{*}_{r}\Gamma)\,,\quad n\longmapsto n\cdot[1]=\Sign(n)\cdot\big[\diag
(\underbrace{1,\ldots,1}_{|n|\;\mbox{\footnotesize terms}},0,0,\ldots)\big]\,,
$$
where $[1]$ denotes, this time, the $K$-theory class of the unit in $C^{*}_{r}\Gamma$\,. It is an easy but instructive
exercise (see e.g. \cite[Ex.~2.11 in Part~2]{Val} or \cite[Ex.~6.1.5]{Val99}) to check that $\nu^{\Gamma}_{0}\iotaX
{B\Gamma}_{*}[1]=[1]$\,. Moreover, the canonical trace $\tau$ on $C^{*}_{r}\Gamma$ induces a map $\tau_{*}\colon
K_{0}(C^{*}_{r}\Gamma)\longrightarrow\bbR$ such that $\tau_{*}[1]=1$\,. This shows for free that $\betaa{0}$ is
injective.

\smallskip

In this paper, we implement the program sketched above in the cases $j=1$ and $j=2$\,, exploiting especially simple
descriptions of $H_{j}(\Gamma;\bbZ)$ available in this range. The main feature is that the construction of $\betat{j}$
is performed in the analytical and unbounded description -- that is, ``à
la Kasparov and Connes'' -- of $K$-homology.
In particular, the proof of the fact that $\betat{j}$ is a group homomorphism is instructively subtle in the analytical
framework and is presented in full details. Roughly speaking, this leads to the \emph{construction of well-defined
$K$-homology classes out of densely defined unbounded symmetric operators on Hilbert spaces, that are not necessarily
self-adjoint}, precisely the subject of Theorem~\ref{thm-main}, hence the connection.


\part{Symmetric unbounded Fredholm modules}

\label{p-II}


\section{Construction of the class $[\mathcal{H},\pi,T]$}

\label{s-KK}


%
In this section, we provide some general information on unbounded symmetric operators
and we construct the promised $K$-homology class $[\mathcal{H},\pi,T]$\,. We also state
a general result, Theorem~\ref{thm-KK}, of which Theorem~\ref{thm-main} is a direct corollary.
The proof will be presented in Section~\ref{s-proofs}.

\smallskip

Recall that for a densely defined closed operator $T$\,, the operator $T^{*}T+1$ is densely defined,
self-adjoint, injective on its domain and surjective, and its inverse satisfies $(T^{*}T+1)^{-1}\in
B(\mathcal{H})$ (see for example \cite[Prop.~A.8.4, p.~511]{TayI}). Therefore, condition~{\bf (ii)}
of Definition~\ref{def-triple} makes sense.

\smallskip

It is well-known that a densely defined closed symmetric operator $T$ is self-adjoint
if and only if $\Ker(T^{*}{-i})$ and $\Ker(T^{*}{+i})$ are trivial. In general, there can
exist none, just one (in case $T$ is already self-adjoint), or uncountably many self-adjoint
extensions of $T$\,. In fact, they are canonically parameterized by the space
$$
\mathcal{U}:=\big\{u\colon\Ker(T^{*}{-i})\longrightarrow\Ker(T^{*}{+i})\,
\big|\,u\;\,\mbox{is a unitary isomorphism}\big\}\,,
$$
equipped with the norm-topology inherited from $B\big(\Ker(T^{*}{-i}),\Ker(T^{*}{+i})\big)$\,.
This means in particular that $T$ possesses self-adjoint extensions if and only
if the \emph{deficiency spaces} $\Ker(T^{*}{-i})$ and $\Ker(T^{*}{+i})$ of $T$ have
the same (possibly infinite) dimension; here dimension, is meant in the sense
of the cardinal of a Hilbert basis. If the \emph{deficiency indices} $\dim
\big(\Ker(T^{*}{-i})\big)$ and $\dim\big(\Ker(T^{*}{+i})\big)$ of $T$ are finite and
equal, say equal to $n$\,, then $\mathcal{U}$ is a principal homogeneous space over
the unitary group in dimension $n$\,, hence $\mathcal{U}$ is homeomorphic to $U(n)$
(non-canonically for $n>0$); in particular, it is Hausdorff and compact. Explicitly,
in the general case, the correspondence is given as follows (provided that $\mathcal{U}$
is non-empty)\,:
$$
\mathcal{U}\ni u\qquad\longleftrightarrow\qquad T_{u}\colon\dom(T_{u})
\longrightarrow\mathcal{H}\,,
$$
where $T_{u}$ is the unbounded operator with domain
$$
\dom(T_{u}):=\big\{\xi+\eta+u(\eta)\,\big|\,\xi\in\dom(T)\;\,\mbox{and}
\;\,\eta\in\Ker(T^{*}{-i})\big\}
$$
and given by the (well-defined) formula
$$
T_{u}\big(\xi+\eta+u(\eta)\big):=T(\xi)+i\eta-iu(\eta)\,.
$$
Finally, letting ``\,$\oplus$\,'' stand for the algebraic direct sum (not necessarily orthogonal),
we point out that $\dom(T^{*})=\dom(T)\oplus\Ker(T^{*}{-i})\oplus\Ker(T^{*}{+i})$ and that every
$T_{u}$ is a restriction of $T^{*}$\,. For the details and proofs, we refer, for instance, to
Reed-Simon \cite[Section~X.1, pp.~135--143]{RSII}.

\smallskip

For the sequel, we suppose that the deficiency indices of $T$ coincide and are
equal to $n<\infty$ (in particular, $\mathcal{U}$ is compact Hausdorff and $C(\mathcal{U})$
is a unital separable $C^{*}$-algebra). We consider the unbounded operator $T\hat{\otimes}1$ on the
Banach space $\mathcal{H}\hat{\otimes}C(\mathcal{U})$\,, viewed as a Hilbert $C^{*}$-module
over $C(\mathcal{U})$ in the obvious way or, in other words, as a (constant) continuous
field of Hilbert spaces over $\mathcal{U}$\,. The operator $T\hat{\otimes}1$ has, as
domain, the image of the algebraic tensor product $\dom(D)\otimes C(\mathcal{U})$ in
$\mathcal{H}\hat{\otimes}C(\mathcal{U})$\,. It admits a canonical extension $\calT$\,,
which is the unbounded operator equal to $T_{u}$ in the fiber over each $u\in
\mathcal{U}$\,. Let us provide an explicit description of $\calT$\,. First,
we use the canonical isomorphism of Hilbert $C(\mathcal{U})$-modules $\mathcal{H}
\hat{\otimes}C(\mathcal{U})\cong C(\mathcal{U},\mathcal{H})$ (see~\cite[p.~27]{Lance})
to identify both Hilbert $C^{*}$-modules. As usual, $C(\mathcal{U},\mathcal{H})$ is
endowed with the $C(\mathcal{U})$-valued scalar product $\left<f\,|\,g\right>:=
\big(u\mapsto\left<f(u)\,|\,g(u)\right>_{\mathcal{H}}\big)$ for $f,\,g\in C(\mathcal
{U},\mathcal{H})$\,, therefore with the topology of uniform convergence. Then,
$T\hat{\otimes}1$ has, as domain, the dense subspace $\big\{f\in C(\mathcal{U},
\mathcal{H})\,\big|\,\Ima(f)\subseteq\dom(T)\big\}$ of $C(\mathcal{U},\mathcal{H})$\,,
and $\calT$ is defined as the $C(\mathcal{U})$-linear operator with domain
$$
\dom(\calT):=\big\{f\in C(\mathcal{U},\mathcal{H})\,\big|\,f(u)\in\dom
(T_{u})\,,\;\forall u\in\mathcal{U}\big\}
$$
and given by
$$
\calT\colon\dom(\calT)\longrightarrow C(\mathcal{U},\mathcal{H})\,,\quad
f\longmapsto\Big(\calT\!f\colon u\mapsto T_{u}\big(f(u)\big)\Big)\,.
$$
Observe that $\calT$ is a restriction of the closed operator $T^{*}\hat{\otimes}1$\,,
consequently, it is well-defined, \ie $\calT\!f$ is continuous for every $f\in
C(\mathcal{U},\mathcal{H})$\,. Note also that $T^{*}\hat{\otimes}1$ is symmetric
if and only if $T$ is self-adjoint, in which case $\mathcal{U}$ is a point and
$\calT=T$\,.

\smallskip

By definition, a densely defined closed operator $\calT$ on a Hilbert $C^{*}$-module
is called \emph{regular} if $\calT^{*}$ is densely defined and $\calT^{*}\calT+1$ has
dense range (of course, if $\calT$ is symmetric, only the latter property
is significant). For $\calT$ regular, the operator $\calT^{*}$ is regular, and
$\calT^{*}\calT$ is densely defined, self-adjoint and regular (see \cite[Lem.~9.1, Cor.~9.6 and Prop.~9.9]{Lance}).

\begin{Thm}
\label{thm-KK}
Let $(\mathcal{H},\pi,T)$ be a symmetric unbounded Fredholm module of degree~$i$
over a separable $C^{*}$-algebra $A$\,. Let $\mathcal{U}$ and $\calT$ be as
constructed above. Then, the operator $\calT$ has a dense domain, is self-adjoint,
and is regular, in the sense that $\calT^{2}+1$ has dense range. Moreover, the triple
$\big(\mathcal{H}\hat{\otimes}C(\mathcal{U}),\pi\hat{\otimes}1,\calT\big)$
determines, in the unbounded picture of Kasparov's $K\!K$-theory, a well-defined
class
$$
\left<\mathcal{H},\pi,T\right>:=\big[\mathcal{H}\hat{\otimes}C(\mathcal{U}),\pi
\hat{\otimes}1,\calT\big]\;\in\;\KK_{i}\big(A,C(\mathcal{U})\big)\,.
$$
In particular, $\mathcal{U}$ being path-connected, for every choice of points $u,v
\in\mathcal{U}$\,, the corresponding evaluation maps $(\ev_{\!u})_{*},\,(\ev_{\!v})_
{*}\colon\KK_{i}(A,C(\mathcal{U}))\longrightarrow\KK_{i}(A,\bbC)$ yield the same
$K$-homology class, \ie
$$
[\mathcal{H},\pi,T_{u}]=[\mathcal{H},\pi,T_{v}]\;\in\;\KK_{i}(A,\bbC)=K^{i}(A)\,.
$$
We can therefore denote this class unambiguously by $[\mathcal{H},\pi,T]\in K^{i}(A)$\,.
\end{Thm}

The proof is the subject of Section~\ref{s-proofs} below.


\section{Proof of Theorems \ref{thm-main} and \ref{thm-KK}}

\label{s-proofs}


%
In the present section, we establish Theorem~\ref{thm-KK}. Clearly, Theorem~\ref{thm-main}
is merely a part of it, so, we will not need to say more about its proof.

\smallskip

\PrfOf{Theorem~\ref{thm-KK}}
During the proof, we keep notation as in Section~\ref{s-KK}. The proof consists in two
steps. First, we have to show that the triple $\big(\mathcal{H}\hat{\otimes}C(\mathcal{U}),
\pi\hat{\otimes}1,\calT\big)$ indeed is a Kasparov triple in Baaj-Julg's unbounded
description of $\KK$-theory. The second step is simply the observation that the compact
Hausdorff space $\mathcal{U}$ being path-connected, the evaluation maps $(\ev_{\!u})_{*}$
$(\ev_{\!v})_{*}$ as in the statement yield, as is well-known, the same $K$-homology class.
So, we can focus exclusively on the first step. According to Baaj-Julg~\cite{BJ}, we have
to show that
\begin{itemize}
    \item [(1)] the domain $\dom(\calT)$ is dense;
    \item [(2)] the operator $\calT$ is self-adjoint;
    \item [(3)] the operator $\calT$ is regular;
    \item [(4)] the operator $[\pi(b)\hat{\otimes}1,\calT]$ is densely defined and
    bounded, for every $b\in\mathcal{B}$\,;
    \item [(5)] the operator $(\pi(a)\hat{\otimes}1)(\calT^{2}+1)^{-1}$ is compact, for
    every $a\in A$\,.
\end{itemize}
In~(5), compactness is meant in the sense of Hilbert $C^{*}$-modules. Let us now establish these
properties. (Concerning~(4), see also Remark~\ref{rem-black} below.)

\smallskip

\noindent
(1) The domain $\dom(\calT)$ contains $\dom(T\hat{\otimes}1)$\,, so, it is dense.

\smallskip

\noindent
(2) By definition, we have
$$
\dom(\calT^{*})=\big\{g\in C(\mathcal{U},\mathcal{H})\,\big|\,\exists h
\in C(\mathcal{U},\mathcal{H})\;\,\mbox{st.}\;\left<\calT\!f\,|\,g\right>=
\left<f\,|\,h\right>\,,\;\forall f\in\dom(\calT)\big\}\,.
$$
The condition $\left<\calT\!f\,|\,g\right>=\left<f\,|\,h\right>$ amounts to
having $\left<T_{u}f(u)\,|\,g(u)\right>_{\mathcal{H}}=\left<f(u)\,|\,h(u)\right>_
{\mathcal{H}}$ for every $u\in\mathcal{U}$\,. If $g\in\dom(\calT)$\,, then
$\left<T_{u}f(u)\,|\,g(u)\right>_{\mathcal{H}}=\left<f(u)\,|\,T_{u}g(u)\right>_
{\mathcal{H}}$ for every $u\in\mathcal{U}$\,, which means that $\left<\calT
f\,|\,g\right>=\left<f\,|\,\calT\!g\right>$\,. This shows that $\calT
\subseteq\calT^{*}$\,. We pass to the reverse inclusion. Since for every
$\xi\in\dom(T_{u})$\,, there exists $f\in\dom(\calT)$ with $f(u)=\xi$ (take for
$f$ the constant map), the condition for $g$ to be in $\dom(\calT^{*})$ implies
that $\left<T_{u}\xi\,|\,g(u)\right>_{\mathcal{H}}=\left<\xi\,|\,h(u)\right>_{\mathcal
{H}}$ for every  $\xi\in\dom(T_{u})$\,. Since $T_{u}$ is self-adjoint, by definition,
this means that $g(u)\in\dom(T_{u})$ and $h(u)=T_{u}g(u)$\,. Since this has to hold
for every $u\in\mathcal{U}$\,, we see that $g\in\dom(\calT)$\,, so that
$\dom(\mathcal{T^{*}})\subseteq\dom(\calT)$\,.

\smallskip

\noindent
(3) By~\cite[Lem.~9.8]{Lance}, to show that $\calT$ is regular, we just have to check that
$\calT+i$ and $\calT-i$ are surjective. Let $u\in\mathcal{U}$\,; since $T_{u}$ is
self-adjoint, the operator $T_{u}{\pm i}$ is surjective, and $(T_{u}{\pm i})^{-1}$ is a bounded
operator with range equal to $\dom(T_{u})$\,. So, we can consider the $C(\mathcal{U})$-linear operator
$$
\mathcal{S}_{\pm}\colon C(\mathcal{U},\mathcal{H})\longrightarrow C(\mathcal
{U},\mathcal{H})\,,\quad f\longmapsto\big(u\mapsto(T_{u}{\pm i})^{-1}f(u)\big)\,.
$$
In fact, we have to prove that $\mathcal{S}_{\pm}$ is really well-defined, namely
that the function $\mathcal{S}_{\pm}f\colon u\longmapsto(T_{u}{\pm i})^{-1}f(u)$ is
continuous. We do this just below and assume it for a while. Since $\Ima((T_{u}\pm
i)^{-1})=\dom(T_{u})$ for every $u\in\mathcal{U}$\,, we see that
$$
\dom(\calT\mathcal{S}_{\pm})=\big\{f\in C(\mathcal{U},\mathcal{H})\,\big|\,
(T_{u}{\pm i})^{-1}f(u)\in\dom(T_{u})\,,\;\forall u\in\mathcal{U}\big\}=C(\mathcal
{U},\mathcal{H})\,.
$$
Since, obviously, $(\calT{\pm i})\mathcal{S}_{\pm}=1$ on $\dom(\calT
\mathcal{S}_{\pm})$\,, the operator $\calT{\pm i}$ is surjective. So, let
us establish the continuity of $\mathcal{S}_{\pm}f$\,. Since by assumption $T$
admits at least one self-adjoint extension, the operator $T{\pm i}$ is injective
on its domain $\dom(T)$\,, therefore, the operator $(T{\pm i})^{-1}\colon\Ima
(T{\pm i})\longrightarrow\dom(T)\subseteq\mathcal{H}$ is well-defined on its
domain $\Ima(T{\pm i})$\,. It will be crucial for us to observe that $\Ima(T{\pm i})$
is closed, as follows from \cite[Prop.~X.2.5\,(c)]{Conw}; as a side-remark, note
also that $\Ima(T{\pm i})$ is the whole of $\mathcal{H}$ if and only if $T$ is self-adjoint,
see \cite[Thm.~VIII.2, pp.~256--257]{RSI}. Fix a point $u\in\mathcal{U}$\,. Consider
the closed subspace $\mathcal{H}_{u}:=\big\{\eta+u(\eta)\,\big|\,\eta\in\Ker(T^{*}{-i})
\big\}$ of $\mathcal{H}$ and the operator
$$
R_{u}\colon\mathcal{H}_{u}\longrightarrow\Ker(T^{*}{-i})\oplus\Ker(T^{*}{+i})\subseteq
\mathcal{H}\,,\quad\eta+u(\eta)\longmapsto i\eta-iu(\eta)\,,
$$
where the direct sum is an algebraic one, that is, of mere vector spaces.
For $\eta\in\Ker(T^{*}{-i})$\,, we compute that $(R_{u}+i)(\eta+u(\eta))=2i\eta\in
\Ker(T^{*}{-i})$ and that $(R_{u}-i)(\eta+u(\eta))=-2iu(\eta)\in\Ker(T^{*}{+i})$\,. So,
we can view $R_{u}{\pm i}$ as an operator with codomain $\Ker(T^{*}{\mp i})$\,, \ie
$$
R_{u}{\pm i}\colon\mathcal{H}_{u}\longrightarrow\Ker(T^{*}{\mp i})\,.
$$
As a consequence, $T_{u}{\pm i}$ decomposes as an algebraic direct sum of two operators, as follows\,:
$$
T_{u}{\pm i}=(T{\pm i})\oplus(R_{u}{\pm i})\colon\underbrace{\dom(T)\oplus\mathcal{H}_{u}}_{=\dom
(T_{u})}\longrightarrow\underbrace{\Ima(T{\pm i})\oplus\Ker(T^{*}{\mp i})}_{=\mathcal{H}}\,.
$$
(The last direct sum is orthogonal, but we will not need this fact.) From the above explicit
computation of $R_{u}{\pm i}$\,, we deduce that
$$
{\textstyle (R_{u}{\pm i})^{-1}\colon\Ker(T^{*}{\mp i})\longrightarrow\mathcal{H}_{u}
\subseteq\mathcal{H}\,,\quad\xi\longmapsto\pm\frac{1}{2i}\big(\xi+u^{\pm 1}(\xi)\big)\,.}
$$
So, we can write the operator $(T_{u}{\pm i})^{-1}\colon\mathcal{H}\longrightarrow\dom
(T_{u})\subseteq\mathcal{H}$ as the direct sum
$$
(T{\pm i})^{-1}\oplus(R_{u}{\pm i})^{-1}\colon\Ima(T{\pm i})\oplus\Ker(T^{*}{\mp i})\longrightarrow
\dom(T)\oplus\mathcal{H}_{u}\subseteq\mathcal{H}\,.
$$
Finally, letting $P_{\pm}\colon\mathcal{H}\onto\Ima(T{\pm i})$ and $Q_{\pm}\colon\mathcal{H}\onto\Ker(T^{*}{\mp i})$
denote the ortho\-gonal projections (recall that $\Ima(T{\pm i})$ is closed\,!), we see that
$$
{\textstyle \mathcal{S}_{\pm}f\colon u\longmapsto(T{\pm i})^{-1}P_{\pm}f(u)\pm\frac{1}{2i}
\Big(Q_{\pm}f(u)+u^{\pm 1}\big(Q_{\pm}f(u)\big)\Big)\,.}
$$
Observe that the function $u\longmapsto u^{\pm 1}\big(Q_{\pm}f(u)\big)$ is the composition
$$
\arraycolsep1pt
\renewcommand{\arraystretch}{1.2}
\begin{array}{rcccl}
  \mathcal{U} & \longrightarrow & \mathcal{U}\times\Ker(T^{*}{\mp i}) & \longrightarrow &
  \Ker(T^{*}{\pm i})\subseteq\mathcal{H} \\
  u & \longmapsto & (u,Q_{\pm}f(u)) & & \\
   & & (v,\xi) & \longmapsto & v^{\pm 1}(\xi)\,,
\end{array}
$$
where the two indicated maps are continuous (for the latter, recall that $\mathcal{U}$
is equipped with the norm-topology). It follows that $\mathcal{S}_{\pm}f$ is continuous,
as was to be shown.

\smallskip

\noindent
(4) Fix an element $b\in\mathcal{B}$\,. Since $[\pi(b)\hat{\otimes}1,T\hat{\otimes}1]\subseteq[\pi(b)\hat
{\otimes}1,\calT]\subseteq [\pi(b)\hat{\otimes}1,T^{*}\hat{\otimes}1]$\,, the result follows.

\smallskip

\noindent
(5) Finally, we fix an element $a\in A$\,. By assumption, the operator $T^{*}T+1$ has dense range and $\pi(a)(T^{*}T+1)^
{-1}$ is compact. Since $T_{u}^{2}+1$ is an extension of $T^{*}T+1$\,, it
follows that $\pi(a)(T_{u}^{2}+1)^{-1}$ coincides with the compact operator
$\pi(a)(T^{*}T+1)^{-1}$\,. Consequently, we have
$$
(\pi(a)\hat{\otimes}1)(\calT^{2}+1)^{-1}=\pi(a)(T^{*}T+1)^{-1}\hat{\otimes}
1\;\in\;\mathcal{K}(\mathcal{H})\hat{\otimes}C(\mathcal{U})=\mathcal{K}\big(\mathcal{H}
\hat{\otimes}C(\mathcal{U})\big)
$$
(see~\cite[p.~10]{Lance} for the final equality), and the proof is complete.
\qed

\begin{Rem}
\label{rem-black}
Following Blackadar's treatment of the Baaj-Julg results, we did not require $[\pi(b)\hat{\otimes}1,
\mathcal{T}\,]$ to have domain containing $\dom(\mathcal{T})$ for every $b\in\mathcal{B}$\,, but merely
to have dense domain, see~\cite[pp.~163--165]{Bla}.
\end{Rem}


\section{Generalization of Theorem \ref{thm-main} and proof of Theorem \ref{thm-main2}}

\label{s-general}


%
We start this section with some observations from which we derive a generalization of
Theorem~\ref{thm-main}.

\begin{Obs}
\label{rem-improve}
\mbox{}\\
\vspace*{-1.2em}
\begin{itemize}
    \item [(1)] The assumption that $T$ is closed is not really essential, since
    otherwise one can simply replace it by its closure $\bar{T}=T^{**}$\,, and
    then check/require properties {\bf (i)}, {\bf (ii)} and {\bf (iii)} of Definition~\ref{def-triple}
    for the closure.
    \item [(2)] The condition, weaker than {\bf (i)} of Definition~\ref{def-triple}, saying that
    the deficiency indices of $T$ coincide, but are not necessarily finite is enough to define
    $[\mathcal{H},\pi,T]\in K^{i}(A)$ unambiguously. Indeed, the space $\mathcal{U}$
    is always path-connected, so, we replace it everywhere by a path $\mathcal{P}_
    {\!uv}$ connecting two arbitrary points $u$ and $v$ in $\mathcal{U}$\,. For this,
    note that $\mathcal{P}_{\!uv}$ is a non-empty, compact Hausdorff and path-connected
    space, and that the map $\mathcal{P}_{\!uv}\times\Ker(T^{*}{\mp i})\longrightarrow\Ker(T^{*}
    {\pm i})$ taking $(v,\xi)$ to $v^{\pm 1}(\xi)$ is also continuous. The proof is `less canonical'
    in this case (since we are constrained to make a choice for the path $\mathcal{P}_{\!uv}$).
    \item [(3)] Let $T$ be a densely defined closed symmetric operator on a Hilbert space
    $\mathcal{H}$\,, and let $\pi\colon A\longrightarrow B(\mathcal{H})$ be a $*$-homomorphism.
    Then the following property -- which does \emph{not} involve $T^{*}$ -- implies {\bf (iii)}
    of Definition~\ref{def-triple}\,:
    \begin{itemize}
       \item [{\bf (iii$^{\prime}$)}] $[\pi(b),T]$ and $[\pi(b^{*}),T]$ are densely defined
       and $[\pi(b),T]$ is bounded for every $b$ in a dense subspace $\mathcal{B}$ of $A$\,.
    \end{itemize}
    Indeed, to show that {\bf (iii)} of~\ref{def-triple} follows from {\bf (iii$^{\prime}$)},
    we first note that $[\pi(b),T^{*}]$ is densely defined and also closable, since its adjoint
    satisfies
    $$
    \qquad\quad\;\;[\pi(b),T^{*}]^{*}=(\pi(b)T^{*}-T^{*}\pi(b))^{*}\supseteq T^{**}\pi(b)^
    {*}-\pi(b)^{*}T^{**}=-[\pi(b^{*}),T]\,,
    $$
    so, is densely defined \cite[Prop.~X.1.6\,(b) \& Ex.~X.1.1, pp.\ 305 \& 308]{Conw}. Using the
    inclusion $T\subseteq T^{*}$\,, we get $[\pi(b),T^{*}]^{*}\subseteq[\pi(b),T]^{*}$\,. Since
    by assumption $[\pi(b),T]$ is densely defined and bounded, $[\pi(b),T]^{*}\in B(\mathcal{H})$\,.
    Therefore, $[\pi(b),T^{*}]^{*}$ is densely defined and bounded, so that the closure of $[\pi(b),
    T^{*}]$ satisfies $[\pi(b),T^{*}]^{**}\in B(\mathcal{H})$\,, showing that $[\pi(b),T^{*}]$ is,
    indeed, densely defined and bounded.
    \item [(4)] Let $[\mathcal{H},\pi,D]\in\KK_{i}(A,\bbC)$ be an unbounded
    Fredholm module in the usual sense, \ie with $D$ self-adjoint. Let $T$ be a densely defined
    closed symmetric restriction of $D$\,. Then, the deficiency indices of $T$ are automatically
    equal, so that $T$ satisfies {\bf (i)} provided one of them is finite; this happens exactly
    when the quotient $\dom(T^{*})/\dom(T)$ is finite dimensional. Furthermore, $T$ necessarily
    verifies {\bf (ii)}, since then $(T^{*}T+1)^{-1}=(D^{2}+1)^{-1}$\,, and $\pi(a)(D^{2}+1)^{-1}$
    is compact for every $a\in A$ by assumption.
\end{itemize}
\end{Obs}

These observations combined with Theorem~\ref{thm-main} lead us directly to the following statement.

\begin{Thm}
\label{thm-gen}
Let $A$ be a separable $C^{*}$-algebra. Suppose given the following data\,:
\begin{itemize}
    \item [(a)] an integer $i\in\{0,1\}$\,;
    \item [(b)] a Hilbert space $\mathcal{H}$\,;
    \item [(c)] a $*$-representation $\pi\colon A\longrightarrow B(\mathcal{H})$ of $A$\,;
    \item [(d)] a densely defined symmetric operator $T$ on $\mathcal{H}$ with domain $\dom(T)$\,.
\end{itemize}
These data are required to fulfill, firstly, the two conditions
\begin{itemize}
    \item [{\bf (i)}] the deficiency indices of $T^{**}$ coincide (as cardinals);
    \item [{\bf (ii)}] the operator $\pi(a)(T^{*}T^{**}+1)^{-1}$ is compact for every $a\in A$\,;
\end{itemize}
and, secondly, one of the following two conditions\,:
\begin{itemize}
    \item [{\bf (iii)}] the operator $[\pi(b),T^{**}]$ is densely defined, and $[\pi(b),T^{*}]$ is
    densely defined and bounded for every $b$ in some norm-dense subspace $\mathcal{B}$ of $A$\,;
    \item [{\bf (iii$^{\prime}$)}] $[\pi(b),T^{**}]$ and $[\pi(b^{*}),T^{**}]$ are densely defined
    and $[\pi(b),T^{**}]$ is bounded for every $b$ in a dense subspace $\mathcal{B}$ of $A$\,.
\end{itemize}
Thirdly, if $i=0$\,, we moreover require $\mathcal{H}$ to be $\bbZ/2$-graded, $\pi$ to preserve the
grading, and $T$ to reverse it. Then, there exists, in the unbounded picture of
analytical $K$-homology,
a well-defined $K$-homology class
$$
[\mathcal{H},\pi,T]\in K^{i}(A)\,,
$$
that is canonical, and coincides with the usual class in case $T$ is self-adjoint. More precisely, given an
arbitrary self-adjoint extension $\widetilde{T}$ of $T$ -- and at least one such extension exists --, one
has
$$
[\mathcal{H},\pi,T]=[\mathcal{H},\pi,\widetilde{T}]\;\in\;K^{i}(A)\,,
$$
independently of the choice of $\widetilde{T}$\,.\qed
\end{Thm}

\begin{Def}
By extension, we call a triple $(\mathcal{H},\pi,T)$ satisfying the hypotheses of Theorem~\ref{thm-gen}
a \emph{symmetric unbounded Fredholm module of degree $i$} over $A$\,.
\end{Def}

We pass to the proof of Theorem~\ref{thm-main2}.

\smallskip

\PrfOf{Theorem~\ref{thm-main2}}
We proceed somehow as in~\ref{rem-improve}\,(4). Since $T_{1}$ is self-adjoint, its
restriction $T$ is symmetric with closure $T^{**}$ satisfying $T\subseteq T^{**}\subseteq T_{1}$
and being symmetric. By hypothesis (a), $T$ is densely defined, therefore, so is $T^{**}$\,. In
particular, the densely defined closed symmetric operator $T^{**}$ admits at least one self-adjoint
extension, namely $T_{1}$\,, so that its deficiency indices coincide (see Section~\ref{s-KK}). This
shows that $T$ satisfies {\bf (i)} of Theorem~\ref{thm-gen}. Moreover, one has
$$
(T^{*}T^{**}+1)^{-1}=(T_{1}^{2}+1)^{-1}\,.
$$
By assumption that $(\mathcal{H},\pi,T_{1})$ is a (usual) Fredholm module, $\pi(a)(T_{1}^{2}+1)^{-1}$
is compact for every $a\in A$\,. This implies {\bf (ii)} of \ref{thm-gen} for $T$\,. By hypothesis (b),
we have that both the operators $[\pi(b),T^{**}]$ and $[\pi(b^{*}),T^{**}]$ are densely defined and
bounded for every $b$ in $\mathcal{B}$ (recall that $\mathcal{B}$ is $*$-closed). This implies {\bf
(iii$^{\prime}$)} of \ref{thm-gen} for $T$\,. All in all, we have a symmetric unbounded Fredholm module
$(\mathcal{H},\pi,T)$\,, and, applying Theorem~\ref{thm-gen} to it twice ($T_{2}$ is a self-adjoint
extension of $T$ as well), we get the equalities
$$
[\mathcal{H},\pi,T_{1}]=[\mathcal{H},\pi,T]=[\mathcal{H},\pi,T_{2}]
$$
in $K^{i}(A)$\,, as desired. It remains to prove that condition~(b$^{\prime}$) implies condition~(b). First,
the operators $[\pi(b'),T^{**}]$ and $[\pi(b'^{*}),T^{**}]$ are densely defined and bounded for every $b'$
in the subspace $\mathcal{B'}$\,. Secondly, for every $b'\in\mathcal{B'}$\,, the operator $[\pi(b'),T^{**}]$
is a restriction of $[\pi(b'),T_{1}]$\,, and is therefore also bounded, by choice of $\mathcal{B'}$\,. So,
we get condition~(b) with $\mathcal{B}:=\mathcal{B}'$\,. The proof is now complete.
\qed

%


\part{Proofs of the ``addition formulae''}

\label{p-III}


\section{Topological proof of the ``addition formulae''}

\label{s-mfld-top-proof}


This section is subdivided into three subsections. In the first one, we establish a
useful general principle that will allow us to reduce the proofs (both topological
and analytical) of Theorems \ref{thm-add-formula-1} and \ref{thm-add-formula-2} to
the verification of one equality that embodies the pure substance of the ``addition
formulae'', without any extraneous ornament. In the other two subsections, one
for each treated dimension, we prove these theorems in the topological setting.


\subsection{A general principle}

\label{ss-gen-princ}


\mbox{}\\
\vspace*{-.8em}

We present here a general, but easy, principle on homology (and related) theories, that will
be used on several occasions in the sequel, even for the analytical proofs. To state it, we
call a functor $F(-)$ from the category of CW-complexes to the category of abelian groups
\emph{additive} if, given two maps $f_{1}\colon X_{1}\longrightarrow X$ and $f_{2}\colon
X_{2}\longrightarrow X$ of CW-complexes, one has a natural isomorphism
$$
F(X_{1}\amalg X_{2})\cong F(X_{1})\oplus F(X_{2})
$$
such that, using it as an identification,
$$
(f_{1}\amalg f_{2})_{*}(x_{1},x_{2})=(f_{1})_{*}(x_{1})+(f_{2})_{*}(x_{2})\;\in\;F(X)\,,
$$
for every $x_{1}\in F(X_{1})$ and $x_{2}\in F(X_{2})$\,, where $f_{*}$ stands for $F(f)$
whenever $f$ is a map between CW-complexes. For example, an additive homology theory with
compact supports, like integral homology or $K$-homology, is an additive functor. The point
for us is that the assignment
$$
X\longmapsto\RKK_{*}(X,\bbC)
$$
is straightforwardly seen to be an additive functor in our sense, \emph{without using the identification
of $\RKK_{*}(X,\bbC)$ with $K_{*}(X)$}, so that, later, our analytical proofs will really be purely
and strictly analytical.

\begin{Lem}
\label{lem-gen-princ}
Let $F(-)$ be an additive functor as defined above. Let $M_{1}$ and $M_{2}$ be connected
oriented manifolds of the same dimension $n>0$\,. Let $D_{i}$ ($i=1,2$) be a `small' embedded open
disk in $M_{i}$\,, whose boundary inside $M$ contains the base-point, and form the oriented connected
sum $M_{1}\# M_{2}$ by gluing $M_{1}\!\smallsetminus D_{1}$ and $M_{2}\!\smallsetminus D_{2}$ along
their boundaries. Consider the obvious maps
$$
j\colon M_{1}\amalg M_{2}\longrightarrow M_{1}\vee M_{2}\qquad\mbox{and}\qquad
p\colon M_{1}\# M_{2}\longrightarrow M_{1}\vee M_{2}\,,
$$
given by identification of the base-points and pinching the boundary $\partial D_{1}\approx\partial D_{2}$
to a point, respectively. Let $X$ be a pointed CW-complex and let $f_{i}\colon M_{i}\longrightarrow X$
be a continuous map, which, on $D_{i}$\,, is constant and equal to the base-point of $X$\,. Consider
the connected sum $f_{1}\# f_{2}\colon M_{1}\# M_{2}\longrightarrow X$\,. Finally, suppose given three
elements
$$
x_{1}\in F(M_{1})\,,\quad x_{2}\in F(M_{2})\quad\mbox{and}\quad x\in F(M_{1}\# M_{2})
$$
satisfying the compatibility condition
$$
j_{*}(x_{1},x_{2})=p_{*}(x)\;\in\;F(M_{1}\vee M_{2})\,.
$$
Then, one has the equality
$$
(f_{1}\# f_{2})_{*}(x)=(f_{1})_{*}(x_{1})+(f_{2})_{*}(x_{2})\;\in\;F(X)\,.
$$
\end{Lem}

\Prf
We have the commutative diagram
\begin{displaymath}
 \begin{diagram}[height=2em,width=3em]
    & & & & M_{1}\# M_{2} \\
    & & & & \dTo_{p} \\
   M_{1}\amalg M_{2} & & \rTo^{j} & & M_{1}\vee M_{2} \\
    & \rdTo_{f_{1}\amalg f_{2}} & & \ldTo_{f_{1}\vee f_{2}} & \\
    & & X & & \\
 \end{diagram}
\end{displaymath}
Noticing that $f_{1}\# f_{2}=(f_{1}\vee f_{2})\circ p$\,, we compute
$$
\arraycolsep1pt
\renewcommand{\arraystretch}{1.4}
\begin{array}{rcl}
    (f_{1}\# f_{2})_{*}(x) & = & (f_{1}\vee f_{2})_{*}\,p_{*}(x)  \\
     & = & (f_{1}\vee f_{2})_{*}\,j_{*}(x_{1},x_{2})  \\
     & = & (f_{1}\amalg f_{2})_{*}(x_{1},x_{2})  \\
     & = & (f_{1})_{*}(x_{1})+(f_{2})_{*}(x_{2})\,,
\end{array}
$$
and this completes the proof.
\qed

\smallskip

As a prototypical illustration of Lemma~\ref{lem-gen-princ}, we deduce the following simple
example on the homology of manifolds.

\begin{Ex}
\label{ex-hmlgy-mfld}
Keep notation as in Lemma~\ref{lem-gen-princ}, but assume $M_{i}$ ($i=1,2$) to be closed and denote
by $[M_{i}]\in H_{n}(M_{i};\bbZ)$ its fundamental class. Then, in the group $H_{n}(X;\bbZ)$\,, one has
$$
(f_{1}\# f_{2})_{*}[M_{1}\# M_{2}]=(f_{1})_{*}[M_{1}]+(f_{2})_{*}[M_{2}]\,.
$$
Indeed, the map $j_{*}\colon H_{n}(M_{1}\amalg M_{2};\bbZ)\longrightarrow H_{n}(M_{1}\vee M_{2};\bbZ)$
satisfies the compatibility condition
$$
j_{*}\big([M_{1}],[M_{2}]\big)=p_{*}[M_{1}\# M_{2}]
$$
(as a computation using suitable triangulations shows), so, Lemma~\ref{lem-gen-princ} applies to give
the result. Now, suppose that $M_{1}=\Sigma_{g_{1}}$ and $M_{2}=\Sigma_{g_{2}}$ are closed oriented surfaces.
As we have noticed, the Euler characteristic is \emph{not} additive with respect to the connected sum.
Note that this amounts to saying that the corresponding compatibility condition is \emph{not} satisfied,
as next indicated\,:
$$
j_{*}\big(\chi(\Sigma_{g_{1}})\!\cdot\![1],\chi(\Sigma_{g_{2}})\!\cdot\![1]\big)\neq p_{*}\big(\chi(\Sigma_{g_{1}}
\#\Sigma_{g_{2}})\!\cdot\![1]\big)\;\in\;H_{0}(\Sigma_{g_{1}}\!\vee\Sigma_{g_{2}};\bbZ)\,,
$$
where $[1]$ stands for the prescribed generator of the zeroth homology group of any connected CW-complex.
\end{Ex}

The next lemma is a slight variation of Lemma~\ref{lem-gen-princ}; for the modified connected
sum ``\,$\natural$\,'', we refer to Figures 1 and 2 in Section~\ref{s-appl-mfld} and to the
statement of Theorem~\ref{thm-add-formula-2}.

\begin{Lem}
\label{lem-gen-princ-var}
Keep the same notation and hypotheses as in Lemma~\ref{lem-gen-princ}, but with $M_{1}=
\Sigma_{g_{1}}$ and $M_{2}=\Sigma_{g_{2}}$ being closed oriented surfaces, and $D_{i}$ being
a handle given as a small open tubular neighbourhood of a suitable non-retractable embedded
circle $C_{i}$ ($i=1,2$). Let $\Sigma_{g_{1}}\cup_{S^{1}}\Sigma_{g_{2}}$ denote the CW-complex
obtained as the union of $\Sigma_{g_{1}}$ and $\Sigma_{g_{2}}$ with the circles $C_{1}$
and $C_{2}$ pointwise identified in an orientation-preserving way. Then, the equality
$$
(f_{1}\natural f_{2})_{*}(x)=(f_{1})_{*}(x_{1})+(f_{2})_{*}(x_{2})\;\in F(X)
$$
holds for $x_{1}\in F(\Sigma_{g_{1}})$\,, $x_{2}\in F(\Sigma_{g_{2}})$ and $x\in F(\Sigma_
{g_{1}}\natural\Sigma_{g_{2}})$ satisfying the compatibility condition $j_{*}(x_{1},x_{2})
=p_{*}(x)$ in $F(\Sigma_{g_{1}}\cup_{S^{1}}\Sigma_{g_{2}})$\,, where $j$ and $p$ stand for
the obvious identification and pinching maps
$$
j\colon \Sigma_{g_{1}}\amalg\Sigma_{g_{2}}\longrightarrow\Sigma_{g_{1}}\cup_{S^{1}}\Sigma_
{g_{2}}\qquad\mbox{and}\qquad p\colon\Sigma_{g_{1}}\natural\Sigma_{g_{2}}\longrightarrow
\Sigma_{g_{1}}\cup_{S^{1}}\Sigma_{g_{2}}\,.
$$
\end{Lem}

\Prf
This time, we have the commutative diagram
\begin{displaymath}
 \begin{diagram}[height=2em,width=3em]
    & & & & \Sigma_{g_{1}}\natural\Sigma_{g_{2}} \\
    & & & & \dTo_{p} \\
   \Sigma_{g_{1}}\amalg\Sigma_{g_{2}} & & \rTo^{j} & & \Sigma_{g_{1}}\cup_{S^{1}}\Sigma_{g_{2}} \\
    & \rdTo_{f_{1}\amalg f_{2}} & & \ldTo_{f_{1}\cup_{S^{1}}f_{2}} & \\
    & & X & & \\
 \end{diagram}
\end{displaymath}
Noticing that $f_{1}\natural f_{2}=(f_{1}\cup_{S^{1}}f_{2})\circ p$\,, we can perform
a similar computation as to establish Lemma~\ref{lem-gen-princ}.
\qed

\smallskip

The next result will turn useful on several occasions later on.

\begin{Prop}
\label{prop-sum-var}
Let $X$ be a pointed CW-complex. For $i=1,2$\,, let $f_{i}\colon\Sigma_{g_{i}}\longrightarrow X$
be a pointed continuous map, that is constant in a small neighbourhood of a handle of $\Sigma_{g_{i}}$
($g_{i}\geq 0$). Consider the modified connected sum (along the two given handles) $f_{1}\natural
f_{2}\colon\Sigma_{g_{1}}\natural\Sigma_{g_{2}}\longrightarrow X$ of these two maps. Then, one has
the equality
$$
\chi(\Sigma_{g_{1}}\natural\Sigma_{g_{2}})=\chi(\Sigma_{g_{1}})+\chi(\Sigma_{g_{2}})
$$
and, in homology, one has
$$
(f_{1}\natural f_{2})_{*}[\Sigma_{g_{1}}\natural\Sigma_{g_{2}}]=(f_{1})_{*}[\Sigma_{g_{1}}]+
(f_{2})_{*}[\Sigma_{g_{2}}]\;\in\;H_{2}(X;\bbZ)\,.
$$
\end{Prop}

\Prf
The first equality is obvious, since $\chi(\Sigma_{g})=2-2g$ for any $g\geq 0$\,.
By the general principle~\ref{lem-gen-princ-var}, it suffices to check that
$$
j_{*}\big([\Sigma_{g_{1}}],[\Sigma_{g_{2}}]\big)=p_{*}[\Sigma_{g_{1}+g_{2}-1}]\;\in\;H_{2}
(\Sigma_{g_{1}}\cup_{S^{1}}\Sigma_{g_{2}};\bbZ)\,.
$$
One can either decree it to be visible and hence obvious, or, for instance, draw concrete triangulations
on both surfaces, respecting the prescribed circles and handles, and then determine the corresponding
triangulations for $\Sigma_{g_{1}}\natural\Sigma_{g_{2}}$ and $\Sigma_{g_{1}}\cup_{S^{1}}\Sigma_{g_{2}}$\,,
and finally deduce explicitly the maps $j_{*}$ and $p_{*}$ on the $H_{2}$-level to get the result.
\qed


\subsection{Topological proof of Theorem \ref{thm-add-formula-1}}

\label{ss-top-prf-add-form-1}


\mbox{}\\
\vspace*{-.8em}

We begin by stating a result that will be needed for the topological proof of
Theorem~\ref{thm-add-formula-2} as well (and also in our application to the Baum-Connes
Conjecture).

\begin{Lem}
\label{lem3}
For a $2$-dimensional CW-complex $X$\,, there are canonical and natural isomorphisms
(so-called \emph{integral Chern characters})
$$
ch_{\even}^{\bbZ}\colon K_{0}(X)\stackrel{\cong}{\longrightarrow}H_{0}(X;\bbZ)\oplus H_{2}
(X;\bbZ)\quad\mbox{and}\quad ch_{\odd}^{\bbZ}\colon K_{1}(X)\stackrel{\cong}{\longrightarrow}
H_{1}(X;\bbZ)\,.
$$
They are compatible with the usual Chern character $ch_{*}\colon K_{*}(X)\longrightarrow
H_{*}(X;\bbQ)$ and the map in homology corresponding to the change of coefficients $\bbZ\hookrightarrow
\bbQ$\,.
\end{Lem}

The compatibility in question is the obvious one; in case of doubt, see the diagrams in the proof
of Lemma~\ref{lem4} below. See \cite[Lem.~4]{BBV} or \cite[Prop.~2.1]{MM2} for a proof of the lemma.
The proof actually reveals more. Indeed, if $X^{[1]}$ denotes the $1$-skeleton of $X$\,, consider the
exact sequences both in $K$-homology and ordinary homology, associated with the cofibre sequence of
spaces $X^{[1]}\longrightarrow X\longrightarrow X/X^{[1]}$\,. The isomorphisms in the lemma are the unique
isomorphisms making the following diagram with exact rows commute\,:
\begin{displaymath}
 \begin{diagram}[height=2em,width=2.3em]
     & & \cong\bbZ & & & & & & & & & \\
    0 & \rTo & \overbrace{K_{0}(X^{[1]})} & \rTo & K_{0}(X) & \rTo & \widetilde
    {K}_{0}(X/X^{[1]}) & \rTo & K_{1}(X^{[1]}) & \rTo & K_{1}(X) & \rTo & 0 \\
     & & \dEqto & & \dDashto & & \dTo_{\cong} & & \dTo_{\cong} & & \dDashto & \\
    0 & \rTo & \underbrace{H_{0}(X)} & \rTo & H_{\even}(X) & \rTo & H_{2}(X/X^{[1]}) &
    \rTo & H_{1}(X^{[1]}) & \rTo & H_{1}(X) & \rTo & 0 \\
     & & \cong\bbZ & & & & & & & & & \\
 \end{diagram}
\end{displaymath}

\noindent
where $H_{*}$ stands for \emph{integral} homology. (Note that $X^{[1]}$ and $X/X^{[1]}$ are homotopy
equivalent to a bouquet of circles and of $2$-spheres, respectively). We present now an example,
that we state as a lemma for later reference.

\begin{Lem}
\label{lem-ch-S1}
For the class $[D]\in K_{1}(S^{1})$ of the Dirac operator over the circle, one has $ch_{\odd}^{\bbZ}[D]
={-[S^{1}]}$ in $H_{1}(S^{1};\bbZ)$\,, where $[S^{1}]$ is the fundamental homology class of $S^{1}$ corresponding
to the selected orientation; in particular, $K_{1}(S^{1})\cong\bbZ$ is generated by $[D]$\,.
\end{Lem}

\Prf
This is well-known (up to the minus sign\,!), but since no proof seems to be available in the literature,
we provide one here. Let $(e_{n})_{n\in\bbZ}$ be the trigonometric basis of the Hilbert space $L^{2}(S^{1})$\,.
The phase $F$ of $D$\,, \ie the operator appearing in the polar decomposition $D=F\cdot|D|$\,, is given by
$F(e_{n})=e_{n}$ if $n>0$\,, $F(e_{0})=0$ and $F(e_{n})=-e_{n}$ if $n<0$\,. The homotopy $t\mapsto
F\cdot(D^{*}D)^{t/2}$ between $F$ and $D$ shows that $[F]=[D]$ in $K^{1}(C(S^{1}))$\,. Now, consider the rank
one (hence compact) perturbation $F_{\nnspace\circ}$ of $F$ taking the same values on the $e_{n}$'s, except
that $F_{\nnspace\circ}(e_{0})=e_{0}$\,; of course, $[F_{\nnspace\circ}]=[F]$ in $K^{1}(C(S^{1}))$\,. The operator
$F_{\nnspace\circ}$ is a self-adjoint involution and the corresponding projection $P:=\frac{1+F_{\nnspace\circ}}
{2}$ is the Toeplitz projection on the Hardy space $\mathcal{H}^{2}(S^{1})$ of $S^{1}$\,; indeed, $\mathcal
{H}^{2}(S^{1})$ is defined as the closed span of the set $\{e_{n}\}_{n\geq 0}$ in $L^{2}(S^{1})$\,. By
\cite[2.7.7, 2.7.9, 5.1.6 and pp.~213--214]{HigRoe}, this means that $[D]$\,, as an element of
$\Ext(C(S^{1}))\cong\KK_{1}(C(S^{1}),\bbC)$\,, corresponds to the Toeplitz extension of $C(S^{1})$
by the compact operators on $\mathcal{H}^{2}(S^{1})$ described in~\cite[(2.3.5)]{HigRoe}). Now,
given a unitary $u$ in $C(S^{1})$\,, consider the corresponding Toeplitz operator $T_{\!u}:=PuP$
(see~\cite[Def.~2.7.7]{HigRoe}). Then, for the canonical pairing (\ie the Kasparov product)
$$
\otimes\colon\underbrace{\KK_{1}(\bbC,C(S^{1}))}_{\cong K^{1}(S^{1})}\otimes_{\bbZ}\underbrace{\KK_{1}
(C(S^{1}),\bbC)}_{\cong K_{1}(S^{1})}\longrightarrow\underbrace{\KK_{0}(\bbC,\bbC)}_{\cong\bbZ}\,,\quad
(x,y)\longmapsto x\otimes y\,,
$$
writing $\Wind(u)$ for the winding number of $u$\,, one has
$$
[u]\otimes [D]=\Index(T_{\!u})=-\Wind(u)\,,
$$
where the first equality follows from \cite[Thm.~18.10.2]{Bla}, and the second from
\cite[Thm.~2.3.2]{HigRoe}. By Lemma~\ref{lem3} and by its well-known \emph{co}homological
counterpart (see for instance \cite[Lem.~5.1]{MM2}), one has
$$
K_{1}(S^{1})\stackrel{\!ch_{\odd}^{\bbZ}\!}{\cong}H_{1}(S^{1};\bbZ)\;\,\cong\;\,\bbZ\qquad\mbox{and}
\qquad K^{1}(S^{1})\stackrel{\!ch^{\odd}_{\bbZ}\!}{\cong}H^{1}(S^{1};\bbZ)\;\,\cong\;\,\bbZ\,.
$$
Consider the unitary $u_{\circ}=(z\mapsto z)$ in $C(S^{1})$\,, whose class $[u_{\circ}]$ is the standard
generator of $K^{1}(S^{1})$\,, \ie the one satisfying $ch^{\odd}_{\bbZ}[u_{\circ}]=[S^{1}]$ in $H^{1}(S^{1};
\bbZ)$\,. Then, one gets $[u_{\circ}]\otimes[D]=-1$\,. Altogether, this shows that $[D]$ indeed is a generator
of $K_{1}(S^{1})$ (another approach for this result is one based on the ideas of~\cite{BD}). Now, since
the Chern character in $K$-homology and in $K$-theory of finite CW-complexes is induced by a map
of spectra (in the sense of algebraic topology), for such a space $X$\,, there is a commutative
diagram
\begin{displaymath}
 \begin{diagram}[height=2em,width=3em]
   K^{*}(X)\otimes_{\bbZ}K_{*}(X) & \rTo^{\!\!\!\!\left<\,.\,,.\,\right>_{\!_{K}}} & \bbZ \\
   \dTo^{ch^{*}\otimes ch_{*}} & & \dInto \\
   H^{*}(X;\bbQ)\otimes_{\bbQ}H_{*}(X;\bbQ) & \rTo^{\;\,\left<\,.\,,.\,\right>} & \bbQ \\
 \end{diagram}
\end{displaymath}
where ``\,$\left<\,.\,,.\,\right>$\,'' stands for the Kronecker product, and
``\,$\left<\,.\,,.\,\right>_{\!_{K}}$\,'' denotes the usual pairing between
$K$-theory and $K$-homology (see~\cite{BD}). As is folklore (see however \cite{Atiyah}
and \cite[Section~6]{Kas1}), the diagram
\begin{displaymath}
 \begin{diagram}[height=2em,width=4.5em]
   \KK_{*}(\bbC,C(X))\otimes_{\bbZ}\KK_{*}(C(X),\bbC) & \rTo^{\;\otimes} & \KK_{*}(\bbC,\bbC) \\
   \dTo^{\cong} & &  \dTo_{\cong} \\
   K^{*}(X)\otimes_{\bbZ}K_{*}(X) & \rTo^{\!\!\!\!\left<\,.\,,.\,\right>_{\!_{K}}} & \bbZ \\
 \end{diagram}
\end{displaymath}
does also commute. Since the integral (co)homology of the circle injects inside its rational
(co)homology and since $\left<[S^{1}],[S^{1}]\right>=\left<1,[S^{1}]\cap[S^{1}]\right>=
\left<1,1\right>=1$ by very choice of both orientation classes (see \cite[VII.12.8]{Dold}
for the first equality), it follows that $ch_{\odd}^{\bbZ}[D]={-[S^{1}]}$\,.
\qed

\smallskip

Now, we state a result that might be of independent interest, and to which the proof
of Theorem~\ref{thm-add-formula-2} reduces.

\begin{Prop}
\label{prop-add-form-1}
Consider the obvious identification and pinching maps
$$
j\colon S^{1}\amalg S^{1}\longrightarrow S^{1}\vee S^{1}\qquad\mbox{and}\qquad
p\colon S^{1}\cong S^{1}\# S^{1}\longrightarrow S^{1}\vee S^{1}\,.
$$
Then, in $K$-homology, on has the equality
$$
j_{*}\big([D],[D]\big)=p_{*}[D]\;\in\;K_{1}(S^{1}\vee S^{1})\,.
$$
\end{Prop}

\Prf
By Lemma~\ref{lem3}, we can identify $K_{1}$ with $H_{1}$ for all the (one-dimensional) spaces
in sight in the statement, and therefore, applying Lemma~\ref{lem-ch-S1}, we are reduced to proving
that $j_{*}\big({-[S^{1}]},{-[S^{1}]}\big)=p_{*}\big({-[S^{1}]}\big)$ in $H_{1}(S^{1}\vee S^{1})$\,,
or equivalently that $j_{*}\big([S^{1}],[S^{1}]\big)=p_{*}[S^{1}]$\,, an equality that is obvious
(compare with Example~\ref{ex-hmlgy-mfld}).
\qed

\smallskip

Finally, we can pass to our first proof of Theorem~\ref{thm-add-formula-1}.

\smallskip

\TopPrfOf{Theorem~\ref{thm-add-formula-1}}
Our general principle embodied by Lemma~\ref{lem-gen-princ} implies that the result
follows from Proposition~\ref{prop-add-form-1}, that we have proven with purely topological
methods.
\qed


\subsection{Topological proof of Theorem \ref{thm-add-formula-2}}

\label{ss-top-prf-add-form-2}


\mbox{}\\
\vspace*{-.8em}

As in the preceding subsection, we present a result of independent interest to which
Theorem~\ref{thm-add-formula-2} boils down.

\begin{Prop}
\label{prop-add-form-2}
Let $\Sigma_{g_{1}}$ and $\Sigma_{g_{2}}$ be two closed oriented surfaces of genus $g_{1}$
and $g_{2}$ respectively. Consider the obvious identification and pinching maps
$$
j\colon \Sigma_{g_{1}}\amalg\Sigma_{g_{2}}\longrightarrow\Sigma_{g_{1}}\cup_{S^{1}}\Sigma_
{g_{2}}\quad\mbox{and}\quad p\colon\Sigma_{g_{1}+g_{2}-1}\cong\Sigma_{g_{1}}\natural\Sigma_
{g_{2}}\longrightarrow\Sigma_{g_{1}}\cup_{S^{1}}\Sigma_{g_{2}}\,.
$$
Then, in $K$-homology, on has the equality
$$
j_{*}\big([\bar{\partial}_{g_{1}}],[\bar{\partial}_{g_{2}}]\big)=p_{*}[\bar{\partial}_{g_{1}+g_{2}-1}]
\;\in\;K_{0}(\Sigma_{g_{1}}\cup_{S^{1}}\Sigma_{g_{2}})\,.
$$
\end{Prop}

Before we present the proof, we recall the following fundamental and classical result.

\begin{Lem}
\label{lem-Riemann-Roch}
For a closed oriented surface $\Sigma_{g}$ of genus $g$\,, denote by $\iotaX{\Sigma_{g}}$ the
inclusion of the base-point. Then, letting $[1]$ denote the canonical generator of
the group $K_{0}(pt)\cong\bbZ$\,, one has
$$
K_{0}(\Sigma_{g})\cong\bbZ^{2}\qquad\mbox{and}\qquad K_{1}(\Sigma_{g})\cong\bbZ^{2g}
$$
with $\iotaX{\Sigma_{g}}_{*}[1]$ and $[\bar{\partial}_{g}]$ as generators of $K_{0}(\Sigma_{g})$\,;
furthermore, one has
$$
ch_{\even}^{\bbZ}[\bar{\partial}_{g}]=(1-g)\!\cdot\![1]+[\Sigma_{g}]\;\in\; H_{0}(\Sigma_{g};\bbZ)\oplus
H_{2}(\Sigma_{g};\bbZ)\,,
$$
where $[1]\in H_{0}(\Sigma_{g};\bbZ)\cong\bbZ$ is the canonical generator, and $[\Sigma_{g}]\in H_{2}
(\Sigma_{g};\bbZ)$ is the fundamental class.
\end{Lem}

\Prf
Consider a compact K\"{a}hler manifold $M$\,. Let $ch_{\ev}\colon K_{0}(M)\longrightarrow H_{\even}(M;\bbQ)$
be the usual (``rational-valued'') Chern character. Let $\bar{\partial}_{M}$ be the Dolbeault operator
on $M$\,, and $[\bar{\partial}_{M}]$ its class in $K_{0}(M)$\,. The Riemann-Roch-Hirzebruch Formula for
$M$ (see~\cite[p.~29]{Sha}) says precisely that $ch_{\ev}[\bar{\partial}_{M}]$ is Poincar\'{e}-dual
to the Todd class ${\rm Td}(TM)\in H^{*}(M;\bbQ)$\,. Specializing to $M=\Sigma_{g}$\,, we see that
$ch_{\ev}[\bar{\partial}_{g}]$ is Poincar\'{e}-dual to the rational cohomology class (see~\cite[p.~3]{Sha})
$$
{\textstyle {\rm Td}(T\Sigma_{g})=1+\frac{1}{2}c_{1}(T\Sigma_{g})}\;\in\;H^{\ev}(\Sigma_{g};\bbQ)\,;
$$
here, $T\Sigma_{g}$ is of course viewed as a complex line bundle over $\Sigma_{g}$\,. Since
$$
{\textstyle \frac{1}{2}\big<c_{1}(T\Sigma_{g}),\,[\Sigma_{g}]\big>=\Index(\bar{\partial}_{\Sigma_{g}})=1-g}
$$
(cf.~\cite[p.~27]{Sha}), the desired result concerning $ch_{\even}^{\bbZ}[\bar{\partial}_{g}]$ follows
from Poincar\'{e} duality, the fact that $ch_{\ev}$ and $ch_{\even}^{\bbZ}$ are compatible (see Lemma~\ref{lem3})
for $2$-dimensional spaces, and the fact that the integral homology of $\Sigma_{g}$ is torsion-free,
so that the canonical map from integral homology to rational homology is injective in this case.
The rest follows readily from Lemma~\ref{lem3} and the well-known integral homology of $\Sigma_{g}$\,.
\qed

\smallskip

\PrfOf{Proposition~\ref{prop-add-form-2}}
By means of Lemmas \ref{lem3} and \ref{lem-Riemann-Roch}, it suffices to check that
$$
j_{*}\big((\underbrace{1-g_{1}}_{\!\!\!\!=\frac{1}{2}\chi(\Sigma_{g_{1}})\!\!\!\!})\cdot[1],
(\underbrace{1-g_{2}}_{\!\!\!\!=\frac{1}{2}\chi(\Sigma_{g_{2}})\!\!\!\!})\cdot[1]\big)=p_{*}
\big((\underbrace{1-(g_{1}+g_{2}-1)}_{\!\!\!=\frac{1}{2}\chi(\Sigma_{g_{1}}
\natural\Sigma_{g_{2}})\!\!\!})\cdot[1]\big)
$$
in the group $H_{0}(\Sigma_{g_{1}}\cup_{S^{1}}\Sigma_{g_{2}};\bbZ)$\,, and that $j_{*}\big([\Sigma_{g_{1}}],
[\Sigma_{g_{2}}]\big)=p_{*}[\Sigma_{g_{1}+g_{2}-1}]$ in the group $H_{2}(\Sigma_{g_{1}}\cup_{S^{1}}
\Sigma_{g_{2}}; \bbZ)$\,. The first equality follows from Proposition~\ref{prop-sum-var}, and the second
from the proof of the latter.
\qed

\smallskip

Next, we present our first proof of Theorem~\ref{thm-add-formula-2}.

\smallskip

\TopPrfOf{Theorem~\ref{thm-add-formula-2}}
The modified version of our general principle, stated as Lemma~\ref{lem-gen-princ-var}, implies
that the result is a direct consequence of Proposition~\ref{prop-add-form-1}, whose proof given
above was performed in the topological setting.
\qed

\smallskip

\begin{Rem}
Both these topological proofs, though their relative simplicity, leave an unsatisfactory feeling for
the following reasons. It is not quite explicit \emph{here} what $K$-homology is, and in particular
how the classes $[D]$ and $[\bar{\partial}_{g}]$ are defined (except possibly in the proof of Lemma~\ref{lem-ch-S1}).
Moreover, the Chern character plays a rather mysterious rôle. Reasoning the opposite way, it is rather
pleasant that we did not have to define $K$-homology and these classes explicitly, and the question of
the Chern character is something well-understood and absolutely central in connection with the Atiyah-Singer
Index Theorem and of the Riemann-Roch-Hirzebruch Formula.
\end{Rem}

\begin{Rem}
Keeping notation as in Theorem~\ref{thm-add-formula-2}, we generally have
$$
(f_{1}\# f_{2})_{*}[\bar{\partial}_{g_{1}+g_{2}-1}]\neq(f_{1})_{*}[\bar{\partial}_{g_{2}}]+
(f_{2})_{*}[\bar{\partial}_{g_{2}}]\;\in\;K_{0}(X)\,,
$$
so that there is no ``addition formula'' for Dolbeault operators with respect to the usual
connected sum. Indeed, as the proof of the theorem presented above shows, the reason for
this is precisely that
$$
j_{*}\big((\underbrace{1-g_{1}}_{\!\!\!\!=\frac{1}{2}\chi(\Sigma_{g_{1}})\!\!\!\!})\cdot[1],
(\underbrace{1-g_{2}}_{\!\!\!\!=\frac{1}{2}\chi(\Sigma_{g_{2}})\!\!\!\!})\cdot[1]\big)\neq
p_{*}\big(\underbrace{1-(g_{1}+g_{2})}_{\!\!\!=\frac{1}{2}\chi(\Sigma_{g_{1}}\#\Sigma_{g_{2}})
\!\!\!}\cdot[1]\big)\;\in\;H_{0}(\Sigma_{g_{1}}\vee\Sigma_{g_{2}};\bbZ)\,,
$$
compare with Example~\ref{ex-hmlgy-mfld} and Proposition~\ref{prop-sum-var}. In other words,
the obstruction at the source of the problem is the non-additivity of the Euler characteristic
with respect to the usual connected sum.
\end{Rem}


\section{Analytical proof of the Dirac-type ``addition formula''}

\label{s-mfld-proofs-1}


This section is partitioned into two subsections. In the first, we describe, with the
necessary details, the class defined by the Dirac operator on the circle in $K$-homology,
viewed using the unbounded picture of analytical $K$-homology. In the second, one of the cores of the paper, we prove Theorem~\ref{thm-add-formula-1}
in this setting.


\subsection{Class of the Dirac operator for $S^{1}$ in analytic $K$-homology}

\label{ss-class-D-1}


\mbox{}\\
\vspace*{-.8em}

Consider the Dirac operator on the circle $S^{1}$ (equipped with the standard orientation,
more precisely the standard $\Spinc$-structure), namely
$$
D:=\frac{1}{i}\!\cdot\!\frac{d}{d\theta}\,,
$$
where $\frac{d}{d\theta}$ stands for the distributional derivative with domain
$$
{\textstyle \dom(D):=\big\{\xi\in L^{2}(S^{1})\,\big|\,\frac{d\xi}
{d\theta}\in L^{2}(S^{1})\;\,\mbox{and}\;\,\xi(0)=\xi(1)\big\}}\,.
$$
To be extremely precise, and for later use, let us give some explanations and
recall some basic and well-known facts. First, $[a,b]$ will denote an arbitrary
compact interval (with $a<b$), and $\theta$ (or $\theta_{0}$) a variable in it. We
consider $S^{1}$ as the unit interval $[0,1]$ with $0$ and $1$ identified, and we
view $L^{2}(S^{1})$ as $L^{2}[0,1]$ in the obvious way, namely considering a function
$\xi(e^{2\pi i\theta})$ as a function of the variable $\theta\in[0,1]$\,, denoted
by $\xi(\theta)$ for simplicity. Every class $\xi\in L^{2}[a,b]$ defines a distribution
on the interval $[a,b]$ given by
$$
{\textstyle T_{\xi}\varrho(\theta):=\int_{a}^{b}\xi(\theta)\varrho(\theta)\,d\theta
\qquad\quad(\varrho\in C^{\infty}[a,b])\,.}
$$
The \emph{distributional derivative} of $\xi$ is the distribution
$$
{\textstyle T_{\xi}^{\prime}\varrho(\theta):=-\int_{a}^{b}\xi(\theta)\frac{d\varrho}
{d\theta}(\theta)\,d\theta\qquad\quad(\varrho\in C^{\infty}[a,b])\,.}
$$
One writes $\frac{d\xi}{d\theta}\in L^{2}[a,b]$ if there exists a class $\eta\in
L^{2}[a,b]$ such that $T_{\eta}=T_{\xi}^{\prime}$ (and then, this class is unique).
In this case, we will always consider $\xi$ as being the unique \emph{continuous}
function on $[a,b]$ representing its class in $L^{2}[a,b]$\,; explicitly, it is
given by
$$
{\textstyle \xi(\theta)=\int_{a}^{\theta}\frac{d\xi}{d\theta}(\theta_{0})\,d\theta_{0}\,.}
$$
(This expression is meaningful by the Cauchy-Schwarz inequality.) We now define two kinds
of Sobolev spaces. Firstly, we set
$$
{\textstyle W^{1}[a,b]:=\big\{\xi\in L^2[a,b]\,\big|\,\frac{d\xi}{d\theta}\in L^2[a,b]\big\}\,.}
$$
Secondly, for the circle, we define
$$
{\textstyle W^{1}(S^{1}):=\big\{\xi\in W^{1}[0,1]\,\big|\,\xi(0)=\xi(1)\big\}\,.}
$$
Note that now the domain of $D$ has a transparent meaning, since by definition
it is the latter Sobolev space, \ie $\dom(D)=W^{1}(S^{1})$\,. Observe that $W^{1}[0,1]$
identifies with $AC[0,1]$\,, the space of absolutely continuous functions on $[0,1]$\,,
cf.\ \cite[pp.\ 258 \& 305]{RSI}. We keep this notation in the sequel. The
operator $D$ is self-adjoint on its domain $W^{1}(S^{1})$ (see e.g.\ \cite[Ex.~1 in X.1, p.~141]{RSII}),
as required to be part of the following unbounded Fredholm module\,:
$$
[D]:=[L^2(S^{1}),\mathcal{M},D]\;\in\;\KK_{1}(C(S^{1}),\bbC)\,,
$$
where $\mathcal{M}$ is the $*$-representation of $C(S^{1})$ on $L^2(S^{1})$ by pointwise multiplication.


\subsection{Analytical proof of Theorem \ref{thm-add-formula-1}}

\label{ss2.6}


\mbox{}\\
\vspace*{-.8em}

Recalling the notation introduced in the previous subsection, we can directly move
to the announced proof.

\smallskip

\AnPrfOf{Theorem~\ref{thm-add-formula-1}}
By our general principle~\ref{lem-gen-princ}, and keeping the same notation, it boils
down to proving that
$$
j_{*}\big([D],[D]\big)=p_{*}[D]\;\in\;\KK_{1}(C(S^{1}\vee S^{1}),\bbC)\,.
$$
Of course, this is precisely the content of Proposition~\ref{prop-add-form-1},
that we have already proved, but using topology. So, here, as promised, we
will reprove this in the realm of analytical $K$-homology. Let us now make the unbounded
Fredholm modules $p_{*}[D]$ and $j_{*}\big([D],[D]\big)$ explicit. First,
the pinching map $p\colon S^{1}\longrightarrow S^{1}\vee S^{1}$ is given
by
$$
p(\theta)=
\arraycolsep1pt
\left\{\begin{array}{ll}
    (2\theta)_{1}\,,\; & \mbox{if}\;\,0\leq\theta\leq\frac{1}{2}  \\
    (2\theta-1)_{2}\,,\; & \mbox{if}\;\,\frac{1}{2}\leq\theta\leq 1\,.
\end{array}\right.
$$
Here, $(2\theta)_{1}$ means that we view $2\theta$ as living in the
first copy of $S^{1}$ in $S^{1}\vee S^{1}$\,, and analogously for
$(2\theta-1)_{2}$\,. Then, $p_{*}[D]\in\KK_{1}(C(S^{1}\vee S^{1}),\bbC)$
is described as
$$
p_{*}[D]=\big[L^2(S^{1}),\mathcal{M}',D\big]\,,
$$
where, for $f\in C(S^{1}\vee S^{1})$\,, $\mathcal{M}'(f)=\mathcal{M}(f\circ
p)$\,, the multiplication by $f\circ p$ on $L^{2}(S^{1})$\,, and, as before,
$D$ is $\frac{1}{i}\!\cdot\!\frac{d}{d\theta}$ with domain $W^{1}(S^{1})$\,.
On the other hand,
$$
j_{*}\big([D],[D]\big)=\big[L^2(S^{1})\oplus L^2(S^{1}),\mathcal{M}_{1}\oplus\mathcal
{M}_{2},D\oplus D\big]\,,
$$
where the direct sum is an orthogonal direct sum, and, for $f\in C(S^{1}\vee S^{1})$ and
$i=1,2$\,, we have $\mathcal{M}_{i}(f)=\mathcal{M}(f_{i})$ with $f_{i}$ standing for the
restriction of $f$ to the $i$-th copy of $S^{1}$ in $S^{1}\vee S^{1}$\,. Here and below,
we make the obvious identifications
$$
{\textstyle L^{2}(S^{1})\oplus L^{2}(S^{1})=L^{2}(S^{1}\vee S^{1})=L^{2}
[0,1]\oplus L^{2}[1,2]=L^{2}[0,2]\,.}
$$
Since the domains of the operators play a crucial rôle, let us give
the domain of $D\oplus D$ very explicitly\,:
$$
\renewcommand{\arraystretch}{1.2}
\arraycolsep1pt
\begin{array}{rcl}
  \dom(D\oplus D) & = & W^{1}(S^{1})\oplus W^{1}(S^{1}) \\
   & = & \bigg\{(\xi_{1},\xi_{2})\in L^{2}[0,2]\,\bigg|
   \begin{array}{rcl}
     \xi_{1}\in W^{1}[0,1] & \;\mbox{and}\; & \xi_{1}(0)=\xi_{1}(1) \\
     \xi_{2}\in W^{1}[1,2] & \;\mbox{and}\; & \xi_{2}(1)=\xi_{2}(2) \\
   \end{array}\bigg\}\,. \\
\end{array}
$$

\smallskip

To compare the unbounded Fredholm modules $p_{*}[D]$ and $j_{*}\big([D],[D]\big)$\,, we first have
to compare the corresponding Hilbert spaces. Consider the ``doubling'' unitary
$$
U\colon
\left\{\begin{array}{rcl}
    L^2(S^{1})\oplus L^2(S^{1}) & \stackrel{\cong}{\longrightarrow} & L^2(S^{1})  \\
    (\xi_{1},\xi_{2}) & \longmapsto & \left(\;\theta\mapsto
    \arraycolsep1pt
    \left\{\begin{array}{ll}
        \!\sqrt{2}\!\cdot\!\xi_{1}(2\theta)\,,\; & \mbox{if}\;\,0\leq\theta\leq\frac
        {1}{2}  \\
        \!\sqrt{2}\!\cdot\!\xi_{2}(2\theta-1)\,,\; & \mbox{if}\;\,\frac{1}{2}<\theta\leq 1
    \end{array}\right.\,,\;\mbox{a.e.}\;\right)
\end{array}\right.
$$
The inverse $U^{*}$ of $U$ is given by the formula
$$
{\textstyle (U^{*}\xi)(\theta)=\Big(\frac{1}{\sqrt{2}}\!\cdot\!\xi\big(\frac{\theta}
{2}\big),\frac{1}{\sqrt{2}}\!\cdot\!\xi\big(\frac{\theta+1}{2}\big)\Big)\,,}
$$
for $\xi\in L^2(S^{1})$ and $\theta\in[0,1]$\,. Clearly, $U(\mathcal{M}_{1}\oplus\mathcal{M}_{2})U^{*}=
\mathcal{M}'$ holds, so that we get
$$
j_{*}\big([D],[D]\big)=\big[L^2(S^{1}),\mathcal{M}',U(D\oplus D)U^{*}\big]\,.
$$
It remains to discuss the relationship between the operators $D$ and $U(D\oplus D)U^{*}$\,. The
domain of the latter will be denoted by $\mathcal{E}$ and is given by
$$
\renewcommand{\arraystretch}{1.2}
\arraycolsep1pt
\begin{array}{rcl}
  \mathcal{E} & := & \dom\big(U(D\oplus D)U^{*}\big)=U\big(\dom(D\oplus D)\big)=U\big(W^{1}(S^{1})\oplus
  W^{1}(S^{1})\big)=\qquad\quad \\[.5em]
   & & \qquad\qquad\qquad\quad= \bigg\{(\xi_{1},\xi_{2})\in L^{2}[0,1]\,\bigg|
   \begin{array}{rcl}
     \xi_{1}\in W^{1}[0,\frac{1}{2}] & \;\mbox{and}\; & \xi_{1}(0)=\xi_
     {1}(\frac{1}{2}) \\
     \xi_{2}\in W^{1}[\frac{1}{2},1] & \;\mbox{and}\; & \xi_{2}(\frac{1}
     {2})=\xi_{2}(1) \\
   \end{array}\bigg\}\,, \\
\end{array}
$$
where we identify $L^{2}[0,1]$ with $L^{2}[0,\frac{1}{2}]\oplus
L^{2}[\frac{1}{2},1]$ in the obvious way. It is readily checked that
$U(D\oplus D)U^{*}$ equals $\frac{1}{2i}\!\cdot\!\frac{d}{d\theta}$ on this domain.
It is well-known that an unbounded Fredholm module $[\mathcal{H},\pi,
F]$ is equal to $[\mathcal{H},\pi,\lambda\cdot F]$ for every positive
real number $\lambda>0$ (the triples are operator-homotopic). So, we
have to show that
$$
\big[L^{2}(S^{1}),\mathcal{M}',\underbrace{{\textstyle \frac
{1}{i}\!\cdot\!\frac{d}{d\theta}}\;\,\mbox{on}\;\,W^{1}(S^{1})}_{=D}\big]=\big[L^{2}
(S^{1}),\mathcal{M}',\underbrace{{\textstyle\frac{1}{i}\!\cdot\!\frac{d}{d\theta}}\;\,
\mbox{on}\;\,\mathcal{E}}_{\!\!\!=2U(D\oplus D)U^{*}\!\!\!}\big]
$$
in $\KK_{1}(C(S^{1}\vee S^{1}),\bbC)$\,. The difficulty we alluded to on several occasions
is that the unitary $U$ does \emph{not} map the domain $W^{1}(S^{1})\oplus W^{1}(S^{1})$
to $W^{1}(S^{1})$\,, \ie $\mathcal{E}$ and $W^{1}(S^{1})$ are different. For this reason,
we define a new (dense) domain $\mathcal{D}$\,, contained in $W^{1}(S^{1})$ and in $\mathcal{E}$\,,
and more adapted to the situation, namely
$$
{\textstyle \mathcal{D}:=\left\{\xi\in W^{1}(S^{1})\,\left|\,\xi(0)=
\xi\left(\frac{1}{2}\right)=\xi(1)=0\right.\right\}\,.}
$$
Then $U$ maps $\mathcal{D}\oplus\mathcal{D}$ isometrically onto $\mathcal{D}$\,,
and we have a commutative diagram
\begin{displaymath}
 \begin{diagram}[height=2em,width=3.5em]
    \mathcal{D}\oplus\mathcal{D} & \rTo^{U}_{\cong} & \mathcal{D} \\
    \dTo^{2(D\oplus D)|_{\mathcal{D}\oplus\mathcal{D}}} & & \dTo_{D|_{\mathcal{D}}} \\
    L^{2}(S^{1})\oplus L^{2}(S^{1}) & \rTo^{\;\;\,U}_{\;\;\,\cong} & L^{2}(S^{1}) \\
 \end{diagram}
\end{displaymath}
As already singled out, the operator $D$ is self-adjoint on its domain $W^{1}
(S^{1})$\,, therefore, $2U(D\oplus D)U^{*}$ is also self-adjoint. The operator
$T:=D|_{\mathcal{D}}$ is closed, symmetric, but \emph{not} self-adjoint. In fact,
its adjoint $T^{*}$ is determined as in \cite[Ex.\ in~VIII.2, pp.\ 257--259]{RSI}\,:
it is the operator $T^{*}=\frac{1}{i}\!\cdot\!\frac{d}{d\theta}$ on the domain $\dom(T^{*})=W^{1}[0,1]$\,.
So, we are faced with two genuinely distinct self-adjoint extensions of $T$\,, namely,
one is $D$ with domain $W^{1}(S^{1})$\,, the other one is $2U(D\oplus D)U^{*}$ with domain
$\mathcal{E}$\,. Therefore, we are urged to try to apply Theorem~\ref{thm-main2} to show
that these operators define the same class in analytic $K$-homology. Before we proceed,
as a side-remark, we mention that the deficiency indices of $T$ are both equal to $2$\,,
and we refer to Example~\ref{ex} below for more details on this.

\smallskip

First, $T$ is a closed self-adjoint operator and its domain, $\mathcal{D}$\,, is dense, so that
condition~(a) of Theorem~\ref{thm-main2} is fulfilled. Secondly, to get condition~(b$^{\prime}$),
let us determine a $*$-closed dense subspace $\mathcal{B}'$ of $C(S^{1}\vee S^{1})$ such
that the operator $[\mathcal{M}'(f),T]$ is densely defined and bounded for every $f\in
\mathcal{B}'$\,. To do so, we identify $C(S^{1}\vee S^{1})$ with the $C^{*}$-algebra
$\big\{f\in C(S^{1})\,\big|\,f(0)=f(\frac{1}{2})\big\}$ in the obvious way, namely, still
viewing $S^{1}$ as $[0,1]\big/\partial[0,1]$\,. We correspondingly take for $\mathcal{B}'$
the $*$-closed dense sub-algebra $\big\{f\in C^{\infty}(S^{1})\,\big|\,f(0)=f(\frac{1}{2})
\big\}$\,, in other words,
$$
\mathcal{B}':=\big\{f\in C(S^{1}\vee S^{1})\,\big|\,f\circ p\in C^{\infty}(S^{1})\big\}\,.
$$
Observing that the subspace $\mathcal{M}'(\mathcal{B}')(\mathcal{D})$ is equal to $\mathcal{D}$\,,
the domain of $T$\,, we deduce that $\dom\big([\mathcal{M}'(f),T]\big)=\dom(T)$\,. For
$\xi$ in this domain and $f\in\mathcal{B}'$\,, we have
$$
{\textstyle [\mathcal{M}'(f),T]\xi=\frac{1}{i}\!\cdot\!\left((f\circ p)\!\cdot\!\frac{d\xi}
{d\theta}-\frac{d((f\circ p)\cdot\xi)}{d\theta}\right)=\frac{1}{i}\!\cdot\!\mathcal{M}
\big(\frac{d(f\circ p)}{d\theta}\big)\xi}\,.
$$
As hoped, Theorem~\ref{thm-main2} applies and yields the desired equality of analytic $K$-homology
classes defined by the two given self-adjoint extensions of $T$\,. This completes the proof.
\qed

\smallskip

We thank G.~Skandalis for pointing out to us the r\^{o}le of distinct self-adjoint extensions, while
we were trying to prove Theorem~\ref{Main1} below in an analytical way, which in fact amounts to the
present proof as we will see.

\begin{Rem}
The above proof shows that $(L^{2}(S^{1}),\mathcal{M}',T)$ is a symmetric unbounded Fredholm module,
with $T$ \emph{non}-self-adjoint, and defining a non-trivial analytic $K$-homology class $[L^{2}(S^{1}),
\mathcal{M}',T]\in\KK_{1}(C(S^{1}\vee S^{1}),\bbC)$\,.
\end{Rem}

As a matter of illustration, we would now like to give some more information
on the self-adjoint extensions of the operator $T$ of the preceding proof.

\begin{Ex}
\label{ex}
We keep notation as above. Obviously, one has
$$
{\textstyle \Ker(T^{*}{-i})=\bigg\{(\xi_{1},\xi_{2})\in W^{1}[0,\frac{1}{2}]\oplus
W^{1}[\frac{1}{2},1]\,\bigg|\!\!
\begin{array}{c}
  \xi_{1}(\theta)=\lambda_{1}\cdot e^{-\theta} \\
  \xi_{2}(\theta)=\lambda_{2}\cdot e^{-\theta} \\
\end{array},\;\lambda_{1},\lambda_{2}\in\bbC\bigg\}}
$$
and
$$
{\textstyle \Ker(T^{*}{+i})=\bigg\{(\xi_{1},\xi_{2})\in W^{1}[0,\frac{1}{2}]\oplus
W^{1}[\frac{1}{2},1]\,\bigg|\!\!
\begin{array}{c}
  \xi_{1}(\theta)=\lambda_{1}\cdot e^{\theta} \\
  \xi_{2}(\theta)=\lambda_{2}\cdot e^{\theta} \\
\end{array},\;\lambda_{1},\lambda_{2}\in\bbC\bigg\}}\,.
$$
So, in our situation, we have a `canonical' orthonormal basis for each deficiency space,
namely $\{(e'_{1},e'_{2}),(e''_{1},e''_{2})\}$ for the former and $\{(\varepsilon'_
{1},\varepsilon'_{2}),(\varepsilon''_{1},\varepsilon''_{2})\}$ for the latter,
where $e'_{2}$\,, $e''_{1}$\,, $\varepsilon'_{2}$ and $\varepsilon''_{1}$ are zero
functions, and, letting $\omega_{0}:=\sqrt{\frac{2}{e-1}}$\,,
$$
e'_{1}(\theta):=\omega_{0}\sqrt{e}\cdot e^{-\theta}\,,\;e''_{2}(\theta):=\omega_{0}e\cdot
e^{-\theta}\,,\;\varepsilon'_{1}(\theta):=\omega_{0}\cdot e^{\theta}\;\,\mbox{and}\;\,
\varepsilon''_{2}(\theta):=\frac{\omega_{0}}{\sqrt{e}}\cdot e^{\theta}\,.
$$
This gives an explicit homeomorphism between $\mathcal{U}$ and $U(2)$\,. By direct
computation, one checks that for every $u\in\mathcal{U}$\,, the self-adjoint extension
$T_{u}$ of $T$ is equal to $\frac{1}{i}\!\cdot\!\frac{d}{d\theta}$ on its domain $\dom(T_{u})$\,.
This is no surprise, since every $T_{u}$ is a restriction of $T^{*}$\,, which is also
$\frac{1}{i}\!\cdot\!\frac{d}{d\theta}$ on its domain $W^{1}[0,1]$\,, as we have seen.
One can wonder to which matrix in $U(2)$ does the operator $D$ (resp.\ $2U^{*}(D\oplus D)U$)
correspond to. One obtains
$$
D\;\longleftrightarrow\;
\big(\begin{smallmatrix}
  0 & 1 \\
  1 & 0 \\
\end{smallmatrix}\big)\in U(2)
\qquad\mbox{and}\qquad 2U^{*}(D\oplus D)U\;\longleftrightarrow\;
\big(\begin{smallmatrix}
  1 & 0 \\
  0 & 1 \\
\end{smallmatrix}\big)\in U(2)\,.
$$
In general, to a matrix
$
\big(\begin{smallmatrix}
  \alpha & \beta \\
  \gamma & \delta \\
\end{smallmatrix}\big)\in U(2)
$
with determinant $\Delta$\,, corresponds the operator $\frac{1}{i}\!\cdot\!\frac{d}{d\theta}$
on the domain consisting of the functions $(\xi_{1},\xi_{2})\in W^{1}[0,\frac{1}{2}]\oplus W^{1}
[\frac{1}{2},1]$ satisfying the boundary conditions
$$
{\textstyle
\Big(\!\begin{smallmatrix}
  \xi_{1}(0) \\[.4em]
  \xi_{2}(\frac{1}{2}) \\
\end{smallmatrix}\!\Big)=
\frac{1}{1+(\alpha+\delta)\sqrt{e}+\Delta e}
\Big(\begin{smallmatrix}
  \alpha+(1+\Delta)\sqrt{e}+\delta e & \beta(1-e) \\[.4em]
  \gamma(1-e) & \delta+(1+\Delta)\sqrt{e}+\alpha e \\
\end{smallmatrix}\Big)\cdot
\Big(\!\begin{smallmatrix}
  \xi_{1}(\frac{1}{2}) \\[.4em]
  \xi_{2}(1) \\
\end{smallmatrix}\!\Big)\,.}
$$
\end{Ex}


\section{Analytical proof of the Dolbeault-type ``addition formula''}

\label{s-mfld-proofs-2}


As in the case of the Dirac-type ``addition formula'', this section contains two
subsections. In the first one, we depict the analytic $K$-homology class determined
by the Dolbeault operator for $\Sigma_{g}$\,. In the second, we present the analytical
proof of Theorem~\ref{thm-add-formula-2}; again, this is one of the central parts of
the paper.


\subsection{Class of the Dolbeault operator for $\Sigma_{g}$ in analytic $K$-homology}

\label{ss-class-D-2}


\mbox{}\\
\vspace*{-.8em}

We fix an auxiliary K\"{a}hler structure on $\Sigma_
{g}$\,, \ie we view $\Sigma_{g}$ as a complex curve equipped with a suitably
compatible Riemannian metric. We let $\bar{\partial}_{\Sigma_{g}}:=\bar
{\partial}\oplus\bar{\partial}^{*}$ be the Dolbeault operator, \ie
\begin{displaymath}
 \begin{diagram}[height=2em,width=3.5em]
    L^{2}(\Lambda^{0,0}T^{*}\Sigma_{g})\oplus L^{2}(\Lambda^{0,2}T^{*}\Sigma_{g})
    \supset\dom(\bar{\partial}_{\Sigma_{g}}) & \rTo^{\bar{\partial}\oplus\bar
    {\partial}^{*}} & L^{2}(\Lambda^{0,1}T^{*}\Sigma_{g})\,. \\
 \end{diagram}
\end{displaymath}
Here, $L^{2}(\Lambda^{0,j}T^{*}\Sigma_{g})$ is the Hilbert space of $L^{2}$-forms
of bidegree $(0,j)$ on $T^{*}\Sigma_{g}$ (see for instance \cite[pp.~73--74]{Fried}).
In other words, we view $\Sigma_{g}$ as a K\"{a}hler manifold equipped with the `anti-canonical'
${\rm Spin}^{\!c}$-structure, and $\bar{\partial}_{\Sigma_{g}}$ is $\frac{1}{\sqrt{2}}
D_{\!_{\nabla}}$\,, where $D_{\!_{\nabla}}$ is the Dirac operator corresponding to the
Levi-Civita connection $\nabla$ (for the details, see \cite[pp. 77--81]{Fried}). The
domain of $\bar{\partial}_{\Sigma_{g}}$ is $W^{1}(\Lambda^{0,ev}T^{*}\Sigma_{g})$\,, a
Sobolev space on which the operator $\bar{\partial}_{\Sigma_{g}}$ is self-adjoint (see
\cite[pp.~100--101]{Fried} or \cite[Chap.~20]{BoossW}). To simplify the notation,
we let $[\bar{\partial}_{g}]$ denote the class of the operator $\bar{\partial}_{\Sigma_{g}}$
in $\KK_{0}(C(\Sigma_{g}),\bbC)\cong K_{0}(\Sigma_{g})$\,. Explicitly, $[\bar{\partial}_{g}]$
is given by the unbounded Fredholm module
$$
[\bar{\partial}_{g}]:=\big[L^{2}(\Lambda^{0,*}T^{*}\Sigma_{g}),\mathcal{M},\bar{\partial}_
{\Sigma_{g}}\big]\;\in\;\KK_{0}(C(\Sigma_{g}),\bbC)\,,
$$
where $L^{2}(\Lambda^{0,*}T^{*}\Sigma_{g})=\bigoplus_{j=0}^{2}L^{2}(\Lambda^{0,j}T^{*}
\Sigma_{g})$ is $\bbZ/2$-graded by even and odd degree forms, and $\mathcal{M}$ is the
$*$-representation of $C(\Sigma_{g})$ on $L^{2}(\Lambda^{0,*}T^{*}\Sigma_{g})$ given by pointwise
multiplication. By connectedness of the Teichm\"{u}ller space, the class $[\bar{\partial}_
{g}]$ is independent of the choice of the K\"{a}hler structure. For a later application,
let $\Lip(\Sigma_{g})$ be the $*$-closed dense sub-algebra of $C(\Sigma_{g})$ consisting
of the Lipschitz functions; by Rademacher's Theorem (see \cite[Thm.~11\,A]{Whi}), Lipschitz
functions are differentiable almost everywhere on $\Sigma_{g}$ and we single out that
$[\mathcal{M}(\vartheta),\bar{\partial}_{\Sigma_{g}}]$ is densely defined and bounded for
every function $\vartheta\in\Lip(\Sigma_{g})$\,.


\subsection{Analytical proof of Theorem \ref{thm-add-formula-2}}
\label{ss3.6}


\mbox{}\\
\vspace*{-.8em}

For the analytical proof, we will need the notions and notation introduced
in Subsection~\ref{ss2.6}, and we will apply Theorem~\ref{thm-main2}.

\smallskip

\AnPrfOf{Theorem~\ref{thm-add-formula-2}}
Again, using (the slight variation of) our general principle~\ref{lem-gen-princ-var}, still
with the same notation, it remains prove that
$$
j_{*}\big([\bar{\partial}_{g_{1}}],[\bar{\partial}_{g_{2}}]\big)=p_{*}[\bar{\partial}_{g_{1}+g_{2}-1}]
\;\in\;\KK_{0}(C(\Sigma_{g_{1}}\cup_{S^{1}}\Sigma_{g_{2}}),\bbC)\,.
$$
As in the one-dimensional case, this is exactly Proposition~\ref{prop-add-form-2}, and, this time,
we will establish it while sticking to analysis. We start by carefully describing the two Fredholm
modules under consideration. For sake of readability, we set $g:=g_{1}+g_{2}-1$ and $X:=\Sigma_{g_{1}}
\cup_{S^{1}}\Sigma_{g_{2}}$\,. First, in the group $\KK_{0}(C(X),\bbC)$\,, we have
$$
p_{*}\big[\bar{\partial}_{g}]=\big[L^{2}(\Lambda^{0,*}T^{*}\Sigma_{g}),\mathcal{M}',\bar{\partial}_{\Sigma_
{g}}\big]\,,
$$
where, for a function $\vartheta\in C(X)$\,, the operator $\mathcal{M}'(\vartheta)$ is fiber-wise multiplication
by $\vartheta\circ p\in C(\Sigma_{g})$ on the Hilbert space $L^{2}(\Lambda^{0,*}T^{*}\Sigma_{g})$ of $L^{2}$-sections
of the vector bundle $\Lambda^{0,*}T^{*}\Sigma_{g}$\,; as before, the domain of $\bar{\partial}_{\Sigma_{g}}$ is
$W^{1}(\Lambda^{0,ev}T^{*}\Sigma_{g})$\,. On the other hand, we get
$$
j_{*}\big([\bar{\partial}_{g_{1}}],[\bar{\partial}_{g_{2}}]\big)=\big[L^{2}(\Lambda^{0,*}T^{*}\Sigma_{g_{1}})
\oplus L^{2}(\Lambda^{0,*}T^{*}\Sigma_{g_{2}}),\mathcal{M}_{1}\oplus\mathcal{M}_{2},\bar{\partial}_{\Sigma_
{g_{1}}}\oplus\bar{\partial}_{\Sigma_{g_{2}}}\big]\,,
$$
where the direct sum is an orthogonal and graded one, and, for $\vartheta\in C(X)$ and $i=1,2$\,, we have
$\mathcal{M}_{i}(\vartheta)=\mathcal{M}(\vartheta_{i})$ with $\vartheta_{i}$ standing for the restriction
$\vartheta|_{\Sigma_{g_{i}}}$\,; the domain of $\bar{\partial}_{\Sigma_{g_{1}}}\oplus\bar{\partial}_{\Sigma_
{g_{2}}}$ is the orthogonal direct sum
$$
\dom(\bar{\partial}_{\Sigma_{g_{1}}}\oplus\bar{\partial}_{\Sigma_{g_{2}}})=W^{1}(\Lambda^{0,ev}T^{*}\Sigma_{g_{1}})
\oplus W^{1}(\Lambda^{0,ev}T^{*}\Sigma_{g_{2}})\,.
$$
Now, we would like to determine a grading-preserving unitary isomorphism
$$
U\colon L^{2}(\Lambda^{0,*}T^{*}\Sigma_{g_{1}})\oplus L^{2}(\Lambda^{0,*}T^{*}\Sigma_{g_{2}})\stackrel{\cong}
{\longrightarrow}L^{2}(\Lambda^{0,*}T^{*}\Sigma_{g})\,.
$$
We can modify $\Sigma_{g_{1}}\natural\Sigma_{g_{2}}$ by an orientation-preserving analytic diffeomorphism,
so, we may suppose that the modified connected sum $\Sigma_{g_{1}}\natural\Sigma_{g_{2}}\cong\Sigma_{g}$
is obtained from $\Sigma_{g_{1}}$ and $\Sigma_{g_{2}}$ by gluing the open manifolds $V_{1}:=\Sigma_{g_{1}}\!
\smallsetminus\!C_{1}$ and $V_{2}:=\Sigma_{g_{2}}\!\smallsetminus\!C_{2}$ along the closed manifold $K:=S^{1}
\amalg S^{1}$\,, with $C_{1}$ and $C_{2}$ as in Lemma~\ref{lem-gen-princ-var}, \ie
$$
\Sigma_{g_{1}}\natural\Sigma_{g_{2}}=(\underbrace{\Sigma_{g_{1}}\!\smallsetminus\!C_{1}}_{=V_{1}})\amalg(\underbrace
{S^{1}\amalg S^{1}}_{=K})\amalg(\underbrace{\Sigma_{g_{2}}\!\smallsetminus\!C_{2}}_{=V_{2}})\,.
$$
This way, we can consider $V_{i}$ ($i=1,2$) as an analytic open sub-manifold of both $\Sigma_{g_{i}}$ and $\Sigma_
{g_{1}}\natural\Sigma_{g_{2}}$\,, and the complement in the latter of the union $V_{1}\amalg V_{2}$ is of measure
zero. Moreover, $p$ merely identifies the two copies of $S^{1}$ pointwise. Now, with this in mind, we define $U$
almost everywhere by the formula
$$
U(\omega_{1},\omega_{2}):=
\arraycolsep1pt
\renewcommand{\arraystretch}{1.2}
\left\{\begin{array}{ll}
    \sqrt{2}\!\cdot\!\omega_{1}|_{V_{1}}\,,\; & \mbox{on}\;\,V_{1}\subset\Sigma_{g}  \\
    \sqrt{2}\!\cdot\!\omega_{2}|_{V_{2}}\,,\; & \mbox{on}\;\,V_{2}\subset\Sigma_{g}\,.
\end{array}\right.
$$
The inverse is simply given (almost everywhere) by
$$
{\textstyle U^{*}\,\omega:=\Big(\frac{1}{\sqrt{2}}\!\cdot\!\omega|_{V_{1}},\frac{1}{\sqrt{2}}\!\cdot\!
\omega|_{V_{2}}\Big)\,.}
$$
Is it obvious that $U$ intertwines $\mathcal{M}_{1}\oplus\mathcal{M}_{2}$ and $\mathcal{M}'$\,, \ie
$U(\mathcal{M}_{1}\oplus\mathcal{M}_{2})U^{*}=\mathcal{M}'$\,. It follows that
$$
j_{*}\big([\bar{\partial}_{g_{1}}],[\bar{\partial}_{g_{2}}]\big)=\big[L^{2}(\Lambda^{0,*}T^{*}\Sigma_{g}),
\mathcal{M}',U(\bar{\partial}_{\Sigma_{g_{1}}}\oplus\bar{\partial}_{\Sigma_{g_{2}}})U^{*}\big]\,,
$$
where the operator appearing has domain
$$
\dom\big(U(\bar{\partial}_{\Sigma_{g_{1}}}\oplus\bar{\partial}_{\Sigma_{g_{2}}})U^{*}\big)=
U\big(W^{1}(\Lambda^{0,ev}T^{*}\Sigma_{g_{1}})\oplus W^{1}(\Lambda^{0,ev}T^{*}\Sigma_{g_{2}})\big)\,.
$$
To see what happens at the level of the domains of the unbounded operators involved, we first define a
dense subspace $\widetilde{\mathcal{D}}$ of $L^{2}(\Lambda^{0,ev}T^{*}\Sigma_{g})$ by
$$
\mathcal{D}:=\big\{\omega\in W^{1}(\Lambda^{0,ev}T^{*}\Sigma_{g})\,\big|\;\omega|_{K}=0\big\}=H^{1}_{0}
(\Lambda^{0,ev}T^{*}V_{1})\oplus H^{1}_{0}(\Lambda^{0,ev}T^{*}V_{2})\,.
$$
For the definition of the Sobolev space $H^{1}_{0}(\Lambda^{0,ev}T^{*}V_{i})$\,, for the latter equality and
for the sense to give to the equation $\omega|_{K}=0$\,, we refer to \cite{TayI}, p.~290, Ex.~4.5.2 on
p.~294, and to Prop.~4.4.5 on p.~287, respectively. Similarly, for $i=1,2$\,, we define a dense subspace
$\mathcal{D}_{i}$ in $L^{2}(\Lambda^{0,ev}T^{*}\Sigma_{g_{i}})$ by
$$
\mathcal{D}_{i}:=\big\{\omega_{i}\in W^{1}(\Lambda^{0,ev}T^{*}\Sigma_{g_{i}})\,\big|\;\omega_{i}|_{C_{i}}=0\big\}=
H^{1}_{0}(\Lambda^{0,ev}T^{*}V_{i})\,.
$$
The point is that $U$ maps $\mathcal{D}_{1}\oplus\mathcal{D}_{2}$ isometrically onto $\mathcal{D}$\,,
and, since $V_{1}$ and $V_{2}$ are analytic open sub-manifolds of $\Sigma_{g}$ and since Dolbeault
operators are local (\ie defined locally), there is a commutative diagram
\begin{displaymath}
 \begin{diagram}[height=2em,width=6em]
   \mathcal{D}_{1}\oplus\mathcal{D}_{2} & \rTo^{\!\!\!\!\!\!U}_{\!\!\!\!\!\!\cong} & \mathcal{D} \\
   \dTo^{\bar{\partial}_{\Sigma_{g_{1}}}\oplus\bar{\partial}_{\Sigma_{g_{2}}}|_{\mathcal{D}_{1}\oplus\mathcal{D}_{2}}}
   & & \dTo_{\bar{\partial}_{\Sigma_{g}}|_{\mathcal{D}}} \\
   L^{2}(\Lambda^{0,1}T^{*}\Sigma_{g_{1}})\oplus L^{2}(\Lambda^{0,1}T^{*}\Sigma_{g_{2}}) &\rTo^{\;\;\;\;\;U}_
   {\;\;\;\;\;\cong} & L^{2}(\Lambda^{0,1}T^{*}\Sigma_{g}) \\
 \end{diagram}
\end{displaymath}

So, letting $T:=\bar{\partial}_{\Sigma_{g}}|_{\mathcal{D}}$\,, we are faced with
two self-adjoint extensions of the densely defined symmetric operator $T$\,, namely
$\bar{\partial}_{\Sigma_{g}}$ with domain $W^{1}(\Lambda^{0,ev}T^{*}\Sigma_{g})$ and
$U(\bar{\partial}_{\Sigma_{g_{1}}}\oplus\bar{\partial}_{\Sigma_{g_{2}}})U^{*}$ with domain
$U\big(W^{1}(\Lambda^{0,ev}T^{*}\Sigma_{g_{1}})\oplus W^{1}(\Lambda^{0,ev}T^{*}\Sigma_{g_{2}})
\big)$\,, and we have to show that they define the same $K$-homology class. Again, the instructive
difficulty is that these two domains are distinct. As in the one-dimensional case, we will
now verify that Theorem~\ref{thm-main2} applies to establish the desired $K$-equality.
Condition~(a) of Theorem~\ref{thm-main2} being clearly fulfilled by $T$\,, to get
condition~(b$^{\prime}$), let us determine a $*$-closed dense subspace $\mathcal{B}'$
of $C(X)$ such that the operator $[\mathcal{M}'(\vartheta),T]$ is densely defined and
bounded for every $\vartheta\in\mathcal{B}'$\,. Let us consider the $*$-closed dense
subalgebra $\Lip(X)$ of $C(X)$ consisting of the Lipschitz functions on $X=\Sigma_{g_{1}}
\cup_{S^{1}}\Sigma_{g_{2}}$\,. Given a function $\vartheta\in\Lip(X)$\,, the map $p\colon
\Sigma_{g}\longrightarrow X$ being Lipschitz, we see that the composite satisfies
$\vartheta\circ p\in\Lip(\Sigma_{g})$\,. By the final sentence in Subsection~\ref{ss-class-D-2},
the operator $[\mathcal{M}'(\vartheta),T]$ is indeed densely defined and bounded for
every $\vartheta\in\Lip(X)$\,. Finally, Theorem~\ref{thm-main2} applies and gives the
desired equality of analytic $K$-homology classes defined by the two self-adjoint extensions
of $T$ at hand. This completes the proof.
\qed

\begin{Rem}
This proof shows that the triple $\big(L^{2}(\Lambda^{0,*}T^{*}\Sigma_{g_{1}+g_{2}-1}),\mathcal{M}',T\big)$ is
a symmetric unbounded Fredholm module, with $T$ \emph{non}-self-adjoint, and defining a non-trivial
analytic $K$-homology class $\big[L^{2}(\Lambda^{0,*}T^{*}\Sigma_{g_{1}+g_{2}-1}),\mathcal{M}',T]$
in the group $\KK_{0}\big(C(\Sigma_{g_{1}}\cup_{S^{1}}\Sigma_{g_{2}}),\bbC\big)$\,.
\end{Rem}

\begin{Rem}
The deficiency indices of the symmetric unbounded operator $T$ in the above proof are equal and
countably infinite.
\end{Rem}

\begin{Rem}
Contrarily to the Dirac case (see the commutative diagram in the analytical proof of
Theorem~\ref{thm-add-formula-1} in Subsection~\ref{ss2.6}), in the commutative diagram
with Dolbeault operators in the proof above, there is no constant popping up like the $2$
in the Dirac case. The reason for this is the equality
$$
\underbrace{\Area(\Sigma_{g_{1}}\cup_{S^{1}}\Sigma_{g_{2}})}_{\!\!=\Area(\Sigma_{g_{1}})+
\Area(\Sigma_{g_{2}})\!\!}=\Area(\Sigma_{g})\,,
$$
of areas, whereas, in the Dirac case, with our choice of parametrizations, we have
$$
\underbrace{\Length(S^{1}\vee S^{1})}_{\!\!\!\!\!\!\!\!=\Length(S^{1})+\Length(S^{1})\!\!\!\!\!\!\!\!}=
2\cdot\Length(S^{1})\,.
$$
\end{Rem}


\part{Application to the Baum-Connes Conjecture}

\label{p-IV}


\section{The first homology of a group and the Baum-Connes map}

\label{s-j1}


%
In this section, subdivided into five subsections, we treat our program described
in Section~\ref{s-appl-BC} for the case $j=1$\,.


\subsection{Topological and analytical definitions of $\betaa{1}$}

\label{ss2.1}


\mbox{}\\
\vspace*{-.8em}

We denote the abelianization of the group $\Gamma$ by $\Gamma^{\ab}$ and we identify
it with $H_{1}(\Gamma;\bbZ)$ in the usual way. We write $\gamma^{\ab}$ for the class
of the element $\gamma\in\Gamma$ in the quotient group $\Gamma^{\ab}$\,. We consider
$$
\betaa{1}\colon\Gamma^{\ab}\longrightarrow K_{1}(C^{*}_{r}\Gamma)\,,\quad\gamma^{\ab}\longmapsto
[\gamma]=\big[\diag(\gamma,1,1,\ldots)\big]\,,
$$
the canonical homomorphism induced by the homomorphism $\tbetaa{1}\colon\Gamma
\longrightarrow K_{1}(C^{*}_{r}\Gamma)$ coming from the inclusion of $\Gamma$
into the group of invertible elements in $C^{*}_{r}\Gamma$\,. In the analytical
description of $K\!K$-theory, the class $[\gamma]\in K_{1}(C^{*}_{r}\Gamma)\cong\KK_{1}
(\bbC,C^{*}_{r}\Gamma)$ is given via the equality
$$
[\gamma]=\alpha^{\gamma}_{*}\big[\mathcal{E},\textstyle{\frac{d}{dx}}\big]\;\in\;\KK_{1}
(\bbC,C^{*}_{r}\Gamma)\,,
$$
where the notation is as follows. First, $[\mathcal{E},\frac{d}{dx}]=[\mathcal{E},\pi_{\circ},
\frac{d}{dx}]\in\KK_{1}(\bbC,C^{*}_{r}\bbZ)\cong\bbZ$ is the `standard' generator, with
$\mathcal{E}$ denoting the separation-completion of the algebra $C_{c}^{\infty}(\bbR)$ of
compactly supported smooth complex-valued functions on the real line with respect to the
$C^{*}_{r}\bbZ$-valued scalar product determined by
$$
\left<\xi_{1}|\xi_{2}\right>(n):=\left<\xi_{1}|\varrho(n)\cdot\xi_{2}\right>_{L^{2}(\bbR)}
=\int_{\bbR}\overline{\xi_{1}(x)}e^{-2\pi i nx}\xi_{2}(x)\,dx\,,
$$
for $\xi_{1},\xi_{2}\in C_{c}^{\infty}(\bbR)$ and $n\in\bbZ$\,, where $\varrho$ is the
action of $\bbZ$ on $C_{c}^{\infty}(\bbR)$ by point-wise multiplication by integer powers
of the function $e^{-2\pi ix}$\,; and $\pi_{\circ}\colon\bbC\longrightarrow\mathcal{L}_{C^{*}_
{r}\bbZ}(\mathcal{E})$ is the unit, \ie the $*$-homomorphism taking $\lambda\in\bbC$ to
$\lambda\cdot\id_{\mathcal{E}}$\,; compare with \cite[Section~4.2 in Part~2]{Val}. Second,
$\alpha^{\gamma}_{*}$ stands for the map
$$
\alpha^{\gamma}_{*}\colon\KK_{1}(\bbC,C^{*}_{r}\bbZ)\longrightarrow\KK_{1}(\bbC,C^{*}_{r}\Gamma)
$$
induced by the composition of $*$-homomorphisms
$$
C^{*}_{r}\bbZ=C^{*}\bbZ\stackrel{\smash[t]{\widehat{\alpha^{\gamma}}}}{\longrightarrow}C^{*}\Gamma\stackrel
{\,\lambda_{\Gamma}}{\onto}C^{*}_{r}\Gamma
$$
defined using \emph{maximal} group-$C^{*}$-algebras, where the first indicated $*$-homomorphism
is induced by the obvious group homomorphism determined by $\gamma$\,, namely
$$
\alpha^{\gamma}\colon\bbZ\longrightarrow\Gamma\,,\quad n\longmapsto\gamma^{n}\,;
$$
the map $\lambda_{\Gamma}$ is the canonical epimorphism. It is also possible to describe $[\gamma]$
directly as an unbounded Kasparov element (in the sense of Baaj-Julg~\cite{BJ}), namely,
$$
[\gamma]=[\mathcal{E}',\textstyle{\frac{d}{dx}}]=[\mathcal{E}',\pi_{\circ}',\textstyle{\frac{d}{dx}}]
\;\in\;\KK_{1}(\bbC,C^{*}_{r}\Gamma)\,,
$$
where $\mathcal{E}'$ is the separation-completion of $C_{c}^{\infty}(\bbR)$ with respect to the
$C^{*}_{r}\Gamma$-valued scalar product determined by
$$
\left<\xi_{1}|\xi_{2}\right>(\gamma'):=\sum_{\!\!\!n\in(\alpha^{\gamma})^{-1}(\gamma')\!\!\!}\left<
\xi_{1}|\varrho(n)\cdot\xi_{2}\right>_{L^{2}(\bbR)}\,,
$$
for $\xi_{1},\xi_{2}\in C_{c}^{\infty}(\bbR)$ and $\gamma'\in\Gamma$\,, and where $\pi_{\circ}'
\colon\bbC\longrightarrow\mathcal{L}_{C^{*}_{r}\Gamma}(\mathcal{E}')$ is the unit.


\subsection{Topological definition of $\betat{1}$}

\label{ss2.2}


\mbox{}\\
\vspace*{-.8em}

We begin by constructing a homomorphism $\betat{1}\colon\Gamma^{\ab}\longrightarrow
K_{1}(B\Gamma)$ in such a way that $\nu^{\Gamma}_{1}\circ\betat{1}=\betaa{1}$\,.
This was previously done by Natsume~\cite{Nat}, under the extra assumption that
$\Gamma$ is a torsion-free group. Our definition of $\betat{1}$ will be in two steps\,:
first, we define a (set-theoretic\,!) map $\tbetat{1}\colon\Gamma\longrightarrow
K_{1}(B\Gamma)$\,; next, we prove that $\tbetat{1}$ is a group homomorphism.
Since the target group is abelian, this will imply that $\tbetat{1}$ factors
through the desired homomorphism $\betat{1}$\,.

\smallskip

To define $\tbetat{1}$\,, we notice that since $\pi_{1}(B\Gamma)=\Gamma$\,,
every element $\gamma\in\Gamma$ defines (up to homotopy) a pointed continuous
map $\gamma\colon S^{1}\longrightarrow B\Gamma$\,. Keeping notation as in
Lemma~\ref{lem3}, we let
$$
[S^{1}]_{K}:=(ch_{\odd}^{\bbZ})^{-1}\big([S^{1}]\big)\;\in\;K_{1}(S^{1})
$$
be the (unique) $K$-homology class with integral Chern character given by
the fundamental class (the usual orientation, and even $\Spinc$-structure, is
fixed on $S^{1}$). Letting $D:=\frac{1}{i}\!\cdot\!\frac{d}{d\theta}$ be
the Dirac operator, see Section~\ref{s-mfld-proofs-1}, by Lemma~\ref{lem-ch-S1},
we have
$$
[S^{1}]_{K}={-[D]}=[-D]=[i\!\cdot\!{\textstyle\frac{d}{d\theta}}]\;\in\;K_{1}(S^{1})\,.
$$
By functoriality, we get a homomorphism $\gamma_{*}\colon K_{1}(S^{1})\longrightarrow
K_{1}(B\Gamma)$ and we set $\tbetat{1}(\gamma):=\gamma_{*}[S^{1}]_{K}$\,, for $\gamma\in
\Gamma$\,, so that
$$
\betat{1}\colon\Gamma^{\ab}\longrightarrow K_{1}(B\Gamma)\,,\quad\gamma^{\ab}
\longmapsto\gamma_{*}[S^{1}]_{K}={-\gamma_{*}[D]}\,.
$$


\subsection{Analytical definition of $\betat{1}$}

\label{ss2.3}


\mbox{}\\
\vspace*{-.8em}

We describe the map $\betat{1}\colon\Gamma^{\ab}\longrightarrow K_{1}(B\Gamma)$ analytically, using
the unbounded picture for $K$-homology, see~\cite{BJ}. The element $[S^{1}]_{K}$ is then described as the
unbounded Fredholm module (see Subsection~\ref{ss-class-D-1})
$$
[S^{1}]_{K}=[-D]=[L^2(S^{1}),\mathcal{M},-D]\;\in\;\KK_{1}(C(S^{1}),\bbC)\,,
$$
where $\mathcal{M}$ is the $*$-representation of $C(S^{1})$ on $L^2(S^{1})$ by pointwise multiplication.

\smallskip

If $\gamma\in\Gamma$ corresponds to a map $\gamma\colon S^{1}\longrightarrow
B\Gamma$\,, and if $X$ is an arbitrary compact subspace of $B\Gamma$ containing
$\gamma(S^{1})$\,, as for example $\gamma(S^{1})$ itself, then $\betat{1}
(\gamma^{\ab})$ is described by image of the Fredholm module
$$
\gamma_{*}[L^2(S^{1}),\mathcal{M},-D]\in\KK_{1}(C(X),\bbC)
$$
(where $\gamma$ is viewed as a map from $S^{1}$ to $X$) under the homomorphism
$$
\KK_{1}(C(X),\bbC)\longrightarrow\RKK_{1}(B\Gamma,\bbC)\,.
$$
induced by the inclusion (recall that $\RKK_{1}(B\Gamma,\bbC)$ is by definition the colimit,
over the compact subspaces $Y$ of $B\Gamma$\,, of the abelian groups $\KK_{1}(C(Y),\bbC)$).
Assume moreover that the map $\gamma\colon S^{1}\longrightarrow X$ is Lipschitz; up to homotopy,
one can always make this assumption on the map $\gamma$ (with $\gamma(S^{1})$ as suitable $X$).
Then, letting $\gamma^{*}\colon C(X)\longrightarrow C(S^{1})$ take a function $f$ to $f\circ\gamma$\,,
we see that the $*$-closed subalgebra $\Lip(X)$ of $C(X)$ is dense and $\gamma^{*}\Lip(X)$ verifies
$\gamma^{*}\Lip(X)\subseteq\Lip(S^{1})$ and consists therefore of functions that are differentiable
almost everywhere by Rademacher's Theorem (see \cite[Thm.~11\,A]{Whi}), so that
$$
\gamma_{*}[L^2(S^{1}),\mathcal{M},-D]=[L^2(S^{1}),\mathcal{M}\circ\gamma^{*},-D]\;\in\;\KK_{1}
(C(X),\bbC)\,.
$$


\subsection{Properties of $\betat{1}$}

\label{ss2.4}


%
\begin{Thm}
\label{Main1}
The map $\tbetat{1}\colon\Gamma\longrightarrow K_{1}(B\Gamma)$ is a
group homomorphism. Consequently, the map
$$
\betat{1}\colon\Gamma^{\ab}\longrightarrow K_{1}(B\Gamma)\,,\quad\gamma^{\ab}
\longmapsto\gamma_{*}[S^{1}]_{K}=-\gamma_{*}[D]
$$
is a well-defined group homomorphism.
\end{Thm}

This will be proved in Subsection~\ref{ss2.5} below. Before the proof, assuming
Theorem~\ref{Main1} for a while, we deduce some consequences.

\begin{Rem}
We claim that \emph{$\varphi_{1}^{\Gamma}\circ\tbetat{1}$ is zero on torsion elements
of $\Gamma$}\,, where, recall, $\varphi_{1}^{\Gamma}$ denotes the canonical map
$K_{1}(B\Gamma)\longrightarrow K_{1}^{\Gamma}(\bEG)$\,.  Indeed, if $\gamma\in\Gamma$
has order $n\geq 1$\,, the map $\gamma_{*}\colon K_{1}(S^{1})\longrightarrow K_{1}(B\Gamma)$
factorizes as
\begin{displaymath}
 \begin{diagram}[height=2em,width=3em]
    K_{1}(S^{1}) & & \rTo^{\gamma_{*}} & & K_{1}(B\Gamma) \\
     & \rdTo & & \ruTo_{B\!\incl_{*}} & \\
     & & K_{1}(B\bbZ/n) & & \\
 \end{diagram}
\end{displaymath}
where $\incl\colon\bbZ/n\hookrightarrow\Gamma$ takes $1$ to $\gamma$\,. On the
other hand, the diagram
\begin{displaymath}
 \begin{diagram}[height=2em,width=3.5em]
    K_{1}(B\bbZ/n) & \rTo^{\!\!\!B\!\incl_{*}} & K_{1}(B\Gamma) \\
    \dTo^{\varphi_{1}^{\bbZ/n}} & & \dTo_{\varphi_{1}^{\Gamma}} \\
    K_{1}^{\bbZ/n}(\bE\bbZ/n) & \rTo^{\!\underline{E}\!\incl_{*}} & K_{1}^{\Gamma}(\bEG) \\
 \end{diagram}
\end{displaymath}
commutes. However, one can take $\bE\bbZ/n=pt$\,, so that $K_{1}^{\bbZ/n}(\bE\bbZ/n)=0$\,.
Our claim follows. This observation is elaborated on in~\cite{MM2}.
\end{Rem}

\begin{Prop}
\label{prop2}
Let $ch_{\odd}\colon K_{1}(B\Gamma)\longrightarrow\bigoplus_{n=1}^{\infty}H_
{2n+1}(\Gamma;\bbQ)$ be the odd Chern character. Then $(ch_{\odd}\otimes\id_{\bbQ})
\circ(\betat{1}\otimes\id_{\bbQ})=\id_{H_{1}(\Gamma;\bbQ)}$ holds, in particular,
$\betat{1}$ is rationally injective.
\end{Prop}

\Prf
Fix an element $\gamma\in\Gamma$\,, and denote by $\alpha^{\!\gamma}\colon\bbZ
\longrightarrow\Gamma$ the homomorphism taking $1$ to $\gamma$\,. Note that
the pointed continuous map $\gamma\colon B\bbZ=S^{1}\longrightarrow B\Gamma$
considered earlier is merely $B\alpha^{\!\gamma}$\,. Due to the naturality of
the Chern character in $K$-homology, we have a commutative diagram
\begin{displaymath}
 \begin{diagram}[height=2em,width=4em]
    K_{1}(S^{1})\otimes\bbQ & \rTo^{\!\!\alpha^{\!\gamma}_{*}} & K_{1}
    (B\Gamma)\otimes\bbQ \\
    \dTo^{ch_{\odd}\otimes\id_{\bbQ}}_{\cong} & & \dTo_{ch_{\odd}\otimes\id_{\bbQ}}^{\cong} \\
    H_{1}(S^{1};\bbQ) & \rTo^{\!\!\!\!\alpha^{\!\gamma}_{*}} & H_{\odd}(B
    \Gamma;\bbQ) \\
 \end{diagram}
\end{displaymath}
Then, dropping ``\,$\otimes\id_{\bbQ}$\,'' and ``\,$\otimes 1$\,'' from the notation,
we compute
$$
ch_{\odd}\,\betat{1}(\gamma^{\ab})=ch_{\odd}\,\alpha^{\!\gamma}_{*}[S^{1}]_{K}=\alpha^
{\!\gamma}_{*}\,ch_{\odd}[S^{1}]_{K}=\alpha^{\!\gamma}_{*}[S^{1}]=\gamma^{\ab}\,,
$$
where we have used the fact that for $S^{1}$\,, the usual Chern character takes $[S^{1}]_{K}$
to the fundamental class $[S^{1}]$ in rational homology (see Lemma~\ref{lem-ch-S1}).
\qed

\begin{Thm}
\label{thm1}
The equality $\betaa{1}=\nu^{\Gamma}_{1}\circ\betat{1}$ holds.
\end{Thm}

\Prf
Clearly it is enough to prove that $\tbetaa{1}=\nu^{\Gamma}_{1}\,\tbetat{1}$\,.
As in the proof of Proposition~\ref{prop2}, we fix $\gamma\in\Gamma$ and write
$\alpha^{\!\gamma}\colon\bbZ\longrightarrow\Gamma$ for the corresponding homomorphism.
Consider the diagram
\begin{displaymath}
 \begin{diagram}[height=2.2em,width=3.5em]
    \bbZ & & & \rTo^{\alpha^{\!\gamma}} & & & \Gamma \\
     & \rdTo^{\betaa{1}} & & & & \ldTo^{\tbetaa{1}} & \\
    \dTo^{\betat{1}} & & K_{1}(C^{*}_{r}\bbZ) & \rTo^{\alpha^{\!\gamma}_{*}} &
    K_{1}(C^{*}_{r}\Gamma) & & \dTo_{\tbetat{1}} \\
     & \ruTo^{\nu_{1}^{\bbZ}} & & & & \luTo^{\nu_{1}^{\Gamma}} & \\
    K_{1}(S^{1}) & & & \rTo^{\alpha^{\!\gamma}_{*}} & & & K_{1}(B
    \Gamma) \\
 \end{diagram}
\end{displaymath}
We have $\alpha^{\!\gamma}_{*}\,\betaa{1}=\tbetaa{1}\alpha^{\!\gamma}$
by obvious reasons, $\tbetat{1}\alpha^{\!\gamma}=\alpha^{\!\gamma}_{*}\,
\betat{1}$ by definition of $\tbetat{1}$\,, and $\alpha^{\!\gamma}_{*}\,
\nu_{1}^{\bbZ}=\nu^{\Gamma}_{1}\,\alpha^{\!\gamma}_{*}$ by naturality of the
Novikov assembly map when the source group is $K$-amenable (see~\cite[Cor.~1.3 in Part~2]{Val}).
By diagram chasing, one sees that the desired equality $\tbetaa
{1}=\nu_{1}^{\Gamma}\,\tbetat{1}$ follows from the analogous result
for $\bbZ$\,, namely from $\betaa{1}=\nu^{\bbZ}_{1}\,\betat{1}$\,, which
in turn is a consequence of the well-known fact that the Baum-Connes
Conjecture holds for the group $\bbZ$ (see \cite[12.5.9]{HigRoe},
\cite[Section~4 in Part~2]{Val} or \cite[Ex.~6.1.6]{Val99} for a
direct proof).
\qed

\smallskip

We have already mentioned in Section~\ref{s-appl-BC} that $\betaa{1}\colon\Gamma^
{\ab}\longrightarrow K_{1}(C^{*}_{r}\Gamma)$ is rationally injective
(see~\cite{BV,ElN}).


\subsection{Proof of Theorem \ref{Main1}}

\label{ss2.5}


\mbox{}\\
\vspace*{-.8em}

We treat the topological and the analytical settings together.
Consider two elements $\gamma_{1},\gamma_{2}\in\Gamma$\,, viewed as (homotopy
classes of) pointed continuous maps $S^{1}\longrightarrow B\Gamma$\,.
By definition of $K$-homology with compact supports and of $\RKK$-groups, both
$K_{1}(B\Gamma)$ and $\RKK_{1}(B\Gamma,\bbC)$ are defined as the colimit of
$K_{1}(Y)$ and $\KK_{1}(C(Y),\bbC)$\,, respectively, with $Y$ running over the
compact subspaces of $B\Gamma$\,. Letting $X:=\gamma_{1}(S^{1})\cup\gamma_{2}(S^{1})$\,,
a compact subspace of $B\Gamma$\,, it is therefore enough to check that the equality
$$
(\gamma_{1}\gamma_{2})_{*}[S^{1}]_{K}=(\gamma_{1})_{*}[S^{1}]_{K}+(\gamma_{2})_{*}[S^{1}]_{K}
$$
holds in $K_{1}(X)$ and $\KK_{1}(C(X),\bbC)$ respectively, where $\gamma_{1}\gamma_{2}$
stands for the product-loop. Up to homotopy, we may assume that $\gamma_{1}$ and $\gamma_{2}$
are constant on a neighbourhood of the base-point of $S^{1}$\,. The key-point that allows
to connect the present situation with what has been done so far, is that the product-loop
is nothing but the composition of maps
$$
\gamma_{1}\gamma_{2}=\gamma_{1}\#\gamma_{2}\colon\underbrace{S^{1}\# S^{1}}_{=S^{1}}\stackrel{p}
{\longrightarrow}S^{1}\vee S^{1}\stackrel{\,\gamma_{1}\vee\gamma_{2}}{\longrightarrow}X\,,
$$
where we borrow the notation from Proposition~\ref{prop-add-form-1}, and where we identify
$S^{1}\# S^{1}$ with $S^{1}$\,, as indicated. Bearing in mind the equality $[S^{1}]_{K}={-[D]}$\,,
what has to be proved is that
$$
(\gamma_{1}\#\gamma_{2})_{*}[D]=(\gamma_{1})_{*}[D]+(\gamma_{2})_{*}[D]\,,
$$
which is precisely the ``addition formula'' for the Dirac operator of Theorem~\ref{thm-add-formula-1}.
This proves Theorem~\ref{Main1} both from the topological and from the analytical viewpoint
on $\betat{1}$\,.
\qed

\begin{Rem}
We have spent some time on the analytical proof, because it illustrates a difficulty that, apparently,
went unnoticed so far. A detailed and explicit treatment of this difficulty is in fact one of the central
themes in these notes. We also point out that the second named author provides in~\cite{MM2} an abstract
proof of Theorem~\ref{Main1}, which is of purely homotopical nature.
\end{Rem}


\section{The second homology of a group and the Baum-Connes map}

\label{s-j2}


%
The present section is subdivided into six subsections and presents the program of Section~\ref{s-appl-BC}
for the case $j=2$\,.


\subsection{Notation and Zimmermann's result}

\label{ss3.1}


\mbox{}\\
\vspace*{-.8em}

Let $\Sigma_{g}$ be a closed oriented Riemann surface of genus $g\geq 1$\,,
and let $\Gamma_{\!g}=\pi_{1}(\Sigma_{g})$ be its fundamental group;
$\Gamma_{\!g}$ admits the well-known presentation with $2g$ generators
and one relation
$$
\Gamma_{\!g}=\left<a_{1},b_{1},\dots,a_{g},b_{g}\,\left|\;\prod_{i=1}^{g}[a_{i},b_{i}]=1\right.\right>\,.
$$
The free group $\bbF_{g}$ of rank $g$ is isomorphic to the quotient of $\Gamma_{\!g}$ by
the normal subgroup generated by the $a_{i}b_{i}^{-1}$'s $(1\leq i\leq
g)$\,. It follows that every finitely generated group $\Gamma$ is a
quotient of some $\Gamma_{\!g}$ with $g$ big enough.

\smallskip

This remark was exploited by Zimmermann in~\cite{Zim} to give, for $\Gamma$ finitely
generated, a description of $H_{2}(\Gamma;\bbZ)$ in terms of pointed continuous
maps $\Sigma_{g}\longrightarrow B\Gamma$ inducing epimorphisms on fundamental groups.
We would like to avoid this assumption on $\Gamma$\,. It turns out that all the
results and their proofs in Zimmermann's article~\cite{Zim} are valid if one suppresses
the surjectivity assumption everywhere. Let us now explain the statements one obtains
this way. Denote by $S(\Sigma_{g},B\Gamma)$ the set of pointed continuous maps from
$\Sigma_{g}$ to $B\Gamma$ (\emph{not necessarily} inducing epimorphisms on fundamental
groups). Two maps $f_{1},f_{2}\in S(\Sigma_{g},B\Gamma)$ are called \emph{equivalent}
if there exists some orientation-preserving pointed homeomorphism $h$ of $\Sigma_{g}$
such that $f_{2}$ is homotopic to $f_{1}\circ h$\,.

\smallskip

Two maps $f_{1}\in S(\Sigma_{g_{1}},B\Gamma)$ and $f_{2}\in S(\Sigma_
{g_{2}},B\Gamma)$ are \emph{stably equivalent} if there exists closed oriented
Riemann surfaces $\Sigma^{\prime}$ and $\Sigma^{\prime\prime}$ such that
$f_{1}$ and $f_{2}$ become equivalent after being extended homotopically
trivially to the connected sums $\Sigma_{g_{1}}\#\Sigma^{\prime}$ and
$\Sigma_{g_{2}}\#\Sigma^{\prime\prime}$\,. More precisely, denoting
by $y_{0}$ the base-point of $B\Gamma$\,, we require the applications
$f_{1}\# y_{0}$ on $\Sigma_{g_{1}}\#\Sigma^{\prime}$ and $f_{2}\# y_
{0}$ on $\Sigma_{g_{2}}\#\Sigma^{\prime\prime}$ to be equivalent.

\smallskip

Denote by $\Omega(\Gamma)$ the set of stable equivalence classes in
$\coprod_{g\geq 1}S(\Sigma_{g},B\Gamma)$\,, and by $[f]$ the equivalence
class of $f\in S(\Sigma_{g},B\Gamma)$\,. Denote by $[\Sigma_{g}]
\in H_{2}(\Sigma_{g};\bbZ)$ the fundamental class of $\Sigma_{g}$\,.
The following result of Zimmermann~\cite{Zim} will be crucial.

\begin{Thm}[Zimmermann]
\label{thm2}
For an arbitrary discrete group $\Gamma$\,, the map
$$
Z_{\Gamma}\colon\Omega(\Gamma)\longrightarrow H_{2}(\Gamma;\bbZ)\,,\quad[f]
\longrightarrow f_{*}[\Sigma_{g}]\qquad\quad(\mbox{for}\;\,f\in S(\Gamma_{\!g},
B\Gamma))
$$
is a well-defined bijection (here, $f_{*}$ denotes $H_{2}(f;\bbZ)$).
\end{Thm}

Transferring the group structure of $H_{2}(\Gamma;\bbZ)$ to $\Omega
(\Gamma)$ via this bijection $Z_{\Gamma}$\,, we get a group structure
on $\Omega(\Gamma)$ such that
\begin{itemize}
    \item  [(1)] addition corresponds to connected sum (see Remark~\ref{rem-sum} below);

    \item  [(2)] the zero element is for example given by the class of the constant
    map in $S(\Sigma_{g},B\Gamma)$ (with $g\geq 1$ arbitrary);

    \item  [(3)] if $f\in S(\Sigma_{g},B\Gamma)$ is such that the homomorphism
    $\pi_{1}(f)\colon\Gamma_{\!g}\longrightarrow\Gamma$ factorizes through a free
    group, then $[f]$ is the zero element;

    \item  [(4)] for $f\in S(\Sigma_{g},B\Gamma)$\,, the opposite of $[f]$
    is given by $[f\circ h_{-}]$\,, where $h_{-}$ is an orientation-reversing
    pointed homeomorphism of $\Sigma_{g}$\,.
\end{itemize}

From now on, we shall implicitly identify $H_{2}(\Gamma;\bbZ)$ with
$\Omega(\Gamma)$ by the map $Z_{\Gamma}$\,, which has become a group
isomorphism.

\begin{Rem}
\label{rem-sum}
Let $\Gamma$ be a group. Consider $f_{1}\in S(\Sigma_{g_{1}},B\Gamma)$ and $f_{2}\in
S(\Sigma_{g_{2}},B\Gamma)$\,, and their classes in $\Omega(\Gamma)$\,. Up to stable equivalence
and up to homotopy, we can suppose that $f_{1}$ and $f_{2}$ are constant on a handle of
$\Sigma_{g_{1}}$ and $\Sigma_{g_{2}}$ respectively (and therefore also on a small disk).
Then, according to Example~\ref{ex-hmlgy-mfld} and to Proposition~\ref{prop-sum-var},
the class of $f_{1}+f_{2}$ in $\Omega(B\Gamma)$ is represented by the following two maps\,:
$$
f_{1}\# f_{2}\in S(\Sigma_{g_{1}+g_{2}},B\Gamma)\qquad\mbox{and}\qquad f_{1}\natural f_{2}
\in S(\Sigma_{g_{1}+g_{2}-1},B\Gamma)\,,
$$
where we identify $\Sigma_{g_{1}}\#\Sigma_{g_{2}}$ with $\Sigma_{g_{1}+g_{2}}$\,, and $\Sigma_{g_{1}}
\natural\Sigma_{g_{2}}$ with $\Sigma_{g_{1}+g_{2}-1}$\,, as usual.
\end{Rem}

Note that in the whole subsection, we can replace the particular connected CW-complex $B\Gamma$
by an arbitrary connected CW-complex $X$\,.


\subsection{Topological definition of $\betat{2}$}

\label{ss3.2}


\mbox{}\\
\vspace*{-.8em}

Keeping notation as in Lemma~\ref{lem3}, we let
$$
[\Sigma_{g}]_{K}:=(ch_{\even}^{\bbZ})^{-1}\big([\Sigma_{g}]\big)\;\in\;K_{0}(\Sigma_{g})
$$
be the (unique) $K$-homology class with integral Chern character given by the fundamental
class (an orientation, and even an auxiliary K\"{a}hler structure, is fixed on~$\Sigma_{g}$).
For $f\in S(\Sigma_{g},B\Gamma)$\,, we denote by $f_{*}\colon K_{0}(\Sigma_{g})\longrightarrow
K_{0}(B\Gamma)$ the induced map in $K$-homology. Now, we set
$$
\betat{2}\colon H_{2}(\Gamma;\bbZ)\longrightarrow K_{0}(B\Gamma)\,,\quad
[f]\longmapsto f_{*}[\Sigma_{g}]_{K}\qquad(\mbox{for}\;\,f\in S
(\Gamma_{\!g},B\Gamma))\,.
$$
It is not at all obvious that $\betat{2}$ is well-defined, and that it is a
group homomorphism; this will be stated as Theorem~\ref{Main2} below.


\subsection{Analytical definition of $\betat{2}$}

\label{ss3.3}


\mbox{}\\
\vspace*{-.8em}

Bearing in mind the analytical definition of $K$-homology,
it is interesting to express $[\Sigma_{g}]_{K}\in K_{0}(\Sigma_{g})$ in
this setting. This is precisely the subject of the next lemma, which follows
directly from Lemma~\ref{lem-Riemann-Roch}.

\begin{Lem}
\label{lem2}
One has $[\Sigma_{g}]_{K}=[\bar{\partial}_{g}]+(g-1)\cdot\iotaX{\Sigma_{g}}_{*}[1]$\,, where
$\iotaX{\Sigma_{g}}\colon pt\longrightarrow\Sigma_{g}$ is the inclusion of the base-point, and
$[1]$ is the canonical generator of $K_{0}(pt)\cong\bbZ$\,.
\qed
\end{Lem}

Let $X$ be a compact subspace of $B\Gamma$ such that $f(\Sigma_{g})\subseteq X$\,,
as for example $f(\Sigma_{g})$ itself. Now, the $K$-homology generator $\iotaX
{\Sigma_{g}}_{*}[1]$ is given by the Fredholm module
$$
\iotaX{\Sigma_{g}}_{*}[1]=[\bbC,\ev_{\Sigma_{g}},0]\in\KK_{0}(C(\Sigma_{g}),\bbC)\,,
$$
where $\ev_{\Sigma_{g}}\colon C(\Sigma_{g})\longrightarrow\bbC$ is evaluation at the base-point
of the surface $\Sigma_{g}$\,. Fix a map $f\in S(\Gamma_{\!g},B\Gamma)$\,. In the analytic framework,
$\betat{2}[f]$ is the image of the element
$$
f_{*}\big[L^{2}(\Lambda^{0,*}T^{*}\Sigma_{g}),\mathcal{M},\bar{\partial}_{\Sigma_{g}}
\big]+(g-1)\cdot[\bbC,\ev_{\!X},0]
$$
(where $f$ is viewed as a map from $\Sigma_{g}$ to $X$) under the homomorphism
$$
\KK_{0}(C(X),\bbC)\longrightarrow\RKK_{0}(B\Gamma,\bbC)
$$
induced by the inclusion of $X$ in $B\Gamma$\,, where $\ev_{\!X}$ is evaluation at the base-point of~$X$\,.
Suppose $f$ is Lipschitz; up to homotopy, one can always assume this is the case, with $f(\Sigma_{g})$ as
suitable $X$\,. Then, letting $f^{*}\colon C(X)\longrightarrow C(\Sigma_{g})$ take a function $\vartheta$
to $\vartheta\circ f$\,, we see that the $*$-closed subalgebra $\Lip(X)$ of $C(X)$ is dense and that $f^{*}\Lip(X)$
verifies $f^{*}\Lip(X)\subseteq\Lip(\Sigma_{g})$ and consists therefore of functions that are differentiable
almost everywhere by Rademacher's Theorem again (see \cite[Thm.~11\,A]{Whi}); as a consequence,
$$
f_{*}\big[L^{2}(\Lambda^{0,*}T^{*}\Sigma_{g}),\mathcal{M},\bar{\partial}_{\Sigma_{g}}
\big]=\big[L^{2}(\Lambda^{0,*}T^{*}\Sigma_{g}),\mathcal{M}\circ f^{*},\bar{\partial}_
{\Sigma_{g}}\big]\,.
$$


\subsection{Properties of $\betat{2}$}

\label{ss3.4}


\mbox{}\\
\vspace*{-.8em}

\begin{Thm}
\label{Main2}
The following map is a well-defined group homomorphism\,:
$$
\betat{2}\colon H_{2}(\Gamma;\bbZ)\longrightarrow K_{0}(B\Gamma)\,,\quad[f]\longmapsto
f_{*}[\bar{\partial}_{g}]+(g-1)\cdot\iotaX{BG}_{*}[1]\,,
$$
for $f\in S(\Sigma_{g},B\Gamma)$\,, where $\iotaX{B\Gamma}$ stands for the inclusion
of the base-point of $B\Gamma$\,.
\end{Thm}

We postpone the proof to Subsection~\ref{ss3.5} below, and derive, here, some of its consequences.
We also point out that~\cite{MM1} contains a purely homotopical proof of the theorem.

\begin{Prop}
\label{prop4}
Let $ch_{\ev}\colon K_{0}(B\Gamma)\longrightarrow\bigoplus_{n=0}^{\infty}
H_{2n}(\Gamma;\bbQ)$ be the even Chern character. Then, one has $(ch_{\ev}\otimes
\id_{\bbQ})\circ(\betat{2}\otimes\id_{\bbQ})=\id_{H_{2}(\Gamma;\bbQ)}$\,.
\end{Prop}

\Prf
Let $[f]\in H_{2}(\Gamma;\bbZ)$ be represented by $f\in S(\Sigma_{g},B\Gamma)$\,.
By naturality of the Chern character, we have a commutative diagram
\begin{displaymath}
 \begin{diagram}[height=2em,width=3.5em]
   K_{0}(\Sigma_{g})\otimes_{\bbZ}\bbQ & \rTo^{f_{*}} & K_{0}(B\Gamma)
   \otimes_{\bbZ}\bbQ \\
   \dTo^{ch_{\ev}\otimes\id_{\bbQ}}_{\cong} & & \dTo_{ch_{\ev}\otimes\id_{\bbQ}}^{\cong} \\
   H_{\even}(\Sigma_{g};\bbQ) & \rTo^{f_{*}} & H_{\even}(B\Gamma;\bbQ) \\
 \end{diagram}
\end{displaymath}
Then, dropping ``\,$\otimes\id_{\bbQ}$\,'' and ``\,$\otimes 1$\,'' from the notation,
one computes
$$
ch_{\ev}\,\betat{2}[f]=ch_{\ev}\,f_{*}[\Sigma_{g}]_{K}=f_{*}\,ch_{\ev}[\Sigma_{g}]_{K}=f_{*}
[\Sigma_{g}]=[f]\,,
$$
where the last equality follows from the identification given by Theorem~\ref{thm2}.
\qed


\subsection{Proof of Theorem \ref{Main2}}

\label{ss3.5}


\mbox{}\\
\vspace*{-.8em}

We first show that $\betat{2}$ is well-defined. We then prove it is a homomorphism.

\smallskip

We start with the topological setting. Fix $f_{1}\in S(\Sigma_{g_{1}},B\Gamma)$ and $f_{2}\in
S(\Sigma_{g_{2}},B\Gamma)$\,. We first show that if $g_{1}>g_{2}$ and if $f_{1}=f_{2}\# y_{0}$\,,
then $(f_{1})_{*}[\Sigma_{g_{1}}]=(f_{2})_{*}[\Sigma_{g_{2}}]$ in the group $H_{2}(B\Gamma;\bbZ)$\,.
To do this, we embed $\Sigma_{g_{1}}$ and $\Sigma_{g_{2}}$ in $\bbR^{3}$ in such a way that
$\Sigma_{g_{1}}$ is contained in a tubular neighbourhood $V$ of $\Sigma_{g_{2}}$ (see Figure~3).

\hspace*{6.5em}\epsfig{file=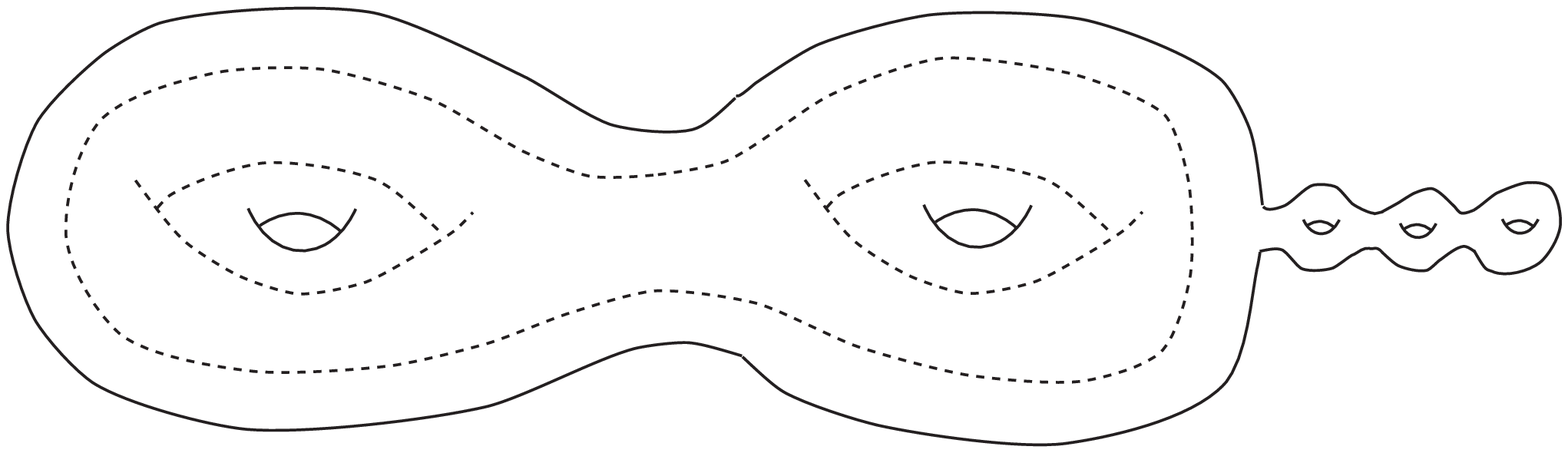,width=21em}

\vspace*{-.8em}

\centerline{Figure 3}

\vspace*{1em}

Identifying $V$ with the total space of the normal bundle of $\Sigma_
{g_{2}}$ yields a projection map $q\colon V\longrightarrow\Sigma_{g_{2}}$\,.
Clearly, the restriction $q|_{\Sigma_{g_{1}}}\colon\Sigma_{g_{1}}\longrightarrow
\Sigma_{g_{2}}$ is a smooth, proper and orientation preserving map; considering
the ``first handle'' (on the left in Figure~3) of $\Sigma_{g_{1}}$ and of $\Sigma_{g_{2}}$
(where $q|_{\Sigma_{g_{1}}}$ is one-to-one and regular), we see that it is of degree one,
so that $(q|_{\Sigma_{g_{1}}})_{*}[\Sigma_{g_{1}}]=[\Sigma_{g_{2}}]$\,. By naturality
and injectivity of the integral Chern character on the $K$-homology of closed oriented Riemann
surfaces (see Lemma~\ref{lem3}), we deduce that $(q|_{\Sigma_{g_{1}}})_{*}[\Sigma_
{g_{1}}]_{K}=[\Sigma_{g_{2}}]_{K}$ in $K_{0}(\Sigma_{g_{2}})$\,. On the other hand,
it is clear that the map $f_{1}=f_{2}\# y_{0}$ is homotopic to $f_{2}\circ q|_
{\Sigma_{g_{1}}}$\,, hence
$$
(f_{1})_{*} [\Sigma_{g_{1}}]_{K}=(f_{2}\circ q|_{\Sigma_{g_{1}}})_{*}[\Sigma_{g_{1}}]_{K}=
(f_{2})_{*}\circ(q|_{\Sigma_{g_{1}}})_{*}[\Sigma_{g_{1}}]_{K}=(f_{2})_{*}[\Sigma_{g_{2}}]_{K}\,.
$$
It remains to check that, if two maps $f_{1},f_{2}\in S(\Sigma_{g},B\Gamma)$ are equivalent,
then $(f_{1})_{*}[\Sigma_{g}]_{K}=(f_{2})_{*}[\Sigma_{g}]_{K}$ in $K_{0}(B\Gamma)$\,. This
follows from the fact that orientation-preserving homeomorphisms of $\Sigma_{g}$ induce the
identity on $K_{0}(\Sigma_{g})$ (again, this can be checked using the integral Chern character
and Lemma~\ref{lem3}). This shows that $\betat{2}$ is a well-defined map.

\smallskip

Now, we prove that $\betat{2}$ is a group homomorphism still in the topological setting.
We fix $f_{1}\in S(\Sigma_{g_{1}},B\Gamma)$ and $f_{2}\in S(\Sigma_{g_{2}},B\Gamma)$\,.
Using the first description of the sum in Remark~\ref{rem-sum}, we must show that
$$
(f_{1}\# f_{2})_{*}[\Sigma_{g_{1}+g_{2}}]_{K}=(f_{1})_{*}[\Sigma_{g_{1}}]_{K}+(f_
{2})_{*}[\Sigma_{g_{2}}]_{K}
$$
holds in $K_{0}(B\Gamma)$\,. We can now exploit Lemma~\ref{lem3} to reduce the proof to
showing the homological equality
$$
(f_{1}\# f_{2})_{*}[\Sigma_{g_{1}+g_{2}}]=(f_{1})_{*}[\Sigma_{g_{1}}]+(f_{2})_{*}[\Sigma_{g_{2}}]
$$
in $H_{2}(B\Gamma;\bbZ)$\,. This is a special case of Example~\ref{ex-hmlgy-mfld} (which was based
on the general principle~\ref{lem-gen-princ}). This completes the proof in the topological setting.
\qed

\smallskip

We move now to the analytical framework and present the corresponding proof of Theorem~\ref{Main2}.
We first have to show that the map
$$
\betat{2}\colon H_{2}(\Gamma;\bbZ)\longrightarrow K_{0}(B\Gamma)\,,\quad[f]\longmapsto f_{*}[\bar
{\partial}_{g}]+(g-1)\cdot\iotaX{BG}_{*}[1]\,,
$$
for $f\in S(\Sigma_{g},B\Gamma)$\,, is well-defined. The proof is subdivided into six steps.

\smallskip

(1) If $q_{0}\colon\Sigma_{g}\longrightarrow pt$ denotes the constant map, then $(q_{0})_{*}[\bar
{\partial}_{g}]=(1-g)\!\cdot\![1]$ holds in $K_{0}(pt)\cong\bbZ$\,. Indeed, the operator $\bar{\partial}_{g}$
has $1-g$ as index, see~\cite[p.~27]{Sha}.

\smallskip

(2) The group $K_{0}(\Sigma_{g_{1}}\!\vee\Sigma_{g_{2}})$ is isomorphic to $\bbZ^{3}$ with the elements
$\iotaX{\Sigma_{g_{1}}\!\vee\Sigma_{g_{2}}}_{*}[1]$\,, $[\bar{\partial}_{g_{1}}]$ and $[\bar{\partial}_{g_{2}}]$
as generators (using the obvious identifications), where $\iotaX{\Sigma_{g_{1}}\!\vee\Sigma_{g_{2}}}$
stands for the inclusion of the base-point, see Lemmas \ref{lem3} and \ref{lem-Riemann-Roch}.

\smallskip

(3) Let $x_{0}$ be the base-point of $\Sigma_{g_{1}}$ and consider the ``crunching'' map
$$
q:=\id_{\Sigma_{g_{1}}}\!\!\vee\;x_{0}\colon\Sigma_{g_{1}}\!\vee\Sigma_{g_{2}}\longrightarrow
\Sigma_{g_{1}}\,,\quad x\longmapsto
\left\{
\arraycolsep1pt
\renewcommand{\arraystretch}{1.2}
\begin{array}{ll}
  x\,,\; & \mbox{if}\;\,x\in\Sigma_{g_{1}} \\
  x_{0}\,,\; & \mbox{if}\;\,x\in\Sigma_{g_{2}}\,. \\
\end{array}\right.
$$
Then, under the identifications of (2), $q_{*}[\bar{\partial}_{g_{1}}]=[\bar{\partial}_{g_{1}}]$ and
$q_{*}[\bar{\partial}_{g_{2}}]=(1-g_{2})\cdot\iotaX{\Sigma_{g_{1}}}_{*}[1]$ hold in the group $K_{0}
(\Sigma_{g_{1}})$\,, as follows from (1) for the latter equality.

\smallskip

(4) Let $p\colon\Sigma_{g_{1}}\#\Sigma_{g_{2}}\longrightarrow\Sigma_{g_{1}}\!\vee\Sigma_{g_{2}}$ be the
pinching map that ``contracts'' the identification circle in $\Sigma_{g_{1}}\#\Sigma_{g_{2}}$ to the
base-point of $\Sigma_{g_{1}}\!\vee\Sigma_{g_{2}}$\,. Then, under the identifications of (2), the following
equality holds\,:
$$
p_{*}[\bar{\partial}_{g_{1}+g_{2}}]=[\bar{\partial}_{g_{1}}]+[\bar{\partial}_{g_{2}}]-\iotaX
{\Sigma_{g_{1}}\!\vee\Sigma_{g_{2}}}_{*}[1]\;\in\;K_{0}(\Sigma_{g_{1}}\!\vee\Sigma_{g_{2}})\,.
$$
This equality is the ``tricky'' part of the present proof (and it is precisely here that the proof
becomes of analytical nature properly speaking -- of course, this can also be directly established in
the topological framework, using the integral Chern character of Lemma~\ref{lem3} and Lemma~\ref{lem-Riemann-Roch},
thus yielding a second topological proof of the well-definiteness). Let $K$ be a small closed neighbourhood
of the base-point $x_{0}$ in $\Sigma_{g_{1}}\!\vee\Sigma_{g_{2}}$ ($K$ is contractible), and let
$K':=p^{-1}(K)$ be the corresponding closed tubular neighbourhood of the identification
circle $p^{-1}(x_{0})$ in $\Sigma_{g_{1}}\#\Sigma_{g_{2}}$ ($K$ is homotopy equivalent to $S^{1}$).
Let $U$ and $U'$ be the complements of $K$ and $K'$ in $\Sigma_{g_{1}}\!\vee\Sigma_{g_{2}}$
and $\Sigma_{g_{1}}\#\Sigma_{g_{2}}$\,, respectively. We can assume that the map $p|_{U'}\colon
U'\longrightarrow U$ is an isometry. The short exact sequences of $C^{*}$-algebras
$$
0\longrightarrow C_{0}(U')\stackrel{i}{\longrightarrow}C(\Sigma_{g_{1}}\#\Sigma_{g_{2}})\stackrel{r}
{\longrightarrow}C(K')\longrightarrow 0
$$
and
$$
0\longrightarrow C_{0}(U)\stackrel{j}{\longrightarrow}C(\Sigma_{g_{1}}\!\vee\Sigma_{g_{2}})\stackrel{s}
{\longrightarrow}C(K)\longrightarrow 0
$$
give rise to the following commutative diagram with exact rows\,:
\begin{displaymath}
 \begin{diagram}[height=2em,width=3em]
    0 & \rTo & K_{0}(K') & \rTo^{\!\!\!\!i_{*}} & K_{0}(\Sigma_{g_{1}}\#\Sigma_{g_{2}}) & \rTo^{\;\;\,r_{*}} &
    K_{0}(U') & & \\
     & & \dTo_{(p|_{K'})_{*}}^{\cong} & & \dTo_{p_{*}} & & \dTo_{(p|_{U'})_{*}}^{\cong} & & \\
    0 & \rTo & K_{0}(K) & \rTo^{\!\!\!\!j_{*}} & K_{0}(\Sigma_{g_{1}}\!\vee\Sigma_{g_{2}}) & \rTo^{\;\;\,s_{*}}
    & K_{0}(U) & \rTo & 0 \\
 \end{diagram}
\end{displaymath}
Note that $K_{0}(K)$ and $K_{0}(K')$ are both isomorphic to $K_{0}(pt)\cong\bbZ$\,, and it is for
this reason that $(p|_{K'})_{*}$ is an isomorphism and that both $i_{*}$ and $j_{*}$ are injective.
Now, since, for each $g$\,, $\bar{\partial}_{g}$ is a symmetric elliptic operator on a Riemannian
manifold, Proposition~\cite[Prop.~10.8.8]{HigRoe} (which is of purely analytical nature) can be
applied, and we have
$$
\arraycolsep1pt
\renewcommand{\arraystretch}{1.2}
\begin{array}{rclcl}
  s_{*}\circ p_{*}[\bar{\partial}_{g_{1}+g_{2}}] & = & p_{*}\circ i_{*}[\bar{\partial}_{g_{1}+g_{2}}] &
  \qquad & \mbox{(by commutativity of the diagram)} \\
  & = & p_{*}[\bar{\partial}_{g_{1}+g_{2}}|_{U'}] & \qquad & \mbox{(by \cite[Prop.~10.8.8]{HigRoe})} \\
  & = & p_{*}[\bar{\partial}_{g_{1}}|_{U'}]+p_{*}[\bar{\partial}_{g_{2}}|_{U'}] & \qquad & \mbox{(by
  the local description of $\bar{\partial}_{g}$)} \\
  & = & [\bar{\partial}_{g_{1}}|_{U}]+[\bar{\partial}_{g_{2}}|_{U}] & \qquad & \mbox{(since $p|_{U'}$
  is an isometry)} \\
  & = & s_{*}[\bar{\partial}_{g_{1}}]+s_{*}[\bar{\partial}_{g_{2}}]& \qquad & \mbox{(by \cite[Prop.~10.8.8]{HigRoe})}\,. \\
\end{array}
$$
Therefore, it follows that
$$
p_{*}[\bar{\partial}_{g_{1}+g_{2}}]-[\bar{\partial}_{g_{1}}]-[\bar{\partial}_{g_{2}}]\in\Ker(s_{*})=
\Ima(j_{*})=\bbZ\cdot\iotaX{\Sigma_{g_{1}}\!\vee\Sigma_{g_{2}}}[1]\,.
$$
The determination of the corresponding integer $\lambda$ (which we have to show is $-1$) amounts to the
determination of the indices, namely
$$
\lambda=\Index(\bar{\partial}_{g_{1}+g_{2}})-\Index(\bar{\partial}_{g_{1}})-\Index(\bar{\partial}_{g_{2}})
=(1-g_{1}-g_{2})-(1-g_{1})-(1-g_{2})=-1\,,
$$
by~\cite[p.~27]{Sha}, and we are done.

\smallskip

(5) By (3) and (4), using the same notation, one has the following equality\,:
$$
(q\circ p)_{*}[\bar{\partial}_{g_{1}+g_{2}}]=[\bar{\partial}_{g_{1}}]-g_{2}\!\cdot\!\iotaX{\Sigma_{g_{1}}}_{*}[1]
\;\in\;K_{0}(\Sigma_{g_{1}})\,.
$$

\smallskip

(6) We now really establish the well-definiteness of $\betat{2}$\,. To verify it, we must show that if
two maps
$$
f_{1}\colon\Sigma_{g_{1}}\longrightarrow B\Gamma\qquad\mbox{and}\qquad f_{2}\colon\Sigma_{g_{1}}\#\Sigma_{g_{2}}
\longrightarrow B\Gamma
$$
are related by the equality $f_{2}=f_{1}\# y_{0}$\,, with $y_{0}$ standing for the base-point of $B\Gamma$\,,
then
$$
(f_{2})_{*}[\bar{\partial}_{g_{1}+g_{2}}]+(g_{1}+g_{2}-1)\cdot\iotaX{B\Gamma}_{*}[1]=(f_{1})_{*}[\bar{\partial}_
{g_{1}}]+(g_{1}-1)\cdot\iotaX{B\Gamma}_{*}[1]
$$
holds in $K_{0}(B\Gamma)$\,. The key observation is that $f_{2}=f_{1}\circ q\circ p$\,, so, by virtue of
(5),
$$
(f_{2})_{*}[\bar{\partial}_{g_{1}+g_{2}}]=(f_{1})_{*}[\bar{\partial}_{g_{1}}]-g_{2}\cdot\iotaX{B\Gamma}_{*}[1]
$$
and we can conclude.

\smallskip

Finally, we show that $\betat{2}$ is a group homomorphism in the analytical setting.
Again, we fix $f_{1}\in S(\Sigma_{g_{1}},B\Gamma)$ and $f_{2}\in S(\Sigma_{g_{2}},B\Gamma)$\,.
Consider the compact subspace $X:=f_{1}(\Sigma_{g_{1}})\cup f_{2}(\Sigma_{g_{2}})$ of $B\Gamma$\,.
Using the second description of the sum in Remark~\ref{rem-sum}, according to Lemma~\ref{lem2},
we must show that
$$
(f_{1}\natural f_{2})_{*}\big((1-(g_{1}+g_{2}-1))\cdot\iotaX{\Sigma_{g_{1}}\!\natural\Sigma_{g_{2}}}_
{*}[1]\big)=(f_{1})_{*}\big((1-g_{1})\cdot\iotaX{\Sigma_{g_{1}}}[1]\big)+(f_{2})_{*}\big((1-g_{2})
\cdot\iotaX{\Sigma_{g_{2}}}_{*}[1]\big)
$$
and that
$$
(f_{1}\natural f_{2})_{*}[\bar{\partial}_{g_{1}+g_{2}-1}]=(f_{1})_{*}[\bar{\partial}_{g_{2}}]+(f_{2})_{*}[\bar{
\partial}_{g_{2}}]
$$
in $\KK_{0}(C(X),\bbC)$\,. For the first equality, it suffices to note that $f_{1}\natural f_{2}$\,,
$f_{1}$ and $f_{2}$ are pointed maps, so that this reduces to an equality of integers. The second is
the content of Theorem~\ref{thm-add-formula-2}, that we have proved both in the topological and in
the analytical settings.
\qed


\subsection{Definition of the map $\betaa{2}$ and connection with $\betat{2}$}

\label{ss3.7}


\mbox{}\\
\vspace*{-.8em}

We now construct the map $\betaa{2}\colon H_{2}(\Gamma;\bbZ)\longrightarrow
K_{0}(C^{*}_{r}\Gamma)$\,. Denote by $C^{*}\Gamma$ the full $C^{*}$-algebra of
the group $\Gamma$\,, and by $\lambda_{\Gamma}\colon C^{*}\Gamma\onto
C^{*}_{r}\Gamma$ the canonical epimorphism. It is well-known that the Novikov
assembly map factors through the $K$-theory of the full $C^{*}$-algebra (see
\cite{Jul1} or \cite[Section~2.3 in Part~2]{Val}), \ie for $i=0,1$\,, there
is a homomorphism
$$
\tilde{\nu}_{i}^{\Gamma}\colon K_{i}(B\Gamma)\longrightarrow K_{i}(C^{*}
\Gamma)\quad\mbox{such that}\quad\nu^{\Gamma}_{i}=(\lambda_{\Gamma})_{*}
\circ\tilde{\nu}^{\Gamma}_{i}\,.
$$
For a map $f\in S(\Sigma_{g},B\Gamma)$\,, we denote by the same symbol the associated group homomorphism
$\pi_{1}(f)\colon\Gamma_{\!g}\longrightarrow\Gamma$\,, and also the corresponding $*$-homomorphism
$C^{*}\big(\pi_{1}(f)\big)\colon C^{*}\Gamma_{\!g}\longrightarrow C^{*}\Gamma$ (the latter being
well-defined thanks to the universal property of the full $C^{*}$-algebra). We define
$$
\betaa{2}\colon H_{2}(\Gamma;\bbZ)\longrightarrow K_{0}(C^{*}_{r}\Gamma)\,,\quad
[f]\longmapsto(\lambda_{\Gamma}\circ f)_{*}\,\tilde{\nu}_{0}^{\Gamma_{\!g}}[\Sigma_
{g}]_{K}\quad(\mbox{for}\;\,f\in S(\Sigma_{g},B\Gamma))\,.
$$

\begin{Thm}
\label{thm3}
The map $\betaa{2}$ is a well-defined group homomorphism satisfying the
equality $\betaa{2}=\nu_{0}^{\Gamma}\circ\betat{2}$\,.
\end{Thm}

\Prf
For $f\in S(\Sigma_{g},B\Gamma)$\,, we have to show that $\nu^{\Gamma}_{0}
\betat{2}[f]=(\lambda_{\Gamma}\circ f)_{*}\,\tilde{\nu}_{0}^{\Gamma_{\!g}}
[\Sigma_{g}]_{K}$ in $K_{0}(C^{*}_{r}\Gamma)$\,. The result will follow since,
by Theorem~\ref{Main2}, the left-hand side only depends on the class
$[f]$ of $f$ in $H_{2}(\Gamma;\bbZ)$\,, and moreover $\betat{2}$ is a
group homomorphism. Now, the map $\tilde{\nu}^{\Gamma}_{i}$ is natural
with respect to \emph{arbitrary} group homomorphisms (and not just injective
ones, see~\cite[Thm.~1.1 in Part~2]{Val}), so that
$$
\arraycolsep1pt
\renewcommand{\arraystretch}{1.2}
\begin{array}{rcl}
  (\lambda_{\Gamma}\circ f)_{*}\,\tilde{\nu}_{0}^{\Gamma_{\!g}}[\Sigma_{g}]_{K} & = &
  (\lambda_{\Gamma})_{*}\,f_{*}\,\tilde{\nu}_{0}^{\Gamma_{\!g}}[\Sigma_{g}]_{K}=(\lambda_
  {\Gamma})_{*}\,\tilde{\nu}_{0}^{\Gamma}f_{*}[\Sigma_{g}]_{K} \\
  & = & \nu^{\Gamma}_{0}f_{*}[\Sigma_{g}]_{K}=\nu^{\Gamma}_{0}\betat{2}[f]\,.
\end{array}
$$
This completes the proof.
\qed

\smallskip

In the unbounded analytical description of $K\!K$-theory in the sense of~\cite{BJ},
the `universal' class $\tilde{\nu}_{0}^{\Gamma_{\!g}}[\Sigma_{g}]_{K}\in K_{0}(C^{*}\Gamma_{\!g})=
\KK_{0}(\bbC,C^{*}\Gamma_{\!g})$ is given by the unbounded Kasparov triple
$$
\tilde{\nu}_{0}^{\Gamma_{\!g}}[\Sigma_{g}]_{K}=[\mathcal{E}_{g},\bar{\partial}_{g}]
=[\mathcal{E}_{g},\pi_{\circ},\bar{\partial}_{g}]\;\in\;\KK_{0}(\bbC,C^{*}\Gamma_{\!g})\,,
$$
where $\mathcal{E}_{g}$ is defined as we next explain and $\pi_{\circ}\colon\bbC\longrightarrow\mathcal{L}_
{C^{*}\Gamma_{\!g}}(\mathcal{E}_{g})$ is the unit. Letting $\widetilde{\Sigma}_{g}$ be the universal
cover of $\Sigma_{g}$\,, $\mathcal{E}_{g}$ is the separation-completion of the algebra $\Gamma_{\!c}
(\Lambda^{0,*}T^{*}\widetilde{\Sigma}_{g})$ of compactly supported smooth sections of
the vector bundle $\Lambda^{0,*}T^{*}\widetilde{\Sigma}_{g}$ over $\widetilde{\Sigma}_{g}$ with
respect to the $C^{*}\Gamma_{\!g}$-valued scalar product determined by
$$
\left<\xi_{1}|\xi_{2}\right>(\sigma):=\left<\xi_{1}|\sigma\cdot\xi_{2}\right>_{L^{2}(\widetilde
{\Sigma}_{g},\Lambda^{0,*}T^{*}\widetilde{\Sigma}_{g})}\,,
$$
for $\xi_{1},\xi_{2}\in\Gamma_{\!c}(\Lambda^{0,*}T^{*}\widetilde{\Sigma}_{g})$
and $\sigma\in\Gamma_{\!g}$ (acting on $\Gamma_{\!c}(\Lambda^{0,*}T^{*}\widetilde
{\Sigma}_{g})$ in the usual way, via deck transformations), compare with D.~Kucerovsky's Appendix
to~\cite{Val}. It follows that for $f\in S(\Sigma_{g},B\Gamma)$\,, we have
$$
\betaa{2}[f]=[\mathcal{E}_{g}',\bar{\partial}_{g}]=[\mathcal{E}_{g}',\pi_{\circ}',\bar{\partial}_{g}]
\;\in\;\KK_{0}(\bbC,C^{*}_{r}\Gamma)\,,
$$
where $\pi_{\circ}'\colon\bbC\longrightarrow\mathcal{L}_{C^{*}_{r}\Gamma}(\mathcal{E}_{g}')$ is
the unit, and $\mathcal{E}_{g}'$ is the separation-completion of $\Gamma_{\!c}
(\Lambda^{0,*}T^{*}\widetilde{\Sigma}_{g})$ with respect to the $C^{*}_{r}\Gamma$-valued scalar
product determined by
$$
\left<\xi_{1}|\xi_{2}\right>(\gamma):=\sum_{\!\!\!\sigma\in\pi_{1}(f)^{-1}(\gamma)\!\!\!}\left<
\xi_{1}|\sigma\cdot\xi_{2}\right>_{L^{2}(\widetilde{\Sigma}_{g},\Lambda^{0,*}T^{*}\widetilde{\Sigma}_{g})}\,,
$$
for $\xi_{1},\xi_{2}\in\Gamma_{\!c}(\Lambda^{0,*}T^{*}\widetilde{\Sigma}_{g})$
and $\gamma\in\Gamma$\,, see \cite[Section~3 in Part~2]{Val}. This provides a purely analytical
description of $\betaa{2}$\,. See also \cite[Section 3]{MOO} for information on $\betaa{2}[f]$
in connection with group homology and algebraic $K$-theory, described therein via an element
$\nu_{2}[\Sigma_{g},f]$ lying in a suitable quotient of $K_{2}^{\rm alg}(\bbZ\Gamma)$\,.


\section{The case of $2$-dimensional groups}

\label{s-two-dim}


%
Recall that we call a group $\Gamma$ \emph{$2$-dimensional} if its classifying space has
the homotopy type of a CW-complex (not necessarily finite) of dimension $\leq 2$\,.

\smallskip

Examples of $2$-dimensional groups abound\,:
\begin{itemize}
    \item  [(1)] Surface groups\,: The Baum-Connes Conjecture was proved for those groups by Kasparov~\cite{Kas2}.

    \item  [(2)] Torsion-free one-relator groups\,: For this class, the Baum-Connes Conjecture was established
    in~\cite{BBV}.

    \item  [(3)] Knot groups\,: By~\cite{BBV}, they also satisfy the Baum-Connes Conjecture.

    \item  [(4)] Groups acting freely co-compactly on a $2$-dimensional
    Euclidean building\,: These groups have Kazhdan's property $(T)$
    (see~\cite{Zuk} for an elegant proof of this fact). For groups
    acting on $\tilde{A}_{2}$-buildings (in particular co-compact
    torsion-free lattices in $\mathop{\rm PGL}_3(F)$\,, with $F$ a local
    field), the Baum-Connes Conjecture is an outstanding result
    of Lafforgue~\cite{Laf}. For other cases (e.g.\ co-compact
    torsion-free lattices in the symplectic group $\mathop{\rm Sp}_4
    (F)$\,, $F$ a local field), the Baum-Connes Conjecture is still
    open. Let us mention however that, in these cases, it is known
    by work of Kasparov and Skandalis~\cite{KaS} that the Novikov
    assembly map $\nu^{\Gamma}_{*}$ is injective.

    \item  [(5)] It was shown by Champetier~\cite{Cha} that there is a
    certain genericity of $2$-dimensional groups among finitely
    presentable groups. Indeed, fix the finite generating set $X$ and
    the number $k$ of relations. Among groups $\Gamma=\left<X\,\left|
    \,r_{1},\ldots,r_{k}\right.\right>$ generated by $X$ and on $k$
    relations, the proportion of $2$-dimensional groups goes to $1$
    as\, $\max\left\{|r_{1}|,\ldots,|r_{k}|\right\}\rightarrow
    +\infty$\, (see \cite[pp.\ 199--200]{Cha}); moreover, for $k=2$\,,
    there is genericity in the stronger sense of Gromov, namely, the
    proportion of $2$-dimensional groups goes to $1$ even as\, $\min
    \left\{|r_{1}|,\,|r_{2}|\right\}\rightarrow+\infty$ (see \cite[Thm.~4.13]{Cha}).

    \item  [(6)] The following result is proved by Wise in~\cite{Wise}. Suppose
    given an arbitrary finitely presentable group $\Gamma$\,. Then, there exists a
    compact negatively curved $2$-dimensional simplicial complex $X$ and a finitely
    generated normal subgroup $N$ of $\pi_{1}(X)$ such that $\pi_{1}(X)/N\cong\Gamma$\,.
    Negative curvature implies that $X$ is acyclic and therefore a model for
    $B\pi_{1}(X)$\,; as a consequence, $\pi_{1}(X)$ is a $2$-dimensional group.
    In particular, any finitely presentable group is a quotient of some
    (finitely presentable) $2$-dimensional group.
\end{itemize}

What is special about $2$-dimensional groups in our context comes from
the canonical ``identification'' between the integral homology of the group
and the $K$-homology of its classifying space, see Lemma~\ref{lem3}.

\begin{Lem}
\label{lem4}
Let $\Gamma$ be a $2$-dimensional group. Then the maps
$$
\betatev\colon\bbZ\oplus H_{2}(\Gamma;\bbZ)\stackrel{\cong}{\longrightarrow}K_{0}(B\Gamma)\,,\quad
(m,[f])\longmapsto m\!\cdot\!\iotaX{B\Gamma}_{*}[1]+\betat{2}[f]
$$
and
$$
\betat{1}\colon H_{1}(\Gamma;\bbZ)\stackrel{\cong}{\longrightarrow}K_{1}(B\Gamma)\,,\quad\gamma^{\ab}
\longmapsto\gamma_{*}[S^{1}]_{K}={-\gamma_{*}[D]}\,,
$$
are isomorphisms, as indicated.
\end{Lem}

\Prf
Since $B\Gamma$ is at most $2$-dimensional, we have first that its integral homology is torsion-free
(so that it injects into its rational homology), and, second, by Lemma~\ref{lem3}, we have commutative
diagrams
\begin{displaymath}
 \begin{diagram}[height=2em,width=2.5em]
   & & K_{0}(B\Gamma) & & & \qquad\quad & & & K_{1}(B\Gamma) & &  \\
   & \ldTo^{ch_{\even}^{\bbZ}}_{\!^{\cong}} & & \rdTo^{ch_{\ev}} & & \qquad\quad & & \ldTo^{ch_{\odd}^{\bbZ}}_{\!^{\cong}}
   & & \rdTo^{ch_{\odd}} & \\
   H_{\even}(\Gamma;\bbZ) & & \rInto & & H_{\even}(\Gamma;\bbQ) & \qquad\quad & H_{1}(\Gamma;\bbZ) & &
   \rInto & & H_{1}(\Gamma;\bbQ) \\
 \end{diagram}
\end{displaymath}

\noindent
By Propositions \ref{prop2} and \ref{prop4}, the maps $\betat{1}$ and $\betat{2}$ are, rationally, right-inverses
of the Chern character $ch_{*}$ in the corresponding degrees. A corresponding result holds for $\betat{0}$\,,
see Section~\ref{s-appl-BC}. By diagram chase, it follows that $\betatev$ and $\betat{1}$ are isomorphisms,
as was to be shown.
\qed

\smallskip

From this lemma and Theorems \ref{thm1} and \ref{thm3}, we immediately
get the following reformulation of the Baum-Connes Conjecture for
$2$-dimensional groups.

\begin{Prop}
\label{prop5}
For a 2-dimensional group $\Gamma$\,, the Baum-Connes Conjecture is
equivalent to the following statement\,: the maps
$$
\betaaev\colon\bbZ\oplus H_{2}(\Gamma;\bbZ)\longrightarrow K_{0}(C^{*}_{r}
\Gamma)\,,\;(m,[f])\longmapsto m\!\cdot\![1]+\betaa{2}[f]
$$
and
$$
\betaa{1}\colon H_{1}(\Gamma;\bbZ)=\Gamma^{\ab}\longrightarrow K_{1}(C^{*}_{r}\Gamma)\,,
\quad\gamma^{\ab}\longmapsto[\gamma]=\big[\diag(\gamma,1,1,\ldots)\big]\,,
$$
are isomorphisms.
\qed
\end{Prop}

We single out one consequence of \emph{surjectivity} of the
Baum-Connes assembly map, consistent with the philosophy that
surjectivity implies analytical results.

\begin{Cor}
\label{cor1}
Let $\Gamma$ be a $2$-dimensional group. Suppose that the assembly
map $\nu_{1}^{\Gamma}\colon K_{1}(B\Gamma)\longrightarrow K_{1}(C^{*}_
{r}\Gamma)$ is onto. Then every element of $\mathop{\rm GL}_
{\infty}(C^{*}_{r}\Gamma)$ lies in the same path-component
as a diagonal matrix $\diag(\gamma,1,1,\ldots)$\,, for
some $\gamma\in\Gamma$\,.
\end{Cor}

\Prf
Since $K_{1}(C^{*}_{r}\Gamma)$ is by definition the group of path-components
of $\mathop{\rm GL}_{\infty}(C^{*}_{r}\Gamma)$\,, the result follows from the
previous one together with the very definition of $\betaa{1}$\,.
\qed

\begin{Rems}
\mbox{}\\
\vspace*{-1.2em}
\begin{itemize}
    \item  [(1)] Suppose that $\Gamma$ is a $2$-dimensional group. One may rephrase
    the previous corollary by saying that the quotient group $K_{1}(C^{*}_{r}\Gamma)
    \big/\big<[\gamma]\,\big|\,\gamma\in\Gamma\big>$ is zero if and only if $\nu_
    {1}^{\Gamma}$ is surjective. Now, observe that for an arbitrary discrete group
    $G$\,, the class $[-1]\in K_{1}(C^{*}_{r}G)$ of the diagonal matrix $\diag(-1,
    1,1,\ldots)$ is zero; indeed, this class lies in the image of the canonical
    homomorphism $K_{1}(\bbC)\longrightarrow K_{1}(C^{*}_{r}G)$ and $K_{1}(\bbC)=0$\,.
    In particular, $\nu_{1}^{\Gamma}$ is surjective if and only if the group
    $$
    \itemspace W\!h^{\topo}(\Gamma):=K_{1}(C^{*}_{r}\Gamma)\Big/\big<[\pm\gamma]\,
    \big|\,\gamma\in\Gamma\big>=K_{1}(C^{*}_{r}\Gamma)\Big/\big<[\gamma]\,\big|\,
    \gamma\in\Gamma\big>
    $$
    vanishes. The definition of this quotient is somewhat reminiscent of the
    definition of the \emph{Whitehead group} in algebraic $K$-theory (hence our
    notation)\,:
    $$
    \itemspace W\!h(\Gamma):=K_{1}^{\alg}(\bbZ\Gamma)\Big/\big<[\pm\gamma]\,\big|\,
    \gamma\in\Gamma\big>\,,
    $$
    see e.g.~\cite{Ros}. It follows from~\cite[Thm.~1.1]{MOO} that the map $\betaa{1}$
    factorizes through the algebraic $K$-group $K_{1}^{\alg}(\bbZ\Gamma)$ (for an
    arbitrary group $\Gamma$). Therefore, we can also deduce from this all that for
    our $2$-dimensional group $\Gamma$\,, the following three statements are implied
    by the surjectivity of $\nu_{1}^{\Gamma}$\,:
    \begin{itemize}
        \item [(a)] the canonical map $K_{1}^{\alg}(\bbZ\Gamma)\longrightarrow
        K_{1}(C^{*}_{r}\Gamma)$ is surjective;
        \item [(b)] the canonical map $W\!h(\Gamma)\longrightarrow
        W\!h^{\topo}(\Gamma)$ is surjective;
        \item [(c)] $W\!h^{\topo}(\Gamma)=0$\,.
    \end{itemize}
    It would be of great interest to study these three properties independently
    of the Baum-Connes Conjecture, and for a larger class of groups.

    \item  [(2)]  Let $\Gamma$ be a discrete group. If $M$ is a closed oriented manifold
    equipped with a continuous map $M\longrightarrow B\Gamma$\,, then all higher signatures of $M$
    coming via $f$ from classes lying in the subring of $H^{*}(\Gamma;\bbQ)$ generated by $H^{j}
    (\Gamma;\bbQ)$ with $j\leq 2$ are oriented homotopy invariants of $M$\,: this is an
    unpublished result of Connes, Gromov and Moscovici (see however~\cite{Gro}); a complete
    proof is now available, see \cite[Cor.~0.3]{VM}. As a corollary, the usual Novikov Conjecture
    in topology holds for a $2$-dimensional group. It is not clear to us whether $\nu^{\Gamma}_{0}$
    is rationally injective for $\Gamma$ a $2$-dimensional group. With no doubt, this would
    constitute a useful result.
\end{itemize}
\end{Rems}

\medskip


{\bf Acknowledgements\,:} Thanks are due to P.~Baum, N.~Higson, P.~Julg and W.~L\"uck for
helpful conversations. We feel considerably indebted to G.~Skandalis for his precious and
generous help.



\end{document}